\documentclass[reqno, 12pt]{amsart}

\usepackage{array}
\usepackage{amsmath}
\usepackage{amsfonts}
\usepackage{amssymb}
\usepackage{stmaryrd}
\usepackage{enumerate}
\usepackage{amsthm}
\usepackage{amsmath, amscd}
\usepackage{xy}
\xyoption{all}
\usepackage{supertabular}
\usepackage{marvosym}
\usepackage{wasysym}
\usepackage{hyperref}
\usepackage{tikz}
\usetikzlibrary{matrix,arrows,backgrounds}%,decorations, pathmorphing, positioning, fit}
	
\newtheorem{theorem}{Theorem}[section]

\newtheorem{lemma}[theorem]{Lemma}
\newtheorem{proposition}[theorem]{Proposition}
\newtheorem{corollary}[theorem]{Corollary}

\newtheorem{conjecture}[theorem]{Conjecture}

\theoremstyle{definition}
\newtheorem{definition}[theorem]{Definition}
\newtheorem{definition-lemma}[theorem]{Definition/Lemma}
\newtheorem{remark}[theorem]{Remark}

\newtheorem{example}[theorem]{Example}

\newcommand{\tw}[1]{\lbrace #1 \rbrace}

\newcommand{\op}[1]{\operatorname{#1}}
\newcommand{\ops}[1]{\operatorname{#1 }}
\newcommand{\leftexp}[2]{{\vphantom{#2}}^{#1}{#2}}

\newcommand{\dbcoh}[1]{\operatorname{D}^{\operatorname{b}}(\operatorname{coh }#1)}

\newcommand{\dbmod}[1]{\operatorname{D}^{\operatorname{b}}(\operatorname{mod }#1)}

\newcommand{\dbgr}[1]{\operatorname{D}^{\operatorname{b}}(\operatorname{gr }#1)}
\newcommand{\dbtors}[1]{\operatorname{D}^{\operatorname{b}}(\operatorname{tors }#1)}
\newcommand{\dbqgr}[1]{\operatorname{D}^{\operatorname{b}}(\operatorname{qgr }#1)}
\newcommand{\dsing}[1]{\operatorname{D}_{\operatorname{sg}}(#1)}

\newcommand{\dbGcoh}[2]{\operatorname{D}^{\operatorname{b}}_{#2}(\operatorname{coh }#1)}
\newcommand{\dGsing}[2]{\operatorname{D}_{\operatorname{sg}}(#1,#2)}

\def\id{\mathop{\mbox{id}}}

\def\N{\mathop{\mathbb{N}}}
\def\Z{\mathbb{Z}}

\def\R{\mathbb{R}}

\def\F{\mathop{\mathcal{F}}}

\def\O{\mathcal{O}}

\def\P{\mathbb{P}}

\def\ra{\rightarrow}

\def\tritime{\text{\Clocklogo}}

\newcommand{\dgrsing}[1]{\operatorname{D}^{\operatorname{gr}}_{\operatorname{sg}}(#1)}

\def\id{\mathop{\mbox{id}}}

\def\N{\mathop{\mathbb{N}}}
\def\Z{\mathbb{Z}}

\def\R{\mathop{\mathbb{R}}}

\def\F{\mathop{\mathcal{F}}}

\def\O{\mathcal{O}}

\def\P{\mathbb{P}}

\def\ra{\rightarrow}

\def\tritime{\operatorname{\Clocklogo}}

\title{Orlov spectra: bounds and gaps}

\author[Ballard]{Matthew Ballard}
\address{
  \begin{tabular}{l}
   Matthew Ballard  \\ 
   \hspace{.1in} University of South Carolina, Department of Mathematics, Columbia, SC, USA \\
   \hspace{.1in} Email: {\bf ballard@math.sc.edu} \\
  \end{tabular}
}

\author[Favero]{David Favero}
\address{
  \begin{tabular}{l}
   David Favero \\
   \hspace{.1in} University of Alberta, Department of Mathematics, Edmonton, AB, Canada \\
   \hspace{.1in} Email: {\bf favero@gmail.com} \\
  \end{tabular}
}

\author[Katzarkov]{Ludmil Katzarkov}
\address{
  \begin{tabular}{l}
   Ludmil Katzarkov \\
   \hspace{.1in} University of Miami, Department of Mathematics, Coral Gables, FL, USA \\ 
   \hspace{.1in} Universit\"at von Wien, Fakult\"at f\"ur Mathematik,  Wien, \"Osterreich \\
   \hspace{.1in} Email: {\bf lkatzark@math.uci.edu} \\
  \end{tabular}
}

\numberwithin{equation}{section}
\begin{document}

%\titlerunning{}
%\authorrunning{Ballard, Favero, and Katzarkov}

\begin{abstract}
 The Orlov spectrum is a new invariant of a triangulated category. It was introduced by D. Orlov building on work of A. Bondal-M. van den Bergh and R. Rouquier. The supremum of the Orlov spectrum of a triangulated category is called the ultimate dimension. In this work, we study Orlov spectra of triangulated categories arising in mirror symmetry. We introduce the notion of gaps and outline their geometric significance. We provide the first large class of examples where the ultimate dimension is finite: categories of singularities associated to isolated hypersurface singularities. Similarly, given any nonzero object in the bounded derived category of coherent sheaves on a smooth Calabi-Yau hypersurface, we produce a new generator by closing the object under a certain monodromy action and uniformly bound this new generator's generation time.  In addition, we provide new upper bounds on the generation times of exceptional collections and connect generation time to braid group actions to provide a lower bound on the ultimate dimension of the derived Fukaya category of a symplectic surface of genus greater than one.
\end{abstract}

\maketitle

%\keywords{Orlov spectrum \and triangulated category}
% \subclass{
% 13A02 \and %Commutative algebra: Graded rings
% 13D05 \and %Commutative algebra: Homological dimension
% 13D09 \and %Commutative algebra: Derived categories
% 14A22 \and %Algebraic geometry: Noncommutative algebraic geometry
% 14E08 \and %Algebraic geometry: Rationality questions
% 14F05 \and %Algebraic geometry: Sheaves, derived categories of sheaves, and related constructions
% 14J33 \and %Algebraic geometry: Mirror symmetry
% 18E30 \and %Category theory: Derived categories, triangulated categories
% 53D37 \and %Differential geometry: Mirror symmetry, symplectic aspects; homological mirror symmetry; Fukaya category
% }

\section{Introduction}

The spectrum of a triangulated category was introduced by D. Orlov in \cite{O4}, building on work of  A. Bondal, R. Rouquier, and M. van den Bergh, \cite{Ro2} \cite{BV}. This categorical invariant, which we shall call the Orlov spectrum, is simply a list of non-negative numbers.  Each number is the generation time of an object in the triangulated category. Roughly, the generation time of an object is the necessary number of exact triangles it takes to build the category using this object. If the triangulated category is of geometric origin, like the bounded derived category of coherent sheaves on a scheme, the Orlov spectrum encodes nontrivial geometric information. In this paper, we study how geometry influences the structure of Orlov spectra and we find geometric meaning in the gaps arising in Orlov spectra. 

Although the (pre-)history of Orlov spectra extends back further, notably to work of A. Neeman, and Bondal-M. Kapranov, the fundamental background for this paper arose in \cite{BV}.  Here, Bondal and van den Bergh layout the foundations, introducing all of the notions necessary to define generation time.  %In this language, they prove the following representability result:
%\begin{theorem}[Bondal-van den Bergh]
% Let $\mathcal T$ be a Karoubi closed (i.e. every projector splits) triangulated category that is Ext-finite. If $\mathcal T$ has a strong generator %(i.e. a generator with finite generation time), then $\mathcal T$ is saturated.
%\end{theorem}
\noindent They apply their new notions to categories arising in algebraic geometry, proving a number of interesting and deep results that tie generators and geometry together.  Let us emphasize the following one:
\begin{theorem}
For a smooth scheme over a field, the bounded derived category of coherent sheaves admits a strong generator (i.e. a generator of finite generation time).
\end{theorem}
   
In \cite{Ro2}, Rouquier expanded on the foundations of \cite{BV}. He studied the minimal generation time amongst all strong generators, i.e. the minimum of the Orlov spectrum.  This notion we shall call the Rouquier dimension of a triangulated category. Rouquier proved many interesting results in \cite{Ro2} concerning the Rouquier dimension. His results had deep applications in both geometry and pure algebra.  Let us emphasize the following theorems which appear in loc.~cit.:
\begin{theorem}
 For a reduced separated scheme of finite type over a field, the Rouquier dimension of derived category of coherent sheaves is bounded below by the Krull dimension. 
\end{theorem}
\begin{theorem}
 For smooth quasi-projective schemes over a field, the Rouquier dimension of the derived category of coherent sheaves is bounded by twice the Krull dimension.
\end{theorem}
\noindent The following generalizes the above-mentioned result of Bondal and van den Bergh:
\begin{theorem}
 For any separated scheme of finite type over a field (not necessarily smooth), the derived category of coherent sheaves admits a strong generator.
\end{theorem}
\noindent Rouquier also contends that the the supremum amongst all generation times, which we shall call the ultimate dimension, should be studied in its own right. 

In \cite{O4}, Orlov utilizes results on the semi-stability of vector bundles on curves to prove the following interesting result:
\begin{theorem}
 The Rouquier dimension of the derived category of coherent sheaves on any smooth algebraic curve is one.
\end{theorem}
\noindent Having proven the one dimensional case, he proposes the following general conjecture:
\begin{conjecture}\label{conj: dimension} 
 For a smooth algebraic variety, $X$, the Krull dimension of $X$ and the Rouquier dimension of $\dbcoh{X}$ are equal.
\end{conjecture}
\noindent This conjecture asserts that Rouquier's notion of dimension of a triangulated category is deeply geometric.  Furthermore, Orlov contends that in order to extract additional, more novel, geometric invariants from the category, one should study all possible generation times - the Orlov spectrum.  With this in mind, he begins the analysis of the Orlov spectrum of a smooth algebraic curve proving:
\begin{theorem}
 The set, $\{1,2\}$, is a subset of Orlov spectrum of the bounded derived category of coherent sheaves on a smooth proper algebraic curve, with equality if and only if the curve is rational.
\end{theorem}
\noindent Orlov then promotes the following questions:
\begin{itemize}
\item Is the Orlov spectrum of the bounded derived category of coherent sheaves on a smooth quasi-projective scheme bounded above? Is it bounded above for a non-smooth scheme?
\item Does the Orlov spectrum of the bounded derived category of coherent sheaves on a (smooth) quasi-projective scheme form an integer interval?
\end{itemize}

Orlov's ideas were developed further by the first two authors in \cite{BF}.  In loc.~cit., the authors prove that generation time for tilting objects reduces to a very simply geometric calculation. They use it to prove Orlov's conjecture in many new cases. Through examples, they illustrate some subtleties encoded in generation time, including how it can vary in certain moduli and its relationship with positivity of the anti-canonical bundle.  

In addition to the papers discussed above, there are other works we should mention. Indeed, study of Orlov spectra, possibly proceeding under other names, seems a common endeavor across different algebraic fields. Rouquier's paper inspired further work in algebra by L. L. Avramov, P. A. Bergh, R.-O. Buchweitz, S. Iyengar, H. Krause, D. Kussin, C. Miller, and S. Oppermann, see \cite{ABIM}, \cite{BIKO}, \cite{KK}, \cite{Opp}, \cite{Opp2}.  Notably, \cite{BIKO} seems closely related to section 4 of this work. D. Benson, J. Carlson, S. Chebolu, J. D. Christensen, M. Hovey, K. Lockridge, Y. Min\'a\u{c}, and G. Puninski, see \cite{Loc}, \cite{HLP}, \cite{HL}, \cite{CCMa}, \cite{CCM09}, \cite{BCCM07}, are inspired by analogs of Freyd's Generating Hypothesis, which, in our language, seeks to determine whether an object has generation time zero.

Even with the wealth of knowledge detailed above, precise descriptions of Orlov spectra for, even simple, categories are still elusive. This paper builds on the growing understanding of the structure of Orlov spectra of categories of geometric origin, particularly categories of interest in mirror symmetry. The novelty of our current work lies in its approach to geometric themes encoded in the Orlov spectrum. Upper bounds on the ultimate dimension are closely tied to the Hochschild homology of the category. Lower bounds on the ultimate dimension are controlled by the complexity of braid groups actions. We expect that these phenomena, properly understood and synthesized, are universal.

We outline a new approach to decode the geometry found in the gaps of Orlov spectra. Gaps are simply the missing numbers in an Orlov spectrum. Their existence is precisely the content of Orlov's second question from above.  No matter their simplicity, the authors expect that gaps are a deep geometric invariant related to monodromy and capturing motivic information in the case of the derived category of coherent sheaves on a smooth proper variety.

Let us highlight our predictions by discussing some of the major themes of this work:

1) We provide the first large class of examples where the Orlov spectrum is bounded above: the category of singularities of an isolated hypersurface singularity. Our bound is expressed in terms of the embedding dimension and the nilpotence of the Tjurina algebra.
\begin{theorem}\label{thm:bounded spec}
 Let $(S, \mathfrak m_S)$ be an isolated hypersurface singularity. The Orlov spectrum of $\dsing{S}$ is bounded by $2(\dim S + 2)\op{LL}(S/(\partial w)) - 1$, where $\op{LL}$ denotes the Loewy length of an algebra.
\end{theorem}
\noindent We also entirely calculate the Orlov spectrum when $(S, \mathfrak m_S)$ is an $A_n$ singularity.  %Notice the reoccurring gaps are roughly quotients of $\op{LL}(S/(\partial w))$.
%\sidenote{I am not clear on what you mean by quotients of the Loewy length.}
\begin{theorem}
 The Orlov spectrum of $\dsing{A_{n-1}}$ is
\begin{displaymath}
  \left\lbrace \left\lceil \frac{\lfloor n/2 \rfloor}{s}  \right\rceil -1 \ : \ s \in \mathbb{N} \right\rbrace
\end{displaymath}
 where $\lfloor \alpha \rfloor$ is the greatest integer less than $\alpha$ and $\lceil \alpha \rceil$ is the least integer greater than $\alpha$.
\end{theorem}
Let us note that the results of \cite{Tak} can be applied to deduce that level of the residue field in $\dsing{S}$ with respect to any nonzero object of $\dsing{S}$ is bounded. This is an important step in the proof of Theorem~\ref{thm:bounded spec}.

2) The most unexpected geometric consequence is the connection of the theory of gaps of Orlov spectra to questions of rationality.  Based on work of Bondal, A. Kuznetsov and Orlov, rationality enters category theory by way of semi-orthogonal decompositions.  We demonstrate that gaps on certain intervals of the spectrum are obstructed by semi-orthogonal decompositions whose components have small Rouquier dimension.
\begin{theorem} 
Suppose $\langle \mathcal A_1 , \ldots , \mathcal A_n \rangle$ is a semi-orthogonal decomposition of $\mathcal T$ and $\mathcal G:= G_1 \oplus \cdots \oplus G_n$ is a generator of $\mathcal T$ with $G_i \in \mathcal A_i$. By performing a series of mutations to dual decompositions, we get a collection of generators. These generators give a sequence of elements in the Orlov spectrum, on which there is no gap greater than any of the generation times of the $G_i$.
\end{theorem}
\noindent Now, in light of the above theorem, we propose the following conjectures:
\begin{conjecture}\label{conj: gaps SO}
Let $X$ be a smooth algebraic variety and $\langle \mathcal A_1 , . . . , \mathcal A_n \rangle$ be a semi-orthogonal decomposition of $\dbcoh{X}$.  The length of any gap in $\dbcoh{X}$ is at most the maximal Rouquier dimension amongst the $\mathcal A_i$.
\end{conjecture}
\begin{conjecture}\label{conj: admissible gaps} 
Let $X$ be a smooth algebraic variety. If $\mathcal A$ is an admissible subcategory of $\dbcoh{X}$, then the length of any gap of $\mathcal A$ is at most the maximal length of any gap of $\dbcoh{X}$. Conversely, if $\mathcal A$ has a gap of length at least $s$, then so does $\dbcoh{X}$.
\end{conjecture}
\noindent These have the following nice corollaries:
\begin{corollary} 
Suppose Conjectures~\ref{conj: dimension} and~\ref{conj: gaps SO} hold.  If $X$ is a smooth variety then any gap of $\dbcoh{X}$ has length at most the Krull dimension of $X$.
\end{corollary}
\begin{corollary} 
Suppose Conjectures~\ref{conj: dimension}, \ref{conj: gaps SO}, and \ref{conj: admissible gaps} hold.  Let $X$ and $Y$ be birational smooth proper varieties of dimension $n$.  The category, $\dbcoh{X}$, has a gap of length $n$ or $n-1$ if and only if $\dbcoh{Y}$ has a gap of the same length i.e. the gaps of length greater than $n-2$ are a birational invariant.
\end{corollary}
\begin{corollary}
Suppose Conjectures~\ref{conj: dimension}, \ref{conj: gaps SO}, and \ref{conj: admissible gaps} hold.  If $X$ is a rational variety of dimension $n$, then any gap in $\dbcoh{X}$ has length at most $n-2$.
\end{corollary}
\noindent The above corollaries outline a new approach to dealing with questions of rationality with enormous potential towards applications.  In particular, based on work of Kuznetsov, we believe this could lead to a proof of non-rationality for a generic cubic fourfold. The mirror interpretation of this framework is discussed in \cite{KP} and \cite{KNS}.

3) The first two themes are related by work of Orlov (see \cite{Orl09}).  For smooth Fano hypersurfaces, the graded category of singularities of their affine cone is a semi-orthogonal component of the derived category of coherent sheaves.  Therefore, the Orlov spectrum of these components is also related to the Loewy length of the Tjurina algebra (which for a homogeneous polynomial is equal to the Milnor algebra) of their defining function.  In this case, this Loewy length is just $(d(n+1)-2n-1)$ by Macaulay's theorem.
\begin{theorem}
Let $f \in k[x_0,\ldots,x_n]$ be a homogeneous polynomial of degree $d$ and $A := k[x_0,\ldots,x_n]/(f)$. Assume that $A$ has an isolated singularity. For any non-zero object, $M$, in $\dgrsing{A}$, the object, $M \oplus M(1) \oplus \cdots \oplus M(d-1)$, is a generator of $\dgrsing{A}$ with generation time at most $2(n+1)(d(n+1)-2n-1) - 1$.
\end{theorem}
\noindent Orlov's work then provides us with the following geometric version:
\begin{corollary}
Let $X$ be a smooth hypersurface of degree $n+1$ in $\mathbb P^n$. Set $\tw{1} := L_{\mathcal O} \circ (-\otimes_{\mathcal O} \mathcal O(1))$ where $L_{\mathcal O}$ is the Seidel-Thomas twist by $\mathcal O$. For any nonzero $E \in \dbcoh{X}$, $E \oplus E\tw{1} \oplus \cdots \oplus E\tw{n}$ is a generator of $\dbcoh{X}$ with generation time bounded by $2n^2(n+1)-1$.
\end{corollary}
\noindent In light of the above discussion, we expect that gaps in the derived category of coherent sheaves on a Fano hypersurface in projective space are related to the structure of the corresponding Milnor algebra.

%\sidenote{The Loewy length is constant and depending only on the degree and dimension but some cubic fourfolds are rational while others are not. I think what you mean to say is that it depends on if there are objects on which these natural transformations do not vanish.}

4) We give a new upper bound on the generation time of any (full, not necessarily strong) exceptional collection. The upper bound comes from studying $A_{\infty}$-enhancements of triangulated categories. It ties in quite nicely with Koszul duality. 

\begin{theorem}
 Let $\mathcal A$ be a cohomologically-finite triangulated $A_{\infty}$-category possessing a (full) exceptional collection $A_1,\ldots,A_n$. The generation time of the dual collection in $H(\mathcal A)$ is bounded above by $\op{LL}_{\infty}(A')-1$ where $A'$ is a minimal $A_{\infty}$-algebra quasi-isomorphic to the $A_{\infty}$-endomorphism algebra of $\bigoplus_{i=1}^n A_i$. If the $A_{\infty}$-endomorphism algebra of $\bigoplus_{i=1}^n A_i$ is formal (quasi-isomorphic to its cohomology), then the generation time of the dual collection is equal to one less than the Loewy length of the cohomology of the $A_{\infty}$-endomorphism algebra of $\bigoplus_{i=1}^n A_i$.
\label{thm:excgentime}
\end{theorem}

Here, $\op{LL}_{\infty}$ is an extension of the notion of Loewy length to minimal $A_{\infty}$-algebras. In addition, we also provide examples to demonstrate how generation time can depend on ``higher homotopy'' information of the endomorphism algebra of an object, and we compute the Orlov spectra of the bounded derived categories of finite-dimensional representations of $A_n$ quivers.

\begin{theorem}
 Let $\text{Q}$ be a quiver such that the underlying graph is a Dynkin diagram of type $A_n$.  The Orlov spectrum of $\dbmod{kQ}$ is equal to the integer interval $\{0, \ldots , n-1\}$.
\end{theorem}

5) From the symplectic perspective, one can consider the Orlov spectrum as a new invariant of a Fukaya category.  Here we see that the correlation with monodromy theory is once again manifest by connecting generation time to braid group actions.  An upper bound, comes from a more well known construction and occurs as follows:
\begin{proposition}
Let $S_1,\ldots, S_n$ be spherical objects in the homotopy category, $\mathcal T$, of a triangulated cohomologically-finite $A_\infty$-category and assume we have $\emph{HH}^0(\mathcal T) = k$.  Suppose there exists a relation among the corresponding spherical twists:
\begin{displaymath}
L_{S_{a_1}} \cdots L_{S_{a_r}} \cong \op{Id}_{\mathcal T}
\end{displaymath}
with $1 \leq a_i \leq n$.  Then $S_1 \oplus \cdots \oplus S_n$ strongly generates $\mathcal T$ with generation time at most $r-1$.
\end{proposition}
\noindent Using a combination of braid relations and geometry, it is also possible to obtain lower bounds on generation time as in the following theorem:
\begin{theorem}
The ultimate dimension of the derived Fukaya category of a symplectic surface of genus $g$ is at least $4g$.
\end{theorem}
\noindent For the elliptic curve we calculate the Orlov spectrum in its entirety.  This result was also attained independently in unpublished work of Orlov.
\begin{theorem}
The Orlov spectrum of the bounded derived category of coherent sheaves on an elliptic curve is $\{ 1,2,3,4 \}$.
\end{theorem}
\noindent We expect derived Fukaya categories of symplectic surfaces to have no gaps, as is the case for a Riemann surface of genus one via homological mirror symmetry.  However, we suspect that there exists symplectic manifolds of real dimension four with large gaps, indicating that gaps of the Orlov spectrum of the derived Fukaya category is a nontrivial invariant of the symplectic motive.

The paper is organized as follows.  In Section 2, we establish our notational conventions and define all necessary mathematical notions revolving around the Orlov spectrum.  We proceed with a discussion of ghost maps, a theory originating in \cite{Kelly}, and used, implicitly and explicitly, by many subsequent authors in connection to this subject.  We illustrate a number of examples occurring in geometry, notably, spherical twists and monodromy of the quintic threefold.  In Section 3, we remind the reader of the basics of semi-orthogonal decompositions and demonstrate how semi-orthogonal decompositions whose components have small Rouquier dimension limit the size of gaps. We then outline how gaps in the Orlov spectrum of the bounded derived category of a variety can be used to answer questions about rationality. Finally, we develop a method, distinct from \cite{BF} and fully general, to bound the generation time of exceptional collections using the Loewy length of the dual collection. We provide a handful of examples to illustrate the utility of the method.  In Section 4, we discuss strong generators for categories of singularities of isolated singularities.  We provide new proofs, from our perspective, and extensions of some of the known results in this area.  We use these ideas to bound the Orlov spectrum of an isolated hypersurface singularity. In Section 5, we give a detailed recap of Orlov's theorem relating graded categories of singularities to bounded derived category of coherent sheaves.  We use our examination of Orlov's theorem and extensions of results from section 4 to study the Orlov spectrum for hypersurfaces in projective space.  Section 6, though connected to the other sections, can certainly be read independently.  Here, we illustrate the relationship between generation time and braid group actions, by means of the derived Fukaya category of a symplectic surface.  We compute the full Orlov spectrum of the elliptic curve and provide a lower bound on the ultimate dimension of derived Fukaya categories of a symplectic surface of higher genus.

\section{Preliminaries}
Throughout this article, $k$ denotes an algebraically-closed field of characteristic zero. All categories will be $k$-linear. For a ring, $R$, $\op{Mod }R$ denotes the category of right $R$-modules and $\op{D}(\op{Mod }R)$ denotes the unbounded derived category of right $R$-modules.  The bounded derived category of right $R$-modules we denote by $\op{D}^{\op{b}}(\op{Mod }R)$. For a Noetherian ring, $R$, $\op{mod }R$ denotes the category of finitely-generated right $R$-modules and $\dbmod{R}$ denotes its bounded derived category. If $X$ is a variety, $\dbcoh{X}$ denotes the bounded derived category of coherent sheaves on $X$.

Let $\mathcal T$ be a triangulated category.  For a full subcategory, $\mathcal I$, of $\mathcal T$ we denote by $\langle \mathcal I \rangle$ the full subcategory of $\mathcal T$ whose objects are isomorphic to summands of finite coproducts of shifts of objects in $\mathcal I$. In other words, $\langle \mathcal I \rangle$ is the smallest full subcategory containing $\mathcal I$ which is closed under isomorphisms, shifting, and taking finite coproducts and summands. For two full subcategories, $\mathcal I_1$ and $\mathcal I_2$, we denote by $\mathcal I_1 \ast \mathcal I_2$ the full subcategory of objects, $B$, such that there is a distinguished triangle, $B_1 \to B \to B_2 \to B_1[1]$, with $B_i \in \mathcal I_i$.  Set $\mathcal I_1 \diamond \mathcal I_2 := \langle \mathcal I_1 \ast \mathcal I_2 \rangle$, $\langle \mathcal I \rangle_0 :=\langle \mathcal I \rangle$, and inductively define 
\begin{displaymath}
 \langle \mathcal I \rangle_n := \langle \mathcal I \rangle_{n-1} \diamond \langle \mathcal I \rangle.
\end{displaymath}
Similarly we define
\begin{displaymath}
 \langle \mathcal I \rangle_{\infty} := \bigcup_{n \geq 0} \langle \mathcal I \rangle_{n}.
\end{displaymath}
For an object, $E \in \mathcal T$, we notationally identify $E$ with the full subcategory consisting of $E$ in writing, $\langle E \rangle_n$.  The reader is warned that, in some of  the previous literature, $\langle \mathcal I \rangle_0 := 0$ and $\langle \mathcal I \rangle_1 := \langle \mathcal I \rangle$. We follow the notation in \cite{BF}.  With our convention, the index equals the number of cones allowed. The operations, $\ast$ and $\diamond$, were introduced in \cite{BV} where their associativity is proven. From associativity, it follows that $\langle  \mathcal I \rangle_n \diamond \langle \mathcal I \rangle_m = \langle \mathcal I \rangle_{n+m+1}$. We will use this fact implicitly.

%\sidenote{What's the difference between these? Kuznetsov wanted it made clear which coproducts are allowed.}
We will need small modifications for the statement and proof of Proposition~\ref{prop:stronggenlocalize}. Let $\bar{\mathcal I}$ denote the smallest subcategory of $\mathcal T$ containing $\mathcal I$ and closed under $\mathcal T$-coproducts of objects of $\mathcal I$. Let $\tilde{\mathcal I}$ denote the smallest subcategory of $\mathcal T$ containing coproducts of the form, $\bigoplus_{a \in A} I$, for a single $I \in \mathcal I$ whenever $\bigoplus_{a \in A} I$ exists in $\mathcal T$. We then set $\langle \bar{\mathcal I} \rangle_0 = \langle \bar{\mathcal I} \rangle$ and
\begin{displaymath}
 \langle \bar{\mathcal I} \rangle_n := \overline{\langle \bar{\mathcal I} \rangle_{n-1} \diamond \langle \bar{\mathcal I} \rangle}.
\end{displaymath} 
We also set $\langle \tilde{\mathcal I} \rangle_0 = \langle \tilde{\mathcal I} \rangle$ and
\begin{displaymath}
 \langle \tilde{\mathcal I} \rangle_n := \widetilde{\langle \tilde{\mathcal I} \rangle_{n-1} \diamond \langle \tilde{\mathcal I} \rangle}.
\end{displaymath}

\begin{definition}
Let $E$ be an object of a triangulated category, $\mathcal{T}$.  If there is an $n$ with $\langle E \rangle_{n} = \mathcal T$, we set
\begin{displaymath}
 \tritime_{\mathcal T}(E):=  \text{min } \lbrace  n \geq 0 \  | \ \langle E \rangle_{n} = \mathcal T \rbrace.
\end{displaymath}
Otherwise, we set $\tritime_{\mathcal T}(E) := \infty$.   We call $\tritime_{\mathcal T}(E)$ the \textbf{generation time} of $E$. When, $\mathcal T$ is clear from context, we omit it and simply write $\tritime(E)$.  If $\langle E \rangle_{\infty}$ equals $\mathcal{T}$, we say that $E$ is a \textbf{generator}. If $\tritime(E)$ is finite, we say that $E$ is a \textbf{strong generator}. The \textbf{Orlov spectrum} of $\mathcal T$, denoted $\op{OSpec}\mathcal T$, is the set
\begin{displaymath}
 \lbrace \tritime(G)  \ | \ G \in \mathcal T, \ \tritime(G) < \infty \rbrace \subset \mathbb{Z}_{\geq 0}.
\end{displaymath}
The \textbf{Rouquier dimension} of $\mathcal T$, denoted $\text{rdim }\mathcal T$, is the infimum of $\ops{OSpec}\mathcal T$, it is defined as $\infty$ when $\ops{OSpec}\mathcal T$ is empty. The \textbf{ultimate dimension} of $\mathcal T$, denoted $\op{udim}\mathcal T$, is the supremum of $\op{OSpec}\mathcal T$ it is defined as $\infty$ when $\ops{OSpec}\mathcal T$ is empty.
\end{definition}

We shall denote the Orlov spectrum, Rouquier dimension, and ultimated dimension of $\dbcoh{X}$ by $\op{OSpec } X$, $\ops{rdim}X$, and $\ops{udim}X$, respectively. It is also convenient to recall the following definition which first appeared in \cite{ABIM}. 
\begin{definition}
Let $E$ be an object of a triangulated category, $\mathcal{T}$.  If there is an $n$ with $A \in \langle E \rangle_{n}$, we set
\begin{displaymath}
 \text{Lvl}^E_{\mathcal T}(A):=  \text{min } \lbrace  n \geq 0 \  | A \in \langle E \rangle_{n} \rbrace.
\end{displaymath}
Otherwise, we set $\text{Lvl}_{\mathcal T}^E(A) = \infty$.  This number is called the \textbf{level of $A$ with respect to $E$}, or simply the level of $A$ when $E$ is implicit.  
\end{definition}

The case where $\mathcal T$ is $\op{D}^{\op{b}}(\op{mod} A)$, the bounded derived category of coherent modules, and $G=A$ is the free module, provides some insight into the formalism above. The following theorem is taken from \cite{KK}; the cases where $A$ is finite-dimensional over $k$ or commutative and essentially of finite type were proven in \cite{Ro2}.

\begin{theorem} 
 Let $A$ be a right-coherent $k$-algebra.  The generation time of $A$, as an object of 	$\emph{D}^{\emph{b}}(\emph{mod }A)$, is equal to the global dimension of $A$.
\label{CKK}
\end{theorem}

\begin{remark}
 Using ideas from \cite{ABIM}, one can extend the notion of global dimension to dg-algebras in a natural manner and check that the analog of Theorem \ref{CKK} holds for dg-algebras. As noted in \cite{Ro2}, for an enhanceable triangulated category, $\mathcal T$, each generator, $G$, allows one to construct an equivalence of $\mathcal T$ with the derived category of perfect dg-modules over the dg-endomorphisms of $G$. In this way, the Orlov spectrum can be viewed as a list of global dimensions of dg-algebras within a derived Morita equivalence class.
\end{remark}

We have the following simple lemma whose proof is left to the reader, see \cite{BF}:
\begin{lemma} 
Let $F: \mathcal T \to \mathcal R$ be an exact functor between triangulated categories.  Let $G$ be an object of $\mathcal T$. If $B \in \langle G \rangle_n$, then $F(B) \in \langle F(G) \rangle_n$.  Moreover, if $F$ commutes with coproducts and $B \in {\langle \overline{G} \rangle}_n$, then $F(B) \in \langle \overline{F(G)} \rangle_n$.
\label{lem:functor invariant}
\end{lemma}

% \begin{proof}
%  Any exact functor commutes with finite coproducts and takes exact triangles to exact triangles so $F \left(\langle G \rangle_n\right) \subset \langle F(G) \rangle_n$. To get the identity, $F \left(\langle \overline{G} \rangle_n\right) \subset \langle \overline{F(G)} \rangle_n$, we need to assume that $F$ respects coproducts.
% \end{proof}

Let $F: \mathcal{T} \to \mathcal{R}$ be an exact functor between triangulated categories. If every object in $\mathcal{R}$ is isomorphic to a direct summand of an object in the image of $F$, we say that $F$ is \textbf{dense}, or has dense image.

\begin{lemma} 
  If $F: \mathcal{T} \to \mathcal{R}$ has dense image and $G$ be a strong generator, then $\tritime(G) \geq \tritime(F(G))$.  In particular, $\emph{dim }\mathcal{T} \geq \emph{dim }\mathcal{R}$. 
\label{density lemma}
\end{lemma}
Again, the proof is left an exercise to the reader, see \cite{BF}.
% \begin{proof}
%  If $G$ is a generator of $\mathcal T$ with minimal generation time $t$, then $\mathcal T = \langle G \rangle_t$. We apply $F$ and use Lemma~\ref{lem:functor invariant} to get $F(\mathcal T) \subset \langle F(G) \rangle_t$. Since every object of $\mathcal R$ is a summand of an object $F(\mathcal T)$, we see that $\mathcal R = \langle F(G) \rangle_t$. Thus, $\dim \mathcal R \leq t$.
% \end{proof}

% \begin{example}
% Let $i: U \to X$ be the inclusion of an open subvariety into an algebraic variety, $X$.  Then $i^*$ is dense.  Hence, $\ops{rdim}\dbcoh{X} \geq \ops{rdim}\dbcoh{U}$.
% \end{example}
% 
% \begin{example} 
% Suppose $i: C \to X$ is an embedding of a smooth curve in a smooth variety, $X$.  Since the category of coherent sheaves on $C$ is hereditary, any complex in $\dbcoh{C}$ is a sum of it's cohomology sheaves, \cite{H}.  Now given any coherent sheaf $\mathcal F$, we have $\mathcal F \cong H^0(i^*i_*(\mathcal F))$, where the functors are derived.  It follows that the derived functor $i^*: \dbcoh{X} \to \dbcoh{C}$ is dense.
% \label{dense curve}
% \end{example}

\begin{example} 
Let $V$ be a vector bundle in $\dbcoh{X}$.  Then the functor $(-\otimes_{\mathcal O} V): \dbcoh{X} \to \dbcoh{X}$ is dense, as any object, $F$, is a summand of $(F \otimes_{\mathcal O} V^\vee) \otimes_{\mathcal O} V$.  
\label{dense bundle}
\end{example}

\begin{example} 
Consider a finite group $\Gamma$ acting on an algebraic variety, $X$, and consider the derived category of coherent sheaves on $X$, $\dbcoh{X}$, and the derived category of $\Gamma$-equivariant coherent sheaves on $X$, $\dbGcoh{X}{\Gamma}$. We have two derived functors: the forgetful functor, $\text{For}:\dbGcoh{X}{\Gamma} \to \dbcoh{X}$, and the inflation functor, $\text{Inf}: \dbcoh{X} \to \dbGcoh{X}{\Gamma}$, where by definition $\text{Inf}(A) = \bigoplus_{g \in \Gamma} g^* A$ with the natural $\Gamma$ action.

Notice that any $A \in \dbcoh{X}$ is a summand of $\text{For}(\text{Inf}(A))$, hence the forgetful functor is dense.  On the other hand, for any $B \in \dbGcoh{X}{\Gamma}$ and each $g \in \Gamma$, we have an isomorphism, $\phi_g : g^*\text{For}(B) \to \text{For}(B)$, in $\dbcoh{X}$ coming from the equivariant structure on $B$. $\text{For}$ is the left adjoint to $\text{Inf}$ with adjunction morphism in $\dbGcoh{X}{\Gamma}$ defined by:
\begin{displaymath}
 \sum_{g \in \Gamma} \phi_g: \text{Inf}(\text{For}(B)) \to B.
\end{displaymath}
The map,
\begin{displaymath}
\frac{1}{|\Gamma|}\bigoplus_{g \in \Gamma} \phi^{-1}_g: B \to \text{Inf}(\text{For}(B)),
\end{displaymath}
provides a splitting of the map above.  Therefore, $B$ is a summand of $\text{Inf}(\text{For}(B))$, and the functor $\text{Inf}$ is also dense.

Hence, for any generator, $G$, of $\dbcoh{X}$, we have:
\begin{displaymath}
\tritime(\text{For}(\text{Inf}(G)) \leq \tritime(\text{Inf}(G)) \leq \tritime(G).
\end{displaymath}
It follows that $\dbcoh{X}$ and $\dbGcoh{X}{\Gamma}$ have the same Rouquier dimension.  Furthermore, for any generator, $G$, of $\dbcoh{X}$ which is invariant under the action of $\Gamma$, we have $\langle G \rangle_0 =  \langle \text{For}(\text{Inf}(G)) \rangle_0$ hence $\tritime(\text{For}(\text{Inf}(G))) = \tritime(G)$ and thus $\tritime(G) = \tritime(\text{Inf}(G))$.
\label{dense equivariant}
\end{example}

\begin{lemma} 
 If $\mathcal{T}$ is a triangulated category with finite Rouquier dimension, then any generator is a strong generator.
 \label{lem:autostrong}
\end{lemma}
The proof is left an exercise to the reader, see \cite{BF}.
% \begin{proof}
%  Let $X$ be a generator of $\mathcal T$.  As $\mathcal{T}$ is finite dimensional, there exists a strong generator, $G$, with $\langle G \rangle_n = \mathcal T$. Since $X$ generates, $G \in \langle X \rangle_t$ for some $t$.  Hence $\langle X \rangle_{(n+1)(t+1)-1} = \mathcal T$ by the associativity of $\diamond$, \cite{BV}.
% \end{proof}

The generation time of an object can be reinterpreted in terms of so called ``ghost maps;'' this reinterpretation turns out to be quite useful both for intuition about generation time and as a means of calculation.

%\sidenote{Kuznetsov didn't like coghost, should we switch to left and right ghost?}
\begin{definition}
Let $\mathcal T$ be a triangulated category, $f$ be a morphism, and $\mathcal I$ be a full subcategory.  We say that $f$ is \textbf{$\mathcal I$ ghost} if, for all $I \in \mathcal I$, the induced map, $\text{Hom}_{\mathcal T}(I,X) \to \text{Hom}_{\mathcal T}(I,Y)$, is zero.  We say that $f$ is \textbf{$\mathcal I$ co-ghost} if, for all $I \in \mathcal I$, the induced map, $\text{Hom}_{\mathcal T}(Y,I) \to \text{Hom}_{\mathcal T}(Y,I)$, is zero.  If $G$ is an object of $\mathcal T$, we will say that $f$ is $G$ ghost if $f$ is $\langle G \rangle_0$ ghost and $f$ is $G$ co-ghost if $f$ is $\langle G \rangle_0$ co-ghost. 
\end{definition}

\begin{remark}
 Note that $\mathcal I$ ghosts and $\mathcal I$ co-ghosts naturally form ideals in $\mathcal T$.
\end{remark}

The following lemmas relate generation time to ghost maps and are a crucial ingredient in our study of Orlov spectra.  Lemma~\ref{lem:ghost} first appeared in \cite{Kelly} and later appeared in many places, for example see \cite{KK} \cite{Ro2}.

\begin{lemma}
 Let $\mathcal T$ be a triangulated category and let $G$ be an object of $\mathcal T$. If there exists a sequence of morphisms, $f_i: X_{i-1} \to X_i$, $1\leq i \leq t$, in $\mathcal T$ where each $f_i$ is $G$ ghost and $f_t \circ \cdots \circ f_1 \not = 0$, then $X_0 \not \in \langle G \rangle_{t-1}$. 
\label{lem:ghost}
\end{lemma}

\begin{proof}
 Let us show that $f_t \circ \cdots \circ f_1$ is ghost for $\langle G \rangle_{t-1}$. For simplicity, set $f^t:= f_t \circ \cdots \circ f_1$. We proceed by induction with the case, $t=1$, clear. Assume we know $f^t$ is $\langle G \rangle_{t-1}$ ghost for $t \leq n-1$, and let us consider the case $t=n$. From the induction hypothesis, $f^{n-1}$ is $\langle G \rangle_{n-2}$ ghost. Let $Y$ be an object of $\mathcal T$ lying in a triangle
\begin{displaymath}
 Z \overset{\alpha}{\to} Y \overset{\beta}{\to} Y_G \to Z[1]
\end{displaymath}
 with $Z \in \langle G \rangle_{n-2}$ and $Y_G \in \langle G \rangle_0$. Take any map $g: Y \to X_0$. As $f^{n-1}$ is $\langle G \rangle_{n-2}$ ghost, the composition $f^{n-1} \circ g \circ \alpha$ vanishes. Thus, we have a map $h: Y_G \to X_{n-1}$ with $f^{n-1} \circ g = h \circ \beta$. As $f_n$ is $\langle G \rangle_0$ ghost, $f_n \circ h \circ \beta = f^n \circ g$ vanishes. Thus, $f^n$ is $\langle G \rangle_{n-2} \ast \langle G \rangle_0$ ghost. It is clear this implies that $f^n$ is $\langle G \rangle_{n-1}$ ghost.

 To finish the proof the lemma, note that, if $X_0$ lies in $\langle G \rangle_{t-1}$, then $f^t \circ \op{id}_{X_0} = f^t$ vanishes. 
\end{proof}

We also have the dual statement whose proof is the same.

\begin{lemma}
 Let $\mathcal T$ be a triangulated category and let $G$ be an object of $\mathcal T$. If there exists a sequence of morphisms, $f_i: X_{i-1} \to X_i$, $1\leq i \leq t$, in $\mathcal T$ where each $f_i$ is $G$ co-ghost and $f_t \circ \cdots \circ f_1 \not = 0$, then $X_t \not \in \langle G \rangle_{t-1}$.
 \label{lem:co-ghost} 
\end{lemma}

The following partial converse seems well-known but first appeared in the literature in S. Oppermann's thesis \cite{Opp2}:

\begin{lemma}
 Let $\mathcal T$ be a triangulated category and let $G$ be an object of $\mathcal T$. Assume that for any object, $X$, of $\mathcal T$ there exists a morphism, $\nu_X: X_G \to X$, with $X_G \in \langle G \rangle_0$ and satisfying the following condition: for any morphism, $g: Y \to X$, with $Y \in \langle G \rangle_0$, there exists a morphism, $h: Y \to X_G$, with $g = \nu_X \circ h$. If $X \not \in \langle G \rangle_{t-1}$, then there exists a sequence of morphisms, $f_i: X_{i-1} \to X_i$, $1\leq i \leq t$, in $\mathcal T$ where each $f_i$ is $G$ ghost, $X_0 = X$ and $f_t \circ \cdots \circ f_1 \neq 0$.
 \label{lem:ghostconverse}
\end{lemma}

\begin{proof}
 Complete $\nu_X: X_G \to X$ to a distinguished triangle
\begin{displaymath}
 X_G \overset{\nu_X}{\to} X \overset{f_1}{\to} X_1 \to X_G[1].
\end{displaymath}
 $f_1$ is $G$ ghost. Now iterate to get triangles
\begin{displaymath}
 \left(X_i\right)_G \overset{\nu_{X_i}}{\to} X_i \overset{f_{i+1}}{\to} X_{i+1} \to \left(X_i\right)_G[1].
\end{displaymath}
 If the composition $f_t \circ \cdots \circ f_1$ vanishes, then repeated application of the octahedral axiom exhibits $X \in \langle G \rangle_{t-1}$.
\end{proof}

We also have the dual statement whose proof is the same.

\begin{lemma}
 Let $\mathcal T$ be a triangulated category and let $G$ be an object of $\mathcal T$. Assume that for any object of $X$ of $\mathcal T$ there exists a morphism, $\nu_X: X \to X_G$, with $X_G \in \langle G \rangle_0$ and satisfying the following condition: for any morphism, $g: X \to Y$, with $Y \in \langle G \rangle_0$, there exists a morphism, $h: X_G \to Y$, with $g = h \circ \nu_X$. If $X \not \in \langle G \rangle_{t-1}$, then there exists a sequence of morphisms, $f_i: X_{i-1} \to X_i$, $1\leq i \leq t$, in $\mathcal T$ where each $f_i$ is $G$ co-ghost and $X_t = X$.
\label{lem:co-ghostconverse}
\end{lemma}

\begin{remark}
 We can replace $G$ by a general subcategory, $\mathcal I$, in each of these statements. However, we should note that it is necessary to assume the existence of an ``$\mathcal I$-approximation'' similar to the hypotheses of Lemmas~\ref{lem:ghostconverse} and~\ref{lem:co-ghostconverse}. If $X$ is a projective variety, then there are no $\op{Perf} X$ ghosts in $\dbcoh{X}$, see \cite{Bal}, and $\op{Perf} X$ is not dense in $\dbcoh{X}$, for a general $X$.
\end{remark}

Recall that a triangulated category, $\mathcal T$, is Ext-finite, if for any pair of objects, $A$ and $B$, of $\mathcal T$, we have
\begin{displaymath}
 \dim_k \left( \bigoplus_{l \in \Z} \op{Hom}_{\mathcal T}(A,B[l]) \right) < \infty.
\end{displaymath}
Combining the previous observations, we get the following corollary, which cannot be called anything other than a lemma:

\begin{lemma} [Ghost/Co-ghost Lemma and Converse] 
Let $\mathcal T$ be a $k$-linear Ext-finite triangulated category and let $G$ and $X_0$ be objects in $\mathcal T$.  The following are equivalent:
\renewcommand{\labelenumi}{\emph{\roman{enumi})}}
\begin{enumerate}
\item one has $X_0 \in \langle G \rangle_n$ and $X_0 \notin \langle G \rangle_{n-1}$;

\item there exists a sequence, %for $1 \leq i \leq n$ there exists a sequence of morphisms $\{ f_i \}$ and objects $X_i$ in $D^b(\mathcal C)$
\begin{displaymath}
\begin{CD}
X_0      @>f_1>> X_1    @>f_2>> \cdots @>f_{n-1}>> X_{n-1} @>f_n>> X_n,   \\
\end{CD}
\end{displaymath}
of maps in $\mathcal T$ such that all the $f_i$ are ghost for $G$ and $f_n \circ \cdots \circ f_1 \neq 0$.  Furthermore there is no such sequence for $n+1$.

\item there exists a sequence, %for $1 \leq i \leq n$ there exists a sequence of morphisms $\{ f_i \}$ and objects $X_i$ in $D^b(\mathcal C)$
\begin{displaymath}
\begin{CD}
X_n      @>f_{n}>> X_{n-1}    @>f_{n-1}>> \cdots @>f_2>> X_{n-1} @>f_1>> X_0,   \\
\end{CD}
\end{displaymath}
of maps in $\mathcal T$ such that all the $f_i$ are co-ghost for $G$ and $f_1 \circ \cdots \circ f_n \neq 0$.  Furthermore there is no such sequence for $n+1$.

\item there exists a sequence, %for $1 \leq i \leq n$ there exists a sequence of morphisms $\{ f_i \}$ and objects $X_i$ in $D^b(\mathcal C)$
\begin{displaymath}
\begin{CD}
X_0      @>f_1>> X_1    @>f_2>> \cdots @>f_{n-1}>> X_{n-1} @>f_n>> X_n,   \\
\end{CD}
\end{displaymath}
of maps in $\mathcal T$ with indecomposable objects, $X_i \in \mathcal T$, such that all the $f_i$ are ghost for $G$ and $f_n \circ \cdots \circ f_1 \neq 0$.  Furthermore, there is no such sequence for $n+1$.

\item there exists a sequence, %for $1 \leq i \leq n$ there exists a sequence of morphisms $\{ f_i \}$ and objects $X_i$ in $D^b(\mathcal C)$
\begin{displaymath}
\begin{CD}
X_n      @>f_{n}>> X_{n-1}    @>f_{n-1}>> \cdots @>f_2>> X_{n-1} @>f_1>> X_0,   \\
\end{CD}
\end{displaymath}
of maps in $\mathcal T$ with indecomposable objects, $X_i \in \mathcal T$, such that all the $f_i$ are co-ghost for $G$ and $f_1 \circ \cdots \circ f_n \neq 0$.  Furthermore, there is no such sequence for $n+1$.
\end{enumerate}
\label{ghost lemma}	
\end{lemma}

\begin{proof}
 $\mathcal T$ satisfies the hypothesis of Lemma~\ref{lem:ghostconverse}. Let $X$ be an object of $\mathcal T$. We set $X_G = \bigoplus_{l \in \Z} \op{Hom}_{\mathcal T}(G,X[l]) \otimes_k G[-l]$ and let $\nu_X: X_G \to X$ be the evaluation map. Similarly, $\mathcal T$ satisfies the hypothesis of Lemma~\ref{lem:co-ghostconverse}. The equivalence of $i),ii),iii)$ is a combination of Lemmas~\ref{lem:ghost}, \ref{lem:co-ghost}, \ref{lem:ghostconverse}, and \ref{lem:co-ghostconverse}. The only difference between $ii)$ and $iv)$ is that the objects are assumed to be indecomposable, their equivalence is clear.  The same goes for $iii)$ and $v)$.
\end{proof}

We have an important special case.  Recall that a hereditary abelian category is one where $\text{Ext}^2(A,B) = 0$ for any two objects, $A$ and $B$.

\begin{lemma}
 Let $\mathcal C$ be a hereditary abelian category with finite dimensional morphism spaces and let $G$ be an object of $\op{D}^{\op{b}}(\mathcal C)$ and $X_0$ be an object of $\mathcal C$. The following are equivalent:
\renewcommand{\labelenumi}{\emph{\roman{enumi})}}
\begin{enumerate}
 \item one has $X_0 \in \langle G \rangle_n$ and $X_0 \notin \langle G \rangle_{n-1}$;
 \item there exists a sequence,
%$s$ and $t$ with $s+t=n$ and for $1 \leq i \leq s$ and $1 \leq j \leq t$ there are morphisms in $\op{D}^{\op{b}}(\mathcal C)$ with indecomposable objects,  in $\mathcal C$
\begin{displaymath}
\begin{CD}
X_0      @>g_1>> \cdots @>g_s>> X_s @>h_1>> Y_1[1] @>h_2>> \cdots@>h_t>> Y_t[1],   \\
\end{CD}
\end{displaymath}
of maps in $\op{D}^{\op{b}}(\mathcal C)$ with $X_i$ and $Y_i$ indecomposable objects of $\mathcal C$, $s+t=n$, and such that all the $f_i$ and $g_i$ are ghost for $G$ and $h_t \circ \cdots \circ g_1 \neq 0$. Furthermore, there is no such sequence for $n+1$. 
 \item there exists a sequence,
%$s$ and $t$ with $s+t=n$ and for $1 \leq i \leq s$ and $1 \leq j \leq t$ there are morphisms in $D^b(\mathcal C)$ with indecomposable objects,  in $\mathcal C$
\begin{displaymath}
\begin{CD}
Y_t[-1]      @>h_t>> \cdots @>h_2>> Y_0 @>g_s>> X_s @>g_{s-1}>> \cdots@>g_1>> X_0,   \\
\end{CD}
\end{displaymath}
of maps in $\op{D}^{\op{b}}(\mathcal C)$ with $X_i$ and $Y_i$ indecomposable objects of $\mathcal C$, $s+t=n$, and such that all the $f_i$ and $g_i$ are co-ghost for $G$ and $g_1 \circ \cdots \circ h_t \neq 0$. Furthermore, there is no such sequence for $n+1$. 
\end{enumerate}
\label{lem:ghosthered}
\end{lemma}
Recall that for a finite dimensional algebra, $A$, with nilradical, $N$, the \textbf{Loewy length}, denoted $\op{LL}(A)$, is smallest $n$ such that $N^n =0$.
\begin{corollary} 
Suppose $\mathcal C$ is a $k$-linear hereditary category with with finite dimensional morphism spaces and finitely many isomorphism classes of indecomposable objects. Let $M_i$ be chosen representatives the isomorphism classes.  Then, $\op{udim}\mathcal T \leq \op{LL}(\mathbf{R}\emph{End}(\oplus M_i))-1$.
\label{cor:finitely many objects bound}
\end{corollary}

There is an important relationship between ghost maps and Serre functors. Let us recall the definition of a Serre functor:
\begin{definition}
A $k$-linear exact autoequivalence, $S$, of $\mathcal{T}$, is called a \textbf{Serre functor} if
for any pair of objects, $X$ and $Y$, of $\mathcal{T}$, there exists an isomorphism of vector spaces,
$$\text{Hom}_{\mathcal T}(Y,X)^\vee \cong \text{Hom}_{\mathcal T}(X, S(Y)),$$
which is natural in $X$ and $Y$.
\end{definition}
A Serre functor, if it exists, is determined uniquely up to natural isomorphism. If $F: \mathcal{T} \ra \mathcal{S}$ is an exact equivalence of triangulated categories possessing Serre functors, then $F$ commutes with those Serre functors \cite{BK}.  Now, recall that a category is called Karoubi closed if all idempotents split. Suppose $\mathcal T$ is a $k$-linear \textbf{Karoubi closed} triangulated category with finite dimensional morphism spaces which admits a Serre functor, $S$.  Let $X$ be an indecomposable object of $\mathcal T$.  In this situation, there is a natural map, $\epsilon_X: X \to S(X)$, corresponding to, 
\begin{displaymath}
\op{Hom}_{\mathcal T}(X,X) \to \op{Hom}_{\mathcal T}(X,X)/\op{Rad}_{\mathcal T}(X,X) \cong k,
\end{displaymath}
where the isomorphism with the base field identifies the image of the identity with $1$.  By definition of a Serre functor, there is also a nondegenerate pairing,
\begin{displaymath}
\text{Hom}_{\mathcal T}(A,B) \otimes_k \text{Hom}_{\mathcal T}(B,S(A)) \to k.
\end{displaymath}
Hence any nonzero morphism, $X \to A$, can be extended to a nonzero morphism, $X \to A \to S(X)$, the total morphism in this situation can be taken to be the natural map described above, see \cite{RV}.

\begin{proposition}\label{prop:extendwithSerre}
 Let $\mathcal T$ be a $k$-linear triangulated Karoubi closed category with finite-dimensional morphism spaces. Assume $\mathcal T$ possesses a Serre functor, $S$. Let $X$ be an indecomposable object in $\mathcal T$ and $f: X \to Y$ a morphism. There exists a morphism, $g: Y \to S(X)$, so that $g \circ f = \epsilon_X$.
\end{proposition}
%\begin{proof}
%Consider the following commutative diagram:
%$$\begin{CD}
%\text{Hom}(Y,S(X))    @>>>   \text{Hom}(X,S(X)) \\
%@|   @|    \\
%\text{Hom}(X,Y)^\vee    @>>>   \text{Hom}(X,X)^\vee.
%\end{CD}$$
%\end{proof}
Given any nonzero ghost sequence, $X \overset{f_1}\to \cdots \overset{f_n}\to X_n$ with $f_n \circ \cdots \circ f_1 \not = 0$, we can extend it to a new sequence, $X \overset{f_1}\to \cdots \overset{f_n}\to X_n \overset{g}\to S(X)$, with $g \circ f_n \circ \cdots f_1 = \epsilon_X$ and where only $g$ is possibly non-ghost.  Concatenating $f_n$ with $g$, we get a ghost sequence of equal length beginning at $X$ and terminating at $S(X)$.

Now, for any map, $G \to X$, consider the following commutative diagram:
$$\begin{CD}
\text{Hom}(X,S(X))    @>>>   \text{Hom}(G,S(X)) \\
@VV{\cong}V   @VV{\cong}V    \\
\text{Hom}(X,X)^\vee    @>>>   \text{Hom}(X,G)^\vee.
\end{CD}$$
By duality, requiring that the image  of $\epsilon_X$ is nonzero in $\text{Hom}(X,G)^\vee$ is equivalent to requiring that $\text{Hom}(X,G) \to \text{Hom}(X,X)$ does not lie in $\text{Rad}(X,X)$.  Hence, $G \to X$ has a section.  Meaning that if $G \to X \overset{\epsilon_X}\to S(X)$ is nonzero than $X$ is a summand of $G$.  One can similarly show that, for any map, $S(X) \to G$, if the composition, $X \to S(X) \to G$, is nonzero, then $G$ is a summand of $S(X)$.  Therefore, $\epsilon_X$ composed with any map besides a sequence of split epimorphisms and/or monomorphisms is zero. In other words, the natural map, $\epsilon_X$, is $G$ ghost and $G$ co-ghost for any object, $G$, of which $X$ is not a summand. Hence, given a ghost sequence whose total map is $\epsilon_X$, it can not be extended any further (although it could be perhaps factored into more maps). 

Ghost maps often have geometric origins. We collect some examples here.

\begin{example}[Central actions as ghosts] 
Let $\mathcal T$ be a triangulated category. The center of $\mathcal T$, denoted $Z(\mathcal T)$, is the space of natural transformations from $\op{Id}_{\mathcal T}$ to $\op{Id}_{\mathcal T}$. Let $x$ be an element of $Z(\mathcal T)$. If $G$ is an object of $\mathcal T$ with $x(G) = 0$, then we say that $x$ annihilates $G$.  For any object, $A \in \mathcal T$, and any morphism, $\alpha: G \to A[i]$, we have the following commutative diagram:  
\begin{displaymath}
\begin{CD}
G    @>x(G)>>   G \\
@V\alpha VV   @      V\alpha VV    \\
A[i]    @>x(A)>>   A[i].
\end{CD}
\end{displaymath}
Since $x(G)=0$ by assumption, $x(A) \circ \alpha =0$ for all $\alpha \in \text{Hom}_{\mathcal T}(G,A[i])$.  In other words, $x(A)$ is $G$ ghost for any object $A \in \mathcal T$. Similarly, $x(A)$ is $G$ co-ghost. If $\mathcal T = \dbcoh{X}$ for a quasi-projective variety, $X$, then $Z(\dbcoh{X}) \cong \Gamma(X,\mathcal O_X)$, \cite{Ro1}.
\label{center ghost}
\end{example}

\begin{example}[Divisors and ghosts]
The choice of a divisor, $i:D \to X$, gives a natural transformation, $\alpha: \op{Id}_{\dbcoh{X}} \to (- \otimes_{\mathcal O} \O(D))$.  Let $H$ be in the essential image of the functor $i_* : \dbcoh{D} \to \dbcoh{X}$.  Then, for any object, $A \in \dbcoh{X}$, $\alpha(A)$ is $H$ ghost by adjunction.
\label{divisors and ghosts}
\end{example}

\begin{example}[Tangent vectors as ghosts]
Let $X$ be a variety of dimension $n$.    Let $G$ be any object of $\dbcoh{X}$ and consider a smooth point, $p$, at which the cohomology sheaves of $G$ are locally free.  It is easily verified that any tangent vector, $\zeta \in \text{Hom}(\O_p, \O_p[1]) \cong T_pX$, is $G$ ghost.  Now take a basis for the tangent space, $\zeta_1, \ldots, \zeta_n$.  The composition is a nonzero ghost sequence: 
\[
\O_p \overset{\zeta_1}\to \O_p[1] \to \cdots \to \O_p[n-1] \overset{\zeta_n[n-1]}\to \O_p[n].
\]
It follows from the Ghost Lemma~\ref{ghost lemma}, that $n \leq \tritime(G)$.  Hence, $n \leq \ops{rdim}X$.  This proof is due to Rouquier and can be found in \cite{Ro2}.
\end{example}

\begin{example}[Cycles and levels]\label{eg:cyclesandlevels}
% Let us combine this type of ghost map with Example~\ref{dense curve}.  Recall,  we have a smooth curve $i: C \to X$ in our smooth variety, $X$.  Since, $i^*$ is dense, it follows that $i_*$ is faithful.  By adjunction, for any $i^*G$ ghost map $A \to B$, we have a $G$ ghost, $i_*A \to i_*B$.  For a point $p \in C$, take any ghost sequence beginning at $\O_p \to A_1 \to \cdots \to A_n$.  By nondegeneracy of the Serre pairing (as mentioned above) we may assume $A_n = O_p[1]$.  Consider the total composition $f: \O_p \to \O_p[1]$.  The pushforward $i_*f$ is nothing more than the tangent vector to the curve $C$ in $\text{Hom}(\O_p, \O_p[1])$.

We can extend the previous example a bit more. Let $i: V \to X$ be a smooth subvariety of $X$. By adjunction, the pushforward of any $\mathbf{L}i^*G$ ghost is $G$ ghost. For any point, $p \in V$, take any $\mathbf{L}i^*G$ ghost sequence, $\O_p \to A_1 \to \cdots \to A_n$.  By nondegeneracy of the Serre pairing (as mentioned above) we may assume $A_n = \O_p[\dim V]$.  Consider the total composition, $f: \O_p \to \O_p[\dim V]$.  The pushforward $i_*f$ is a nonzero element of the top exterior power of $T_pV$ under the isomorphism,
\begin{displaymath}
\text{Hom}_X(\O_p, \O_p[\dim V]) \cong \Lambda^{\dim V} T_pX.
\end{displaymath}

% Now take a collection of smooth subvarieties, $V_1, \ldots , V_s$, intersecting transversally at a point $p \in X$ and whose tangent vectors, $\zeta_j$ form a basis for $\text{Hom}(\O_p, \O_p[1])$.  Denote the inclusion map by $i_j: C_j \to X$.  Now given any generator $G$ of $\der$, by the converse to the ghost lemma, for each $C_j$ we can construct a ghost sequence for $\O_p$ whose length is the level of $\O_p$ with respect to $i_j^*G$.  As noted above, we may assume this ghost sequence terminates at $\O_p[1]$, hence the pushforward of the total composition is nothing more than, ${i_j}_*f_j = \zeta_j$, the tangent vector to $C_j$ in $\text{Hom}_X(\O_p, \O_p[\d])$.  We may then construct a ghost sequence:

Now take a collection of smooth subvarieties, $V_1, \ldots , V_s$, intersecting transversally at a point, $p \in X$. Denote the inclusion maps by $i_j: V_j \to X$.  Let $G$ be a  generator of $\dbcoh{X}$. By the Ghost Lemma~\ref{ghost lemma}, for each $V_j$ we can construct a ghost sequence for $\O_p$ whose length is the level of $\O_p$ with respect to $\mathbf{L}i_j^*G$.  As noted above, we may assume this ghost sequence terminates at $\O_p[\dim V_j]$. Denote the total composition by $f_j: \O_p \to \O_p[\dim V_j]$. The pushforward, ${i_j}_*f_j$, is a nonzero element of $\Lambda^{\dim V_j} T_p V_j \subset \Lambda^{\dim V_j} T_p X$.  We may then construct a ghost sequence:

\[
\O_p \overset{{i_1}_*f_1} \to \O_p[\dim V_1] \to \cdots \to \O_p[n-\dim V_s] \overset{{i_{s}}_*f_{s}[n-\dim V_s]} \to \O_p[n].
\]

Each of the $G$ ghosts, ${i_j}_*f_j$, factors into $\text{Lvl}_{V_j}^{\mathbf{L}{i_j}^*G}(\O_p)$ additional $G$ ghosts.  
Hence we have:
\[
\sum_{j=1}^s  \text{Lvl}^{\mathbf{L}{i_j}^*G}_{V_j}(\O_p) \leq \text{Lvl}^G_X(\O_p). 
\] 
\end{example}

Let us use this example to give a simple proof that the ultimate dimension of $\P^n$ is at least $2n$. 

\begin{proposition}
 $\op{udim } \P^n \geq 2n$.
\end{proposition}

\begin{proof}
 We work by induction. Let $G_n = \mathcal O \oplus \mathcal O_{H_1} \oplus \cdots \oplus \mathcal O_{H_{n-1}} \oplus \mathcal O_p$ where $H_i$ is a linear subspace of $\P^n$ of codimension $i$. The induction hypothesis is that, for any point, $q \in \P^n$, not lying in any $H_i$, the level of $\mathcal O_q$ is at least $2n$. Let us tackle the case of $\P^1$ first. Let $q$ be a point distinct from $p$. The sequence
\begin{displaymath}
 \mathcal O_q \to \mathcal O(-1)[1] \to \mathcal O_q[1]
\end{displaymath}
is a ghost sequence for $\mathcal O \oplus \mathcal O_p$. Hence, $\mathcal O_q \notin \langle \mathcal O \oplus \mathcal O_p \rangle_1$ implying that,
\begin{displaymath}
\text{Lvl}_{\mathcal O \oplus \mathcal O_p}(\mathcal O_q) \geq 2.
\end{displaymath}

Now assume we know the result for $\P^j$ when $j \leq n-1$, and let us work on the case $j=n$. Take any point, $q$, not lying on each $H_i$ so that $G_n$ is free near $q$. Take a hyperplane, $H$, passing through $q$ and intersecting each $H_i$ transversally and a line, $L$, passing through $q$ and intersecting each $H_i$ and $H$ transversally. Restricting $G_n$ to $H$ gives an element of $\langle G_{n-1}\rangle_0$ and restricting to $L$ gives an element of $\langle G_1\rangle_0$. By Example~\ref{eg:cyclesandlevels}, the level of $\mathcal O_q$ is at least $2n$. 
\end{proof}

\begin{remark}
 A more careful analysis reveals that
\begin{displaymath}
 \lbrace n, n+1, \ldots, 2n-1, 2n \rbrace \subset \op{OSpec } \P^n.
\end{displaymath}
 We suspect this is in fact an equality. However, this is only known in the case $n=1$.
\end{remark}

%\sidenote{I like the idea to do a computation based on this example. Do you think that the generation time that inequality is an equality? -- We can get $\{2,3,4\}$ as a subset of $\text{OSpec}(\dbcoh{\mathbb P^2})$}

%\begin{example}
%Let us apply the above example to study the projective space $\P^n$.  For $0 \leq r leq n$ let $D_r$ be a smooth section of $\O(r)$.  For each $s$, %consider the generator of $\emph{D}^\emph{b}(\P^n)$:
%\[
%G_s := \O \oplus \cdots \oplus O(s) \oplus O_{D_{s+1}}(s+1) \oplus \cdots \oplus O_{D_{n}}(n).
%\] 
%Let $p$ be a point in the complement of $D_r$ for all $r$ and consider a collection of $n$ lines $L_1, \ldots, L_n$ all intersecting at the point $p$ and %whose tangent directions span $Ext^1(\O_p, \O_p)$.
%\end{example}

When $\mathcal T$ is Ext-finite, we have a (weakly) universal $G$ ghost from any object, $A \in \mathcal T$: we take the cone over the natural evaluation map
\begin{displaymath}
 \bigoplus_{i \in \Z} \op{Hom}_{\mathcal T}(G[-i],A) \otimes_k G[-i] \overset{\op{ev}_A}{\to} A.
\end{displaymath}
Denote, for the moment, the cone by $L_G(A)$. For a general $\mathcal T$ and $G$, the assignment, $A \to L_G(A)$, cannot necessarily be promoted to an endofunctor. To guarantee good behavior of $L_G$, we can assume that $\mathcal T$ is the homotopy category of a triangulated $A_{\infty}$-category, $\mathcal A$, see Chapter 1, Section 6 of \cite{SeiBook}. In this case, we have a cone construction on $\mathcal A$ which enhances the assignment, $L_G$, and guarantees functoriality. We record the definition of $L_G$ and $R_G$ for further use.

\begin{definition}\label{def:left twist}
Let $\mathcal T$ be an Ext-finite triangulated category that is the homotopy category of a triangulated $A_{\infty}$-category, $\mathcal A$. For any pairs of objects, $G$ and $A$, of $\mathcal A$, we have a natural evaluation map
\begin{displaymath}
 \op{Hom}_{\mathcal A}(G,A) \otimes_k G \overset{\op{ev}_A}{\to} A.
\end{displaymath}
Define $L_G: \mathcal A \to \mathcal A$ as the $A_{\infty}$-endofunctor which takes $A$ to the cone over $\op{ev}_A$. We also use the notation, $L_G: \mathcal T \to \mathcal T$, for the induced exact functor on $\mathcal T$, called the \textbf{left twist} by $G$. There is an natural transformation, $\lambda: \op{Id}_{\mathcal A} \to L_G$, which descends to a natural transformation, $\lambda: \op{Id}_{\mathcal T} \to L_G$. We have an exact triangle in $\mathcal T$, where the slashed arrow denotes a degree one morphism:
\begin{center}
\begin{tikzpicture}[description/.style={fill=white,inner sep=2pt}]
\matrix (m) [matrix of math nodes, row sep=3em, column sep=-0.75em, text height=1.5ex, text depth=0.25ex]
{  A & & L_G(A) \\
   & \bigoplus_{i \in \Z} \op{Hom}_{\mathcal T}(G, A[i])\otimes_k G[-i] & \\ };
\path[->,font=\scriptsize]
(m-1-1) edge node[auto] {$\lambda(A)$} (m-1-3)
(m-1-3) edge node[sloped] {$ | $} (m-2-2)	
(m-2-2) edge node[auto] {$ \op{ev}_A $} (m-1-1);
\end{tikzpicture}
\end{center}
\end{definition}

%There is a corresponding notion of a right twist, or adjoint twist in \cite{SeiBook}.

\begin{definition}\label{def:right twist}
Let $\mathcal T$ be an Ext-finite triangulated category that is the homotopy category of a triangulated $A_{\infty}$-category, $\mathcal A$. For any pairs of objects, $G$ and $A$, of $\mathcal A$, we have a natural co-evaluation map
\begin{displaymath}
 A \overset{\op{coev}_A}{\to} \op{Hom}_{\mathcal A}(A,G)^{\vee} \otimes_k G .
\end{displaymath}
Define $R_G: \mathcal A \to \mathcal A$ as the $A_{\infty}$-endofunctor which takes $A$ to the cone over $\op{coev}_A[-1]$. We also use the notation $R_G: \mathcal T \to \mathcal T$ for the induced exact functor on $\mathcal T$, called the \textbf{right twist} by $G$. There is an natural transformation, $\rho: R_G \to \op{Id}_{\mathcal A}$, which descends to a natural transformation, $\rho: R_G \to \op{Id}_{\mathcal T}$. We have an exact triangle in $\mathcal T$:
\begin{center}
\begin{tikzpicture}[description/.style={fill=white,inner sep=2pt}]
\matrix (m) [matrix of math nodes, row sep=3em, column sep=-0.75em, text height=1.5ex, text depth=0.25ex]
{ R_G(A) && A  \\
   &  \text{Hom}_{\mathcal T}(A, G[i])^\vee \otimes_k G[-i] \\ };
\path[->,font=\scriptsize]
(m-1-1) edge node[auto] {$\rho(A)$} (m-1-3)
(m-1-3) edge node[auto] {$\op{coev}$}(m-2-2)	
(m-2-2) edge node[sloped] {$ | $} (m-1-1);
\end{tikzpicture}
\end{center}
\end{definition}

%\begin{example}[Twist functors and the universal ghost map] 
%The map $\lambda(A): A \to L_G(A)$ is $G$ ghost map. It is weakly universal in the sense that it factors (non-uniquely) through any other ghost map.
%Similarly, $\rho(A): R_G(A) \to A$ is the weakly universal $G$ co-ghost map.

\begin{example}
In \cite{ST}, Seidel and Thomas show that for a spherical object (see Definition~\ref{def:spherical object}) the associated left twist functor is an autoequivalence. For the derived Fukaya category of a symplectic manifold, the left twist functor along a Lagrangian sphere is precisely the autoequivalence given by taking a Dehn twist along this sphere.  Seidel and Thomas also show that certain configurations of spherical objects induce the action of a braid group on the category.  While one twist along a Lagrangian sphere provides a single ghost map, we will see in Section~\ref{sec: Spherical Collections} that words in the braid group induce ghost sequences.
\label{twist ghost}
\end{example}

\begin{example}[Global monodromy of the quintic as a ghost map]
Let $X$ be a quintic hypersurface in $\P^4$, $Y$ be the family of Calabi Yau manifolds
that is mirror to $X$ according to Batyrev's construction. A loop around infinity, in the base $\P^1 \backslash \{0,1,\infty\}$, induces a categorical monodromy, $\tw{1}$. This monodromy is a composition of autoequivalences, $\tw{1} = L_{\O} \circ (-\otimes_{\mathcal O} \O(1))$. 

If we choose a hyperplane section, $H$, we get a natural transformation , $\zeta_H: \op{Id}_{\dbcoh{X}} \to \tw{1}$.  For any object, $A$, the map, $\zeta_H(A)$, is ghost for $\O$ and for the essential image of $\dbcoh{H}$ under inclusion.  If we take any generator, $N$, of $\dbcoh{H}$.  Then $\O \oplus N$ generates $\dbcoh{X}$ and $\zeta_H(A)$ is $\O \oplus N$ ghost (see Examples~\ref{twist ghost} and \ref{divisors and ghosts}).

In Section~\ref{sec:graded}, we will see that $\tw{1}$ is precisely the autoequivalence corresponding to twisting the grading in the associated category of graded singularities.
\label{quintic}
\end{example}

\section{Semi-orthogonal Decompositions, Exceptional Collections and Birational Geometry}

\subsection{Semi-orthogonal Decompositions}

Let $\mathcal T$ be a triangulated category and $\mathcal I$ a full subcategory. Recall that the left orthogonal, $\leftexp{\perp}{\mathcal I}$, is the full subcategory $\mathcal T$ consisting of all objects, $T \in \mathcal T$, with $\op{Hom}_{\mathcal T}(T,I) = 0$ for any $I \in \mathcal I$. The right orthogonal, $\mathcal I^{\perp}$, is defined similarly.

\begin{definition}\label{def:SO}
A \textbf{semi-orthogonal decomposition} of a triangulated category, $\mathcal T$, is a
sequence of full triangulated subcategories, $\mathcal A_1, \dots ,\mathcal A_m$, in $\mathcal T$ such that $\mathcal A_i \subset \mathcal A_j^{\perp}$ for $i<j$ and, for every object $T \in \mathcal T$, there exists a diagram:
\begin{center}
\begin{tikzpicture}[description/.style={fill=white,inner sep=2pt}]
\matrix (m) [matrix of math nodes, row sep=1em, column sep=1.5em, text height=1.5ex, text depth=0.25ex]
{  0 & & T_{m-1} & \cdots & T_2 & & T_1 & & T   \\
   & & & & & & & &  \\
   & A_m & & & & A_2 & & A_1 & \\ };
\path[->,font=\scriptsize]
 (m-1-1) edge (m-1-3) 
 (m-1-3) edge (m-1-4)
 (m-1-4) edge (m-1-5)
 (m-1-5) edge (m-1-7)
 (m-1-7) edge (m-1-9)

 (m-1-9) edge (m-3-8)
 (m-1-7) edge (m-3-6)
 (m-1-3) edge (m-3-2)

 (m-3-8) edge node[sloped] {$ | $} (m-1-7)
 (m-3-6) edge node[sloped] {$ | $} (m-1-5)
 (m-3-2) edge node[sloped] {$ | $} (m-1-1)
;
\end{tikzpicture}
\end{center}
where all triangles are distinguished and $A_k \in \mathcal A_k$. We shall denote a semi-orthogonal decomposition by $\langle \mathcal A_1, \ldots, \mathcal A_m \rangle$.

A case of particular importance is if each $\mathcal A_i$ is equivalent to $\dbmod{k}$ as a triangulated category. Let $A_i$ denote the object in $\mathcal T$ corresponding to $k$ in $\mathcal A_i$. In this case, we call $A_1,\ldots,A_m$ an \textbf{exceptional collection}. If, in addition, $\op{Hom}_{\mathcal T}(A_i,A_j[l]) = 0$ for $l \not = 0$, we say that the exceptional collection, $A_1,\ldots,A_n$, is \textbf{strong}.
\end{definition}

As a warning to the reader.  The notion of exceptional collection which appears here is often called a full exceptional collection in the literature.  The distinction is that our exceptional collections always generate the triangulated category in question.

\begin{remark}
 While not required in the definition, it is easy to see that $T$ uniquely determines the diagram appearing in Definition~\ref{def:SO}. 
\end{remark}

The following lemma is clear from the definition of a semi-orthogonal decomposition:

\begin{lemma}
Suppose $\langle \mathcal A_1, \dots ,\mathcal A_m \rangle$ is a semi-orthogonal decomposition of $\mathcal T$ and, for each $i$, $G_i$ is a strong generator of $\mathcal A_i$.  Then, $\bigoplus_{i=1}^m G_i$ is a strong generator of $\mathcal T$.
\end{lemma}

In this section we will analyze how the generation time behaves when we have generators coming from semi-orthogonal decompositions.

Due to work of Bondal, Kuznetsov, and Orlov, it is widely believed that semi-orthogonal decompositions could play an important role in birational geometry.  We have the following result due to Orlov, see \cite{Orl92}:

\begin{theorem} 
Let $\pi:\hat{X} \to X$ be the blow up of a smooth variety, $X$, along a smooth subvariety, $Y$, of codimension c. Let $E$ denote the exceptional divisor on $\hat{X}$ and $\mathcal O_E(1)$ denote the relative twisting sheaf of $\pi|_E: E \to Y$. Denote the inclusion as $j: E \to \hat{X}$. There is a semi-orthogonal decomposition of $\dbcoh{\hat{X}}$ given by
\begin{displaymath} 
\left\langle \dbcoh{Y}, \ldots, \dbcoh{Y}, \dbcoh{X} \right\rangle.
%\left\langle j_*\left((\pi|_E)^*\dbcoh{Y}\otimes_{\mathcal O}\mathcal O_E(-c+1)\right), \ldots , j_*\left((\pi|_E)^*\dbcoh{Y}\otimes_{\mathcal O}\mathcal O_E(-1)\right), \mathbf{L}\pi^*\dbcoh{X} \right\rangle.
\end{displaymath}
In this decomposition, the category $\dbcoh{Y}$ occurs $c-1$ times under the the following equivalences for $-c+1 \leq l \leq -1$:
\begin{displaymath}  
\dbcoh{Y} \cong j_*\left((\pi|_E)^*\dbcoh{Y}\otimes_{\mathcal O}\mathcal O_E(l)\right),
\end{displaymath}
and the category $\dbcoh{X}$ is equivalent to $\mathbf{L}\pi^*\dbcoh{X}$.
%\begin{displaymath}
%\left\langle j_*\left((\pi|_E)^*\dbcoh{Y}\otimes_{\mathcal O}\mathcal O_E(-c+1)\right), \ldots , j_*\left((\pi|_E)^*\dbcoh{Y}\otimes_{\mathcal O}\mathcal O_E(-1)\right), \mathbf{L}\pi^*\dbcoh{X} \right\rangle.
%\end{displaymath}
% Moreover, $\mathbf{L}\pi^*\dbcoh{X}$ is equivalent to $\dbcoh{X}$ and
%$j_*\left((\pi|_E)^*\dbcoh{Y}\otimes_{\mathcal O}\mathcal O_E(l)\right)$
%is equivalent to $\dbcoh{Y}.$
 \label{blow-up semi-orthogonal}
\end{theorem}

Based on the above theorem, and further work of his own, Kuznetsov has proposed the existence of a categorical analogue to the Clemens-Griffiths component, \cite{Kuz}. Roughly, this is the component of a semi-orthogonal decomposition which is not equivalent to a component of the derived category of a variety of smaller dimension.  In what follows, we hope to suggest that the Orlov spectrum can detect, in some cases, when Kuznetsov's Clemens-Griffiths component is nontrivial.

\begin{definition}
Let $\alpha: \mathcal A \to \mathcal T$ be the inclusion of a full triangulated subcategory of $\mathcal T$. The subcategory, $\mathcal A$, is called \textbf{right admissible} if the inclusion functor, $\alpha$, has a right adjoint $\alpha^!$ and \textbf{left admissible} if it has a left adjoint $\alpha^*$. A full triangulated subcategory is called \textbf{admissible} if it is both right and left admissible.
\end{definition}

The proofs of the following lemmas can be found in \cite{BK}:

\begin{lemma}
Let $\mathcal A$ be a full triangulated subcategory of a triangulated category $\mathcal T$ with Serre functor.  Then the following are equivalent:
\renewcommand{\labelenumi}{\emph{\roman{enumi})}}
\begin{enumerate}
 \item $\mathcal A$ is left admissible
 \item $\mathcal A$ is right admissible
 \item $\mathcal A$ is admissible
\end{enumerate}
\end{lemma}

\begin{lemma}
If $\langle \mathcal A_1, \ldots, \mathcal A_m \rangle$ is a semi-orthogonal decomposition of a triangulated category $\mathcal T$ with Serre functor, then $\mathcal A_i$ is admissible for all $i$.  Furthermore, if $\mathcal T =\langle \mathcal A, \mathcal B \rangle$ is a semi-orthogonal decomposition, then $\mathcal B = \leftexp{\perp}{\mathcal A}$.
 \end{lemma}

Let $\langle \mathcal A_1, \ldots, \mathcal A_m \rangle = \mathcal T$ be a semi-orthogonal decomposition of a triangulated category $\mathcal T$ with Serre functor.  Denote each inclusion functor by $\alpha_i: \mathcal A_i \to \mathcal T$. Let $\lambda_i: \leftexp{\perp}{\mathcal A}_i \to \mathcal T$ denote the inclusion of the left orthogonal and $\rho_i: \mathcal A_i^{\perp} \to \mathcal T$ denote the inclusion of the right orthogonal.  For any $X \in \mathcal T$ we have the following exact triangles,
\begin{equation} \label{decomposition1}
\alpha_i\alpha_i^!X \to X \to \rho_i \rho_i^* X,
\end{equation}
\begin{center}
and
\end{center}
\begin{equation} \label{decomposition2}
\lambda_i \lambda_i^!X \to X \to \alpha_i \alpha_i^* X.
\end{equation}

There is an action of the braid group on $m$ strands on the set of all $m$-term semi-orthogonal decompositions of $\mathcal T$, \cite{BK}. The standard generators are given by either taking right mutations, $\mathbb R_i$, or left mutations, $\mathbb L_i$. Let us recall now the definition,

\begin{displaymath} \mathbb R_i(\mathcal A_{\bullet})_j =
\begin{cases}
\mathcal A_j   & \text{ if } j \neq i-1, i \\
\mathcal A_i   & \text{ if } j = i-1 \\
\leftexp{\perp}\langle \mathcal A_1, \ldots, \mathcal A_{i-2}, \mathcal A_i \rangle \cap \langle \mathcal A_{i+1}, \ldots, \mathcal A_m \rangle^{\perp}  & \text{ if } j = i
 \end{cases}
 \end{displaymath}
 \begin{displaymath}
 \mathbb L_i(\mathcal A_{\bullet})_j =
\begin{cases}
\mathcal A_j   & \text{ if } j \neq i, i+1 \\
\leftexp{\perp}\langle \mathcal A_1, \ldots, \mathcal A_{i-1} \rangle \cap \langle \mathcal A_{i}, \mathcal A_{i+2} \ldots, \mathcal A_m \rangle^{\perp}  & \text{ if } j = i \\
\mathcal A_i   & \text{ if } j = i+1. \\
 \end{cases}
\end{displaymath}
Given a generator, $\mathcal G := G_1 \oplus \cdots \oplus G_m$, with each $G_i \in \mathcal A_i$, we can define new generators,
\begin{displaymath}
\mathbb L_i \mathcal G := G_1\oplus \cdots \oplus G_{i-1} \oplus \rho_i \rho_i^*(G_{i+1}) \oplus G_i \oplus \cdots \oplus G_m,
\end{displaymath}
\begin{center}
and
\end{center}
\begin{displaymath}
\mathbb R_i \mathcal G := G_1\oplus \cdots \oplus G_i \oplus \lambda_i \lambda_i^!( G_{i-1}) \oplus G_{i+1} \oplus \cdots \oplus G_m.
\end{displaymath}
 
Further, let us define: $\mathfrak{L}_i:= \mathbb L_m \cdots \mathbb L_i$ and $\mathfrak{R}_i:= \mathbb R_m \cdots \mathbb R_i$ so that,
\begin{displaymath}
\mathfrak{L}_i \mathcal G  = G_1\oplus \cdots \oplus G_{i-1} \oplus \rho_i \rho_i^*(G_{i+1}) \oplus \cdots \oplus \rho_i \rho_i^*(G_m) \oplus G_i
\end{displaymath}
and
\begin{displaymath}
\mathfrak{R}_i \mathcal G  = G_i \oplus \lambda_i \lambda_i^! (G_1) \oplus \cdots \oplus \lambda_i \lambda_i^!(G_{i-1}) \oplus G_{i+1} \oplus \cdots \oplus G_m.
\end{displaymath}
Finally, set $\mathfrak{L}^D:= \mathfrak L_1 \cdots \mathfrak L_{n-1}$ and $\mathfrak{R}^D:= \mathfrak R_n \cdots \mathfrak R_2$.

\begin{definition}
Given a semi-orthogonal decomposition, $\langle \mathcal A_1, \ldots, \mathcal A_n \rangle$, of a triangulated category, $\mathcal T$, with Serre functor, we define the \textbf{left dual semi-orthogonal decomposition} by,
\begin{displaymath}
\langle \mathcal A_1^{\vee}, \ldots , \mathcal A_n^{\vee} \rangle := \mathfrak L_D \langle \mathcal A_1, \ldots , \mathcal A_n \rangle,
\end{displaymath}
and the \textbf{right dual semi-orthogonal decomposition} by,
\begin{displaymath}
\langle \leftexp{\vee}{\mathcal A}_1, \ldots , \leftexp{\vee}{\mathcal A}_n \rangle: =\mathfrak R_D \langle \mathcal A_1, \ldots, \mathcal A_n \rangle.
\end{displaymath}
\label{def:dualexc}
\end{definition}

The following proposition is clear from the definition of mutation:

\begin{proposition}
 We have the following equalities:
\begin{displaymath}
\mathcal A_i^\vee = \langle \mathcal A_1, \ldots, \mathcal A_i,\mathcal A_{i+1}, \ldots, \mathcal A_n\rangle^\perp
\end{displaymath}
\begin{displaymath}
\leftexp{\vee}{\mathcal A}_i = \leftexp{\perp}\langle \mathcal A_1, \ldots, \mathcal A_i,\mathcal A_{i+1}, \ldots, \mathcal A_n\rangle.
\end{displaymath}
\label{prop:dualorth}
\end{proposition}

\begin{lemma}\label{lem:double dual}
 Let $\mathcal T$ be a triangulated category possessing a Serre functor, $S$, and suppose that $\mathcal T$ has a semi-orthogonal decomposition, $\langle \mathcal A_1,\ldots, \mathcal A_n \rangle$. We have isomorphisms for any $X \in \mathcal A_i$:
\begin{displaymath}
S(X) \cong \mathfrak L_D^2(S_{\mathcal A_i}X) \cong S_{\mathcal A_i^{\vee\vee}}\mathfrak L_D^2(X)
\end{displaymath}
\begin{center}
and
\end{center}
\begin{displaymath}
S^{-1}(X) \cong \mathfrak R_D^2(S^{-1}_{\mathcal A_i}X)\cong S^{-1}_{\leftexp{\vee\vee}{\mathcal A}_i} \mathfrak R_D^2(X).
\end{displaymath} 
\end{lemma}

\begin{proof} 
 Note that the effect of the application of $\mathfrak L_D^2$ is to project $\mathcal A_i$ to $\mathcal A_i^{\perp\perp}$. Similarly, $\mathfrak R_D^2$ is the projection from $\mathcal A_i$ to $\leftexp{\perp\perp}{\mathcal A_i}$. Proposition 3.7 of \cite{BK} states that $\mathfrak L_D^2$ commutes with Serre functors. Similarly, $\mathfrak R_D^2$ commutes with inverses to the Serre functors.
\end{proof}

\begin{definition}
 Let $[a,b]$ denote the integer interval with endpoints $a$ and $b$ in $\Z$.  Despite the usual notation, we do not distinguish between $a \leq b$ and $a \geq b$ i.e. $[a,b] = [b,a]$.   Furthermore, our intervals only contain integers.  Let $I$ be a subset of $\Z$. We say that $I$ has a \textbf{gap} of length $s$ if, for some $a$, $[a,a+s+1] \cap I = \lbrace a, a+s+1 \rbrace$. We say that a triangulated category, $\mathcal T$, has a \textbf{gap} of length $s$ if $\ops{OSpec} \mathcal T$ has a gap of length $s$.
\end{definition}

\begin{theorem} 
Suppose $\langle \mathcal A_1 , \ldots , \mathcal A_n \rangle$ is a semi-orthogonal decomposition of $\mathcal T$ and $\mathcal G:= G_1 \oplus \cdots \oplus G_n$ is a generator of $\mathcal T$ with $G_i \in \mathcal A_i$.  Let $M := \op{max}_i \{ \tritime_{\mathcal A_i}(G_i)\}$. Any gap inside $[\tritime_{\mathcal T}(\mathcal G),\tritime_{\mathcal T}(\mathfrak L_D(\mathcal G))] \cap \ops{OSpec}\mathcal T$ has length at most $M$.  In particular, if $\tritime_{\mathcal A_i}(G_i)$ equals the Rouquier dimension of $\mathcal A_i$ for each $i$, then any gap inside $[\tritime_{\mathcal T}(\mathcal G),\tritime_{\mathcal T}(\mathfrak L_D (\mathcal G))] \cap \ops{OSpec}\mathcal T$ has length at most $\max_i \op{rdim} \mathcal A_i$.  The same statement is true passing to the right dual.
\label{semi-orthogonal gaps}
\end{theorem}

\begin{proof}
Let $\langle \mathcal C_1, \ldots, \mathcal C_n \rangle$ be a semi-orthogonal decomposition and let $\mathcal H = H_1 \oplus \cdots \oplus H_n$ be a generator with $H_i \in \mathcal C_i$. Since $\rho_i^*: \mathcal T \to \mathcal C_i^{\perp}$ is essentially surjective we have, \begin{displaymath}
\tritime_{\mathcal C_i^{\perp}}(\rho_i^* \mathcal H) \leq \tritime_{\mathcal T}(\mathcal H).
\end{displaymath}
Now by definition,
\begin{displaymath}
\mathfrak{L}_i \mathcal H  = H_1\oplus  \cdots \oplus H_{i-1} \oplus \rho_i \rho_i^*(H_{i+1}) \oplus \cdots \oplus \rho_i \rho_i^*(H_n) \oplus H_i
\end{displaymath}
\begin{displaymath}
\text{and}
\end{displaymath}
\begin{displaymath}
\rho_i^*(\mathcal H)  = H_1 \oplus \cdots \oplus H_{i-1} \oplus \rho_i^*(H_{i+1}) \oplus \cdots \oplus \rho_i^*(H_n).
\end{displaymath}
Hence $\mathfrak{L}_i \mathcal H = \rho_i \rho_i^*(\mathcal H) \oplus H_i$.  Therefore, $\mathfrak{L}_i \mathcal H$ generates the left orthogonal of $\mathcal C_i$ in at most $\tritime_{\mathcal T}(\mathcal H)$-steps.  Furthermore, as $H_i$ is a summand of $\mathfrak{L}_i \mathcal H$, $\mathfrak{L}_i \mathcal H$ generates $\mathcal C_i$ in at most $\tritime_{\mathcal C_i}(H_i)$-steps.  Now triangle~\eqref{decomposition1} tells us that $\tritime_{\mathcal T}(\mathfrak{L}_i \mathcal H) \leq \tritime_{\mathcal T}(\mathcal H) + \tritime_{\mathcal C_i}(H_i)+1 \leq \tritime(\mathcal H) +M +1$. We have learned that the generation time increases in increments of at most $M+1$ after application of a single $\mathfrak{L}_i$. A similar argument shows that, after applying the mutation $\mathfrak{R}_i$, the generation time does not increase by more than $M+1$. If we apply $\mathfrak R_D^2 \mathfrak L_D^2$ to $\mathcal G$, we return to $\mathcal{G}$ by Lemma~\ref{lem:double dual}.

Since the generation time must return to its original value and can only increase in increments of at most $M+1$, any gap within the interval with endpoints $\tritime_{\mathcal T}(\mathcal G)$ and $\tritime_{\mathcal T}(\mathfrak{L}_D(\mathcal G))$ has length at most $M$.  The proof for the right dual statement is the same. 
\end{proof}
% Now, as the Serre functor, is an autoequivalence of $\mathcal T$, one has,
% \begin{displaymath}
%  \tritime(\mathfrak{L}_n  \cdots  \mathfrak{L}_1 \mathfrak{L}_n  \cdots  \mathfrak{L}_1 (\mathcal G)) = \tritime(S^{-1} \circ \mathfrak{L}_n  \cdots  \mathfrak{L}_1 \mathfrak{L}_n  \cdots  \mathfrak{L}_1 (\mathcal G)) = \tritime(S_{\mathcal A_1}(G_1) \oplus \cdots \oplus S_{\mathcal A_n}(G_n))
% \end{displaymath}
% by Lemma~\ref{lem:double dual}. Given a semi-orthogonal decomposition $\langle \mathcal C_1,\ldots,\mathcal C_n \rangle$ and a generator $H = H_1 \oplus \cdots \oplus H_n$ with $H_i \in \mathcal C_i$, we have naive bounds on the generation time
% \begin{displaymath}
%  \max\lbrace \tritime_{\mathcal C_1}(H_1), \tritime_{\langle \mathcal C_2,\ldots,C_n \rangle}(H_2 \oplus \cdots \oplus H_n) \rbrace \leq \tritime(H_1 \oplus \cdots \oplus H_n) \leq \tritime_{\mathcal C_1}(H_1) + \tritime_{\langle \mathcal C_2,\ldots,C_n \rangle}(H_2 \oplus \cdots \oplus H_n) + 1
% \end{displaymath}
% So
% \begin{displaymath}
%  |\tritime(G_1 \oplus \cdots \oplus G_n) - \tritime(S_{\mathcal A_1}(G_1) \oplus G_2 \oplus \cdots \oplus G_n)| \leq \tritime(G_1) +1 \leq M +1.
% \end{displaymath}
% Iterating, we have
%  $|\tritime(G_1 \oplus G_2 \oplus \cdots \oplus S_{\mathcal A_i}(G_i) \oplus \cdots \oplus S_{\mathcal A_n}(G_n))-\tritime($
% We see that after applying $\mathfrak{L}_n \cdots \mathfrak{L}_1 \mathfrak{L}_n \cdots \mathfrak{L}_1$ the generation time must return to within $M+1$ of $\tritime(\mathcal G)$.

\subsection{A Conjectural Aside} \label{conjectural}

The ``results'' in this subsection are all purely conjectural. However, nothing from this section will be used for further argument.

Recall from the introduction that the following conjecture appears in \cite{O4}, where it is proven for curves.

\begin{conjecture}%[Orlov] 
For a smooth algebraic variety, $X$, the Krull dimension of $X$ and the Rouquier dimension of $\dbcoh{X}$ are equal.
\label{dimension conjecture}
\end{conjecture}

Now, in light of Theorem~\ref{semi-orthogonal gaps}, let us propose our own conjecture.

\begin{conjecture} 
Let $X$ be a smooth algebraic variety and $\dbcoh{X} =\langle \mathcal A_1 , . . . , \mathcal A_n \rangle$ be a semi-orthogonal decomposition.   The length of any gap in $\dbcoh{X}$ is at most the maximal Rouquier dimension amongst the $\mathcal A_i$.
\label{semi-orthogonal gap conjecture}
\end{conjecture}

\begin{corollary} 
Suppose Conjectures~\ref{dimension conjecture} and~\ref{semi-orthogonal gap conjecture} hold.  If $X$ is a smooth variety then any gap of $\dbcoh{X}$ has length at most the Krull dimension of $X$.
\label{gap bound}
\end{corollary}

Let us propose another conjecture:
\begin{conjecture}  
Let $X$ be a smooth algebraic variety. If $\mathcal A$ is an admissible subcategory of $\dbcoh{X}$, then the length of any gap of $\mathcal A$ is at most the maximal length of any gap of $\dbcoh{X}$. Conversely, if $\mathcal A$ has a gap of length at least $s$, then so does $\dbcoh{X}$.
\label{admissible gap conjecture}
\end{conjecture}

\begin{corollary} 
Suppose Conjectures~\ref{dimension conjecture}, \ref{semi-orthogonal gap conjecture}, and \ref{admissible gap conjecture} hold.  Let $X$ and $Y$ be birational smooth proper varieties of dimension $n$.  The category, $\dbcoh{X}$, has a gap of length $n$ or $n-1$ if and only if $\dbcoh{Y}$ has a gap of the same length i.e. the gaps of length greater than $n-2$ are a birational invariant.
\label{gap birational invariant}
\end{corollary}
\begin{proof}
We may suppose that $Y$ is the blow-up of $X$ along $Z$.  By Theorem~\ref{blow-up semi-orthogonal}, we have a semi-orthogonal decomposition,
\begin{displaymath}
\dbcoh{X} = \langle \dbcoh{X} , \dbcoh{Z}, \ldots, \dbcoh{Z} \rangle.
\end{displaymath}
By Conjecture~\ref{admissible gap conjecture}, if $\dbcoh{X}$ has a gap of length at least $s$, then so does $\dbcoh{Y}$. Now suppose $\dbcoh{Y}$ has a gap of length $s > n-2$.  By Corollary~\ref{gap bound}, the length of any gap in $\text{D}^{\text{b}}(Z)$ is at most $n-2$.  Thus, by Conjecture~\ref{semi-orthogonal gap conjecture}, $\dbcoh{X}$ must have gap $s$ as well.        
\end{proof}

\begin{corollary}
Suppose Conjectures~\ref{dimension conjecture}, \ref{semi-orthogonal gap conjecture}, and \ref{admissible gap conjecture} hold.  If $X$ is a rational variety of dimension $n$, then any gap in $\dbcoh{X}$ has length at most $n-2$.
\end{corollary}
\begin{proof}
It is well known that $\mathbb P^n$ has an exceptional collection. In particular, it has a semi-orthogonal decomposition into categories of Rouquier dimension zero.  By Conjecture~\ref{semi-orthogonal gap conjecture}, $\dbcoh{\mathbb{P}^n}$ has no gaps.  The statement follows from Corollary~\ref{gap birational invariant}.
\end{proof}

In Section~\ref{sec:ungraded} we will see that the category of singularities of an $A_n$-singularity in even dimension has gaps.  In Section~\ref{sec:graded}, we will explore semi-orthogonal decompositions for hypersurfaces in $\mathbb P^n$ and their Orlov spectra.

\subsection{Bounds on generation time for exceptional collections} \label{exceptional collections} 

The following proposition is an immediate consequence of the main theorem of \cite{BF}.

\begin{proposition}
 Let $A_1,\ldots,A_n$ be a strong exceptional collection in an Ext-finite triangulated category, $\mathcal T$, that possesses an enhancement. The generation time of $G_A = A_1 \oplus \cdots \oplus A_n$ is bounded above by
\begin{displaymath}
 \max \ \lbrace i \ | \ \op{Hom}_{\mathcal T}(G_A,S^{-1}(G_A)[i]) \not = 0 \rbrace.
\end{displaymath}
\end{proposition}

In this subsection, we establish a new bound for a general exceptional collection. We require the machinery of triangulated $A_{\infty}$-categories. We will recall the bare necessities and refer the reader to \cite{SeiBook} for a deeper discussion. We also follow the (slightly nonstandard) sign, ordering, and notational conventions found in loc. cit.

Recall that for an $A_{\infty}$-category, $\mathcal A$, the morphism spaces are graded vector spaces and we have multi-compositions. For any sequence of objects, $X_0,\ldots,X_n$, $n > 0$, there is a $k$-linear map
\begin{displaymath}
 m_n: \op{Hom}_{\mathcal A}(X_0,X_1) \otimes_k \cdots \otimes_k \op{Hom}_{\mathcal A}(X_{n-1},X_n) \to \op{Hom}_{\mathcal A}(X_0,X_n)
\end{displaymath}
of degree $2-n$. The ordering of the morphism spaces is as in loc. cit.  These maps satisfy a hierarchy of quadratic relations. The first two of which state that $m_1$ is a differential on each $\op{Hom}_{\mathcal A}(X_0,X_1)$ and $m_2$ is map of complexes. $\mathcal A$ is called \textbf{minimal} if $m_1=0$. 

The homotopy category, $H(\mathcal A)$, of $\mathcal A$ is defined by taking the same objects as $\mathcal A$ but taking morphisms between $X_0$ and $X_1$ to be $H^0(\op{Hom}_{\mathcal A}(X_0,X_1),m_1)$. We also have the graded category where we take the same objects but we take morphisms to be $H^*(\op{Hom}_{\mathcal A}(X_0,X_1),m_1)$. This is denoted by $H^*(\mathcal A)$. If $H^*(\mathcal A)$ has finite dimensional morphisms spaces, i.e. if one has $\dim_k \op{Hom}_{H^*(\mathcal A)}(X,Y) < \infty$ for any pair of objects, $X,Y \in \mathcal A$, then $\mathcal A $ is called \textbf{cohomologically-finite}.

We shall always assume that $\mathcal A$ is strictly unital meaning, for each $A \in \mathcal A$, there is an element, $\op{id}_A \in \op{Hom}_{\mathcal A}(A,A)$, that passes to the identity on $H(\mathcal A)$ and satisfies the following: for any $\phi: B \to A$ and $\psi:A \to B$, we have $m_2(\phi,\op{id}_A) = \phi$,  $m_2(\op{id}_A,\psi) = \psi$ and any multi-composition $m_n(\phi_1,\otimes\cdots\otimes\op{id}_A\otimes\cdots\otimes \phi_{n-1}) =0$ for $n \geq 3$.

A right module over $\mathcal A$ is an $A_{\infty}$-functor from $\mathcal A^{\op{op}}$ to the dg-category of chain complexes of $k$-modules. 
Right modules over $\mathcal A$ form an $A_{\infty}$-category. An $A_\infty$ category, $\mathcal A$, is called \textbf{triangulated}, or often \textbf{pretriangulated} \cite{BK2}, if its essential image, under the Yoneda embedding, in $H(\op{Mod-}\mathcal A)$ is a triangulated category.

Given a generator, $G$, of $H(\mathcal A)$, twisted complexes concretely express how any object in $\mathcal A$ is built from $G$ using cones, shifts, and summands.  They are a useful tool in analyzing generation time.  We recall their definition now, so that may be used in what follows.

First, we additively enlarge to create a new $A_{\infty}$-category. Let $\mathcal B$ be an $A_{\infty}$-category. Its additive enlargement is the $A_{\infty}$-category, $\Sigma \mathcal B$, whose are objects are denoted by
\begin{displaymath}
 \bigoplus_{i \in I} V_i \otimes_k Y_i,
\end{displaymath}
with $I$ a finite set, $V_i$ finite-dimensional graded vector spaces, and $Y_i$ objects of $\mathcal B$. The morphism space in $\Sigma \mathcal B$ between $C := \bigoplus_i V_i \otimes_k Y_i$ and $D := \bigoplus_i W_i \otimes_k Y_i$ is
\begin{displaymath}
 \op{Hom}_{\op{Tw-}\mathcal B}(C,D) := \bigoplus_{i,j} \op{Hom}_k(V_i, W_j) \otimes_k \op{Hom}_{\mathcal B}(Y_i,Y_j)
\end{displaymath}
with the natural associated grading. The multi-compositions in $\Sigma \mathcal B$ are natural linear extensions of those in $\mathcal B$.

A twisted complex over $\mathcal B$ is a pair, $(C,\delta_C)$, where $C$ is an object of $\Sigma \mathcal B$ and where $\delta_C$ is an endomorphism of $C$ in $\Sigma \mathcal B$ of degree one. We require that $\delta_C$ satisfies the following conditions: one, there is a finite decreasing filtration of the $V_i$'s that is preserved under the action of $\delta_C$ and so that the map induced by $\delta_C$ on the associated graded pieces is zero, and, two, the sum
\begin{equation}\label{eqn:MC}
 \sum_{i=1}^{\infty} m_r(\delta_C^{\otimes r}) = 0
\end{equation}
where $m_r$ is the $r$-th composition in $\Sigma B$. Note that finiteness of the sum in Equation \ref{eqn:MC} is a consequence of the first condition on $\delta_C$. We will often suppress the $\delta_C$ from the notation of a twisted complex. Such a twisted complex was called a one-sided twisted complex in \cite{BK2}.

Twisted complexes over $\mathcal B$ form an $A_{\infty}$-category, denoted by $\op{Tw-}\mathcal B$. The graded vector space of morphisms between two twisted complexes $(C,\delta_C)$ and $(D,\delta_D)$ with $C = \bigoplus_i V_i \otimes_k Y_i$ and $D = \bigoplus_i W_i \otimes_k Y_i$ is
\begin{displaymath}
 \op{Hom}_{\op{Tw-}\mathcal B}(C,D) := \bigoplus_{i,j} \op{Hom}_k(V_i, W_j) \otimes_k \op{Hom}_{\mathcal B}(Y_i,Y_j)
\end{displaymath}
with the natural associated grading. 

If we have $n$ twisted complexes, $(C_i,\delta_{C_i})$, $0 \leq i \leq n$, then the $n$-order multi-composition on $\op{Tw-}\mathcal B$ is given by
\begin{equation}
 \phi_1 \otimes \cdots \otimes \phi_n \mapsto \sum_{i_0,\ldots,i_n \geq 0} m_{n+i_0+\cdots+i_n}(\delta_{C_0}^{\otimes i_0} \otimes \phi_1 \otimes \delta_{C_1}^{\otimes i_1} \otimes \cdots \otimes \delta_{C_{n-1}}^{\otimes i_{n-1}} \otimes \phi_n \otimes \delta_{C_n}^{\otimes i_n}).
\label{eqn:m_ntwisted}
\end{equation}
The multi-compositions in $\op{Tw-}\mathcal B$ satisfy the $A_{\infty}$-relations as a result of Equation \ref{eqn:MC}.

We say that $A_1,\ldots,A_n$ is an \textbf{exceptional collection} in $\mathcal A$ if $A_1,\ldots,A_n$ is an exceptional collection in $H(\mathcal A)$. Similarly, $A_1,\ldots,A_n$ is strong in $\mathcal A$ is strong in $H(\mathcal A)$. We will say that $A_1,\ldots,A_n$ is minimal when the $A_{\infty}$-endomorphism algebra of the $A_i$'s is minimal. When $\mathcal A$ has an exceptional collection, we can provide a normalized form for objects of $\mathcal A$.

\begin{definition}
 Let $A$ denote the full subcategory of $\mathcal A$ consisting of $A_1,\ldots,A_n$ and $\op{Tw-}A$ denote the category of twisted complexes over $A$. Let $(C,\delta_C)$ be a twisted complex over $A$ and let
\begin{displaymath}
 C = \bigoplus_{i=0}^n V_i \otimes_k A_i.
\end{displaymath}
 Consider the filtration $F^lC = \bigoplus_{i = l}^n V_i \otimes_k A_i$. We say that $(C,\delta_C)$ is \textbf{normalized} if $\delta_C$ respects the filtration and vanishes on the associated graded pieces, $F^lC/F^{l+1}C$.
\end{definition}

\begin{lemma}
 Let $\mathcal A$ be a cohomologically-finite triangulated $A_{\infty}$-category with $A_1,\ldots, A_n$ an exceptional collection. Every object of $\mathcal A$ is isomorphic to a normalized twisted complex over $A$ in $H(\mathcal A)$.
\label{lem:restrictedtwistedcomplexes}
\end{lemma}

\begin{proof}
 This is essentially Lemma 5.13 of \cite{SeiBook}. For any object, $Y$ of $\mathcal A$, we set $Y_n = Y$ and $Y_{i-1} = L_{A_i} Y_i$. As $Y_0$ lies in the left orthogonal to each of the $A_i$ in $H(\mathcal A)$, it follows that it lies in the left orthogonal to the category generated by $A_1 \oplus \cdots \oplus A_n$, which by assumption is all of $H(\mathcal A)$.  Hence, $Y_{0}$ is acyclic. Choose a basis for $\op{Hom}_{H^*(\mathcal A)}(A_i,Y_i)$ and lift to a subspace, $V_i$, of cycles in $\op{Hom}_{\mathcal A}(A_i,Y_i)$. This lift provides a splitting $V_i \hookrightarrow \op{Hom}_{\mathcal A}(A_i,Y_i) \twoheadrightarrow V_i$. Now, we work backwards to get a normalized twisted complex quasi-isomorphic to $Y$. Since $Y_0$ is trivial in $H(\mathcal A)$, $Y_1$ is quasi-isomorphic to $V_1 \otimes_k A_1$. Now, $Y_2$ is quasi-isomorphic to the cone over the composition of morphisms,
\begin{displaymath}
 V_1 \otimes_k A_1 \to \op{Hom}_{\mathcal A}(A_1,Y_1) \otimes_k A_1 \to Y_1 \to \op{Hom}_{\mathcal A}(A_2,Y_2) \otimes_k A_2[1] \to V_2 \otimes_k A_1[1],
\end{displaymath}
which we denote by $X_2$. As a cone, $X_2$ is a normalized twisted complex with
\begin{displaymath}
 X_2 = V_1 \otimes_k A_1[1] \oplus V_2 \otimes_k A_2[1].
\end{displaymath}
 Applying induction, we see that $X_i$ is the cone over a map from a normalized twisted complex, $X_{i-1}$ of the form
\begin{displaymath}
 X_{i-1} = \bigoplus_{l=1}^{i-1} V_l \otimes_k A_l[i-1],
\end{displaymath}
 to $V_i \otimes_k A_i[i]$. Thus, $X_i$ is a normalized twisted complex quasi-isomorphic to $Y_i$. Setting $C = X_n$ gives the desired twisted complex.
\end{proof}

\begin{remark}
 As noted in \cite{SeiBook}, the $Y_i$ constructed in Lemma \ref{lem:restrictedtwistedcomplexes} fit into a Postnikov tower:
\begin{center}
\begin{tikzpicture}[description/.style={fill=white,inner sep=2pt}]
\matrix (m) [matrix of math nodes, row sep=0.5em, column sep=1em, text height=1.5ex, text depth=0.25ex]
{  Y_n & & Y_{n-1} & \cdots & Y_2 & & Y_1 & & Y_0   \\
   & & & & & & & &  \\
   & V_n \otimes_k A_n & & & & V_2 \otimes_k A_2 & & V_1 \otimes_k A_1 & \\ };
\path[->,font=\scriptsize]
 (m-1-1) edge (m-1-3) 
 (m-1-3) edge (m-1-4)
 (m-1-4) edge (m-1-5)
 (m-1-5) edge (m-1-7)
 (m-1-7) edge (m-1-9)

 (m-3-8) edge (m-1-7)
 (m-3-6) edge (m-1-5)
 (m-3-2) edge (m-1-1)

 (m-1-9) edge node[sloped] {$ | $} (m-3-8)
 (m-1-7) edge node[sloped] {$ | $} (m-3-6)
 (m-1-3) edge node[sloped] {$ | $} (m-3-2)
;
\end{tikzpicture}
\end{center}
 One can also prove Lemma \ref{lem:restrictedtwistedcomplexes} by realizing the diagonal bi-module as a normalized twisted complex over the category of bi-modules consisting of $A_i \boxtimes B_j$ and then convolving. See Proposition $3.8$ of \cite{Kuz09} for a particular example. We thank Kuznetsov for pointing this out.
\end{remark}

As in the case of a triangulated category, there is a left dual collection to $A_1,\ldots,A_n$ in $\mathcal A$. We set
\begin{displaymath}
 B_{n+1-k} := L_{A_1}L_{A_2}\cdots L_{A_{k-1}}(A_k).
\end{displaymath}
It is straightforward to check that $B_1,\ldots,B_n$ descends to the left dual collection to $A_1,\ldots,A_n$ in $H(\mathcal A)$ as defined in Definition~\ref{def:dualexc}.

\begin{lemma}
 Let $\phi: X \to Y$ be a morphism in $H(\mathcal A)$.  Denote the following induced morphisms by:
 \begin{displaymath}
 \phi_i^t :=\op{Hom}_{H(\mathcal A)}(A_i,L_{A_{i+1}}\cdots L_{A_{n}}(\phi)[t]).
 \end{displaymath}
The morphism, $\phi$, is co-ghost for $G_B = \bigoplus_{i=1}^n B_i$ if and only if $\phi_i^t$ vanishes for $1 \leq i \leq n$ and any $t \in \Z$.
\label{lem:dualco-ghost}
\end{lemma}

\begin{proof}
 Take any $B_{n+1-i}$. From Proposition~\ref{prop:dualorth}, $\op{Hom}_{H(\mathcal A)}(A_l,B_{n+1-i}[t])$ is zero for $l \not = i$ for any $t$. Note that, because of this orthogonality,
 \begin{displaymath}
 \op{Hom}_{H(\mathcal A)}(X,B_{n+1-i}[t]) \cong \op{Hom}_{H(\mathcal A)}(X_i,B_{n+1-i}[t]),
\end{displaymath} 
 where $X_i = L_{A_{i+1}} \cdots L_{A_{n}}(X)$. Similarly, the evaluation map $\bigoplus_j \op{Hom}_{H(\mathcal A)}(A_i[j],X_i) \otimes_k A_i[j] \to X_i$ induces an isomorphism,
\begin{gather*}
 \op{Hom}_{H(\mathcal A)}(X_i,B_{n+1-i}[t]) \cong \op{Hom}_{H(\mathcal A)}(\bigoplus_j \op{Hom}_{H(\mathcal A)}(A_i[j],X_i) \otimes_k A_i[j],B_{n+1-i}[t]) \\ \cong \left(\op{Hom}_{H(\mathcal A)}(A_i[t],X_i)\right)^\vee.
\end{gather*}
 The same statement is true for $Y$. We see that the map,
\begin{displaymath}
 \op{Hom}_{H(\mathcal A)}(\phi,B_{n+1-i}[t]): \op{Hom}_{H(\mathcal A)}(Y,B_{n+1-i}[t]) \to \op{Hom}_{H(\mathcal A)}(X,B_{n+1-i}[t]),
\end{displaymath}
 coincides with the map,
\begin{displaymath}
 \left(\phi^{-t}_i\right)^\vee: \op{Hom}_{H(\mathcal A)}(A_i[t],L_{A_n}\cdots L_{A_{i+1}}(Y))^\vee \to \op{Hom}_{H^*(\mathcal A)}(A_i[t],L_{A_n}\cdots L_{A_{i+1}}(X))^\vee,
\end{displaymath}
 under the isomorphisms above. This implies the claim.
\end{proof}

We have the following corollary:

\begin{corollary}
 Assume we have a minimal exceptional collection $A_1,\ldots,A_n$ in $\mathcal A$. Let $X = \bigoplus V_i \otimes_k A_i$ and $Y = \bigoplus W_i \otimes_k A_i$ be twisted complexes.  A cocycle, $\phi \in \op{Hom}_{\op{Tw-}A}(X,Y)$, is co-ghost for $B_{n+1-i}$ if and only if the component $\phi^{ii}: V_i \otimes_k A_i \to W_i \otimes_k A_i$ is zero in $H(\mathcal A)$.
\label{cor:Bco-ghost}
\end{corollary}

%\sidenote{I have no idea what Kuznetsov is trying to say here or the diagram he makes at the bottom of img004.jpg.  I think he doesn't like thinking of Hom as a functor and applying it to a morphism.  He writes that this is bad notation, and did the same thing for the definition of ghosts which I changed.}

\begin{proof}
 Note that, by minimality, $\phi^{ii}$ must be some matrix in $\op{Hom}_k(V_i,W_i)$ tensored with the identity on $A_i$. In particular, it is a cocycle.

 Let $\phi: X \to Y$ be a map of normalized twisted complexes over $A$.  We say that $X$ has length $l$ if $X = \bigoplus_{i=1}^l V_i \otimes_k A_i$. We proceed by induction on the length of the twisted complexes $X$ and $Y$. The case $n=1$ is clear.

 Let us assume we know the claim is true when the lengths of $X$ and $Y$ are less than $n$ and assume we have an exceptional collection of length $n$. For notation, let $X = \bigoplus_{i=1}^n V_i \otimes_k A_i$ and $Y = \bigoplus_{i=1}^n W_i \otimes_k A_i$. Note that the inclusion $V_n \otimes_k A_n \hookrightarrow X$ is a cocycle in $\operatorname{Hom}_{\mathcal A}(V_n \otimes_k A_n,X)$. Let $X_{n-1} = \bigoplus_{i=1}^{n-1} V_i \otimes_k A_i$ with $\delta^{ij}_{X_{n-1}} = \delta^{ij}_X$ for $0 \leq i,j \leq n-1$. The cone over $V_n \otimes_k A_n \hookrightarrow X$ is the twisted complex $X \oplus V_n \otimes_k A_n[1]$ with twisting cochain $\begin{pmatrix} \delta_X & \id_{V_n \otimes_k A_n} \\ 0 & 0 \end{pmatrix}$. The projection $X \oplus V_n \otimes_k A_n[1] \to X_{n-1}$ is a cocycle and induces a quasi-isomorphism of $X_{n-1}$ with $L_{A_n}(X)$. The map $\phi: X \to Y$ induces a commutative diagram,
\begin{center}
\begin{tikzpicture}[description/.style={fill=white,inner sep=2pt}]
\matrix (m) [matrix of math nodes, row sep=3em, column sep=3.5em, text height=1.5ex, text depth=0.25ex]
{ X \oplus V_n \otimes_k A_n[1] & Y \oplus W_n \otimes_k A_n[1]\\
  X_{n-1} & Y_{n-1} \\ };
\path[->,font=\scriptsize]
(m-1-1) edge node[auto] {$\begin{pmatrix} \phi & 0 \\ 0 & \phi^{nn}[1] \end{pmatrix}$} (m-1-2)
        edge (m-2-1)
(m-1-2) edge (m-2-2)
(m-2-1) edge node[auto] {$\phi_{n-1}$} (m-2-2);
\end{tikzpicture}
\end{center}
 where $\phi^{ij}_{n-1} = \phi^{ij}$ for $1 \leq i,j \leq n-1$. For $1 \leq i \leq n-1$, $\phi$ is $B_{n+1-i}$ co-ghost if and only if $\phi_{n-1}$ is $B_{n+1-i}$ co-ghost. Also, $\phi^{ii}$ vanishes if and only $\phi_{n-1}^{ii}$ vanishes. So, to verify the claim in the case that $1 \leq i \leq n-1$, we can pass to $\phi_{n-1}:X_{n-1} \to Y_{n-1}$ and apply the induction hypothesis. When $i=n$, we have the commutative diagram
\begin{center}
\begin{tikzpicture}[description/.style={fill=white,inner sep=2pt}]
\matrix (m) [matrix of math nodes, row sep=3em, column sep=3.5em, text height=1.5ex, text depth=0.25ex]
{ V_n \otimes_k A_n & X  \\
  W_n \otimes_k A_n & Y  \\ };
\path[->,font=\scriptsize]
(m-1-1) edge (m-1-2)
        edge node[auto] {$ \phi^{nn} $} (m-2-1)
(m-1-2) edge node[auto] {$ \phi $} (m-2-2)
(m-2-1) edge (m-2-2);
\end{tikzpicture}
\end{center}
 The inclusions induce isomorphisms,
 \begin{displaymath}
 \op{Hom}_{H(\mathcal A)}(A_n,V_n\otimes_k A_n[t]) \cong \op{Hom}_{H(\mathcal A)}(A_n,X[t])
 \end{displaymath}
  and
  \begin{displaymath}
  \op{Hom}_{H(\mathcal A)}(A_n,W_n\otimes_k A_n[t]) \cong \op{Hom}_{H(\mathcal A)}(A_n,Y[t]).
  \end{displaymath}
  Hence, $\op{Hom}_{H(\mathcal A)}(A_n,\phi[t])=0$ if and only if $\op{Hom}_{H(\mathcal A)}(A_n,\phi^{nn}[t])$ is trivial. Precomposing with the identity on $V_n \otimes_k A_n$ shows that $\op{Hom}_{H(\mathcal A)}(A_n,\phi^{nn}[t])$ vanishes for all $t$ if and only if $\phi^{nn}$ vanishes.
\end{proof}

Let $A_1,\ldots,A_n$ be a minimal exceptional. Abusing notation, let $A$ stand for the endomorphism $A_{\infty}$-algebra of the object $\bigoplus_{i = 1}^n A_i$. Let $I$ be the subspace of $A$ consisting of $\phi \in A$ for which $\op{Hom}_{H(\mathcal A)}(\phi,B_i)$ is zero for each $i$. Let us set $I^1 = I$ and define $I^n$ as the following vector space:
%\begin{displaymath}
% I^n = \op{Span}_{t=2}^{\infty}\op{Image}\left(m_t: \bigoplus_{s_1+\cdots+s_t -t \geq n} I^{s_1} \otimes_k \cdots \otimes_k I^{s_t} \to I\right).
%\end{displaymath}
\begin{displaymath}
\left\langle m_t(i_1,\ldots,i_t) \ : \ i_j \in I^{s_j} \text{ with } 1 \leq s_j \leq n-1, s_1+\cdots+s_t-t \geq n-1, \text{ and } t \geq 2 \right\rangle.
\end{displaymath}

\begin{definition}
 We set $\op{LL}_{\infty}(A) := \min \ \lbrace \ n \ | \ I^n = 0 \ \rbrace$.  We call $\op{LL}_{\infty}(A)$ the \textbf{Loewy length} of $A$.
\end{definition}
In the case that $m_i = 0$ for $i \not =2$, $A$ is an algebra and $I^n$ is the standard $n$-th power of $I$ as an ideal of $A$. So, $\op{LL}_{\infty}(A)$ equals the minimal $n$ for which any product of elements of $I$ of length $n$ is zero.

\begin{proposition}
 Let $\mathcal A$ be a cohomologically-finite triangulated $A_{\infty}$-category possessing an exceptional collection $A_1,\ldots,A_n$. The generation time of $G_B$ in $H(\mathcal A)$ is bounded above by $\op{LL}_{\infty}(A')-1$ where $A'$ is a minimal $A_{\infty}$-algebra quasi-isomorphic to $A$.
\label{prop:excgentime}
\end{proposition}

\begin{proof}
 Let $\phi_i:X_{i-1} \to X_i$, for $1 \leq i \leq s$, be a chain of $G_B$ co-ghosts. By Lemma~\ref{lem:restrictedtwistedcomplexes}, we can assume each $\delta_{X_i}$ has components lying in $I$. By Corollary~\ref{cor:Bco-ghost}, the components of $\phi_i$ must lie in $I$. From the formula in Equation~\ref{eqn:m_ntwisted}, we see that all components of $\phi_s \circ \cdots \circ \phi_1$ lie in $I^s$. If $s \geq \op{LL}_{\infty}(A)$, then $\phi_s \circ \cdots \circ \phi_1$ is zero.
\end{proof}

\begin{proposition}\label{prop: lower Loewy bound}
 Let $\mathcal T$ be a triangulated category and  $A_1,\ldots,A_n$ be an exceptional collection in $\mathcal T$. The generation time of $G_B$ is bounded below by
 \begin{displaymath}
 \op{LL}_{\infty}\left(\bigoplus_{l \in \Z, 1 \leq i \leq n} \op{Hom}_{\mathcal T}(A_i,A_j[l])\right)-1.
\end{displaymath}
\end{proposition}

\begin{proof}
 For $\bigoplus_{l \in \Z, i} \op{Hom}_{\mathcal T}(A_i,A_j[l])$, the ideal $I$ consists of all maps between distinct objects in the exceptional collection.  By the orthogonality properties of the right and left dual, any element of $I$ is ghost for the right dual and co-ghost for the left dual. The Ghost Lemma~\ref{ghost lemma} gives the lower bound.
\end{proof}

\begin{corollary}
 Let $A_1,\ldots,A_n$ be an exceptional collection in $\mathcal A$. Assume that $A$ is formal, i.e. $A$ is quasi-isomorphic to $H(A)$, with $m_i=0$ for $i \not= 2$. The generation time of $G_B$ is equal to the Loewy length of $H(A)$.
\label{cor:gen=Loewy}
\end{corollary}

\begin{proof}
 Note that we can apply Proposition~\ref{prop:excgentime} using $H(A)$ as $A'$, and the upper bound from Proposition~\ref{prop:excgentime} and the lower bound from Proposition~\ref{prop: lower Loewy bound} coincide.
\end{proof}

\begin{example}\label{ex:quiver}
 Let us consider the quiver
\begin{center}
\begin{tikzpicture}[scale=1,level/.style={->,>=stealth,thick}]
	\node (a) at (-4.3,0) {$\bullet$};
	\node (b) at (-2.3,0) {$\bullet$};
	\node (c) at (-.3,0) {$\cdots$};
	\node (d) at (.3,0) {$\cdots$};
	\node (e) at (2.3,0) {$\bullet$};
	\node (f) at (4.3,0) {$\bullet$};
	\draw[level] (a) to (b) node at (-3.3,.325) {$a_1$};
	\draw[level] (b) to (c) node at (-1.3,.325) {$a_2$};
	\draw[level] (d) to (e) node at (1.3,.325) {$a_{n-1}$};
	\draw[level] (e) to (f) node at (3.3,.325) {$a_n$};
\end{tikzpicture}
\end{center}
 with the relation $a_n\cdots a_1 = 0$. Let $A$ denote the path algebra modulo this relation. The right dual collection to the exceptional collection formed by the projective summands of $A$ is the collection of the simple modules, $S_0,\ldots,S_n$ (up to shifting the objects). Let $S := S_0 \oplus \cdots \oplus S_n$. From \cite{Kel} $A^! = \mathbf{R}\op{End}_A(S)$ can be represented by the graded quiver
\begin{center}
\begin{tikzpicture}[scale=1,level/.style={->,>=stealth,thick}]
	\node (a) at (-4.3,0) {$\bullet$};
	\node (b) at (-2.3,0) {$\bullet$};
	\node (c) at (-.3,0) {$\cdots$};
	\node (d) at (.3,0) {$\cdots$};
	\node (e) at (2.3,0) {$\bullet$};
	\node (f) at (4.3,0) {$\bullet$};
	\draw[level] (a) to (b) node at (-3.3,.325) {$b_1$};
	\draw[level] (b) to (c) node at (-1.3,.325) {$b_2$};
	\draw[level] (d) to (e) node at (1.3,.325) {$b_{n-1}$};
	\draw[level] (e) to (f) node at (3.3,.325) {$b_n$};
	\draw[level] (a) .. controls (-3.5,1.5) and (3.5,1.5) .. (f) node at (0,1.5) {$z$};
\end{tikzpicture}
\end{center}
 with each $b_i$ of degree one and $z$ of degree two subject to the relations $b_{i+1}b_i = 0$ with the single multi-composition $m_n(b_n,\ldots,b_1) = z$. We have $\op{LL}_{\infty}(A^!) = 3$. As the right dual differs from the left dual by an application of the Serre functor, we have $\tritime(A) \leq 2$ by Proposition~\ref{prop:excgentime}. Consider the twisted complex over $A^!$
\begin{displaymath}
 (C,\delta_C) = \left( S_1 \oplus \cdots \oplus S_{n-1}, \begin{pmatrix} 0 & b_2 & 0 & \cdots & 0 \\ 0 & 0 & b_3 &  \cdots & 0 \\ & \vdots & & \vdots & \\ 0 & \cdots & 0 & 0 &  b_{n-1} \\ 0 & \cdots & 0 & 0 &  0 \end{pmatrix} \right)
\end{displaymath}
 Then the maps $b_1: S_0 \to C[1]$ and $b_n: C \to S_n[1]$ are $A$ ghost and their composition is nonzero. So $\tritime(A) \geq 2$ and hence $\tritime(A) = 2$. This demonstrates that one can have a strict inequality in Proposition~\ref{prop: lower Loewy bound}. One can also construct examples where the upper bound of Proposition~\ref{prop:excgentime} is strict. 

 One can also apply Corollary~\ref{cor:gen=Loewy} to see that $\tritime(A^!) = n-1$.
\end{example}

%\begin{example}
% Let $X$ be a smooth quadric in $\mathbb{P}^n$ defined by $q$. From \cite{Swan}, we know that $\dbcoh{X}$ has semi-orthogonal decomposition of the form
%\begin{displaymath}
% \langle\dbmod{C_0(q)}, \mathcal O(-n+1),\ldots,\mathcal O \rangle
%\end{displaymath}
% where $C_0(q)$ is the even part of the Clifford algebra on $q$ and $\mathcal O$ is exceptional. If $n$ is odd, then $\op{mod-}C_0(q)$ is equivalent to $\op{mod-}k \oplus \op{mod-}k$. Let $S_0 \oplus S_1$ be a generator of $\op{mod-}k \oplus \op{mod-}k$ with $S_0$ and $S_1$ orthogonal. Consider the generator, $S_0 \oplus S_1 \oplus \mathcal O(-n+1) \oplus \cdots \oplus \mathcal O$, of $\dbcoh{X}$. This is a strong exceptional collection and the Loewy length of the endomorphism algebra is at most $n-1$ due to the orthogonality of $S_0$ and $S_1$. By Corollary~\ref{cor:gen=Loewy}, the dual collection has generation time equal to $\dim X$. This provides a positive answer to Orlov's conjecture, Conjecture~\ref{dimension conjecture}, for smooth quadrics as the case $n$ is odd was covered in \cite{BF}.
%\end{example}

\begin{example}
 In this example, we demonstrate that the supremum of the ultimate dimension over a given birational class is infinite. We first demonstrate the method on $\P^2$ as we get a slightly sharper statement than in the general case. Let $X_1$ be the blow-up of $\mathbb{P}^2$ at a point, $p$, $E_1$ denote the exceptional curve, and $\mathcal O(H)$ the pullback of $\mathcal O(1)$ on $\mathbb{P}^2$. Let $X_2$ denote the blow-up of $X_1$ at a point on $E_1$. Let $E_2$ be the exceptional curve of this blow-up and, abusing notation, let $E_1$ be the total transform of $E_1$, i.e. the union of the strict transform of $E_1$ and $E_2$. Also, set $\mathcal O(H)$ equal to the pullback of $\mathcal O(H)$ on $X_1$. We define $X_n$ inductively as the blow-up of $X_{n-1}$ at a point on the exceptional curve of the blow-up, $X_{n-1} \to X_{n-2}$. We denote by $E_n$ the exceptional curve of the blow-up, $X_n \to X_{n-1}$ and by $E_i$, for $1 \leq i \leq n-1$, the total transforms of the $E_i$ on $X_{n-1}$. We continue to write $\mathcal O(H)$ for the pullback of $\mathcal O(H)$ to $X_{n}$. Consider the object, $G_n = \mathcal O(-2H) \oplus \mathcal O(-H) \oplus \mathcal O \oplus \mathcal O_{E_1} \oplus \cdots \oplus \mathcal O_{E_n}$. From Theorem~\ref{blow-up semi-orthogonal}, $G_n$ is a generator and it is simple to check that $\mathcal O(-2H), \ldots, \mathcal O_{E_n}$ is an exceptional collection. Note that there is a nonzero composition of length $n+2$ in $\op{End}_{X_n}(G_n)$ which corresponds to taking two sections, $s_1,s_2$, of $\mathcal O(1)$ on $\mathbb{P}^2$ not vanishing at $p$, pulling them back to $X_n$, and restricting down the chain
\begin{displaymath}
 \mathcal O(-2H) \overset{\pi^*s_1}{\to} \mathcal O(-H) \overset{\pi^*s_2}{\to} \mathcal O \to \mathcal O_{E_1} \to \mathcal O_{E_2} \to \cdots \to \mathcal O_{E_n}.
\end{displaymath}
 By Proposition~\ref{prop: lower Loewy bound}, the generation time of the dual collection is bounded below by $n+2$. In fact, this is an equality as the exceptional collection consists of $n+3$ objects.  Thus, $n+2 \in \ops{OSpec} X_n$ and $\op{udim}(X_n) \geq n+2$.

On any variety of dimension at least two, by blowing-up points iteratively, one can construct an exceptional collection with arbitrarily high Loewy length.  In doing so, one obtains a generator of an admissible subcategory of some blowup with arbitrarily large generation time.  Extending this generator by the pullback of a generator from the base, gives a generator of some blowup with arbitrarily large generation time.  
\end{example}

\begin{proposition}
Suppose $A_1,\ldots,A_n$ is a strong exceptional collection in a triangulated category, $\mathcal T$, which is the homotopy category of a triangulated $A_{\infty}$-category. Let $r$ be the projective dimension of $\op{End}_{\mathcal T}(G_A)$ and $s$ be the the Loewy length. Then $[r,s]$ is contained in the Orlov spectrum of $\mathcal T$.
\label{prop:interval spectrum}
\end{proposition}
\begin{proof}
The generation time of $G_A$ is $r$ by Theorem~\ref{CKK}. Hence, $r$ is in the Orlov spectrum.  The generator $G_B$ corresponding to the dual collection, $B_1,\ldots,B_n$, has generation time equal to the Loewy length of $\op{End}_{\mathcal T}(G_A)$ by Corollary~\ref{cor:gen=Loewy}. Hence, $s$ is also in the Orlov spectrum. As $\langle A_1,\ldots, A_n \rangle$ is a semi-orthogonal orthogonal decomposition consisting of subcategories of Rouquier dimension zero, the result follows from Theorem~\ref{semi-orthogonal gaps}.   
\end{proof}
%\begin{theorem}
%Let $\emph{Q}$ be a quiver such that the underlying graph is a Dynkin diagram of type ADE.  The dimension spectrum of $\dbmod{kQ}$ is equal to the integer interval $\{0, \ldots , n-1\}$ if $\emph{Q}$ is of type $A_n$ and $\{0, ..., n-2\}$ if $\emph{Q}$ is of type $D_n$ or $E_n$.  
%\end{theorem}

\begin{lemma} 
Let $\emph{Q}$ be a quiver such that the underlying graph is a Dynkin diagram of type $A_n$.  For each isomorphism class of indecomposable objects in $\dbmod{kQ}$ choose a representative, $M_i$.  The Loewy length of the graded algebra $\mathbf{R}\emph{End}_{kQ}(\oplus M_i)$ is $n$.
\label{lem:nilpotence of A_n quiver}
\end{lemma}

\begin{proof}
All such quivers are derived Morita equivalent so we may assume all the arrows point to the right.  Let us denote the right module generated by the $i^{th}$ vertex by $P_i$.  Then one can label the indecomposable objects by $M_{ij}:= P_i/P_{j+1}$, $1 \leq i,j \leq n$ where $M_{in} = P_i$.  If one prefers, this object can be identified with a string of 1-dimensional vector spaces beginning at the $i^{th}$ vertex and ending at the $j^{th}$ vertex with chosen isomorphisms in between.  For $j<n$, the Serre functor $S$ acts on objects which are not projective ($j<n$) by $S(M_{ij}) \cong M_{(i+1)(j+1)}[1]$ (this is merely a computation of Auslander-Reiten translation, see \cite{AR,RV}).

The morphism,
\begin{displaymath}
P_n \to \cdots \to P_1,
\end{displaymath}
is a nontrivial composition of $n-1$ nilpotent elements in $\mathbf{R}\text{End}_{kQ}(\oplus M_i)$. This gives the lower bound.

Now for any nonzero morphism from $M_{ij}$ to $M_{st}$ one has $s \leq i \leq t \leq j$, in order for it to not be an isomorphism, either $s<i$ or $t<j$.  Now, consider a nonzero sequence of morphisms in the nilradical of $\mathbf{R}\text{End}_{kQ}(\oplus M_i)$:
\begin{displaymath}
M_{i_1j_1} \to \cdots \to M_{i_aj_a}.
\end{displaymath}
We have $i_1 \leq i_m \leq j_1 \leq j_m$ for all $m$ and either $i_m$ or $j_m$ decreases.  Thus, the total length of such a sequence is at most $j_1 - j_a + i_1 - i_a$.
Now, let's add a morphism of degree one.  By Proposition \ref{prop:extendwithSerre}, we can assume a sequence of maximal length looks like:
\begin{displaymath}
M_{i_1j_1} \to \cdots \to M_{i_aj_a} \to M_{st}[1] \to \cdots \to M_{(i_1+1)(j_1+1)}[1].
\end{displaymath}
Hence the total length is at most,
\begin{displaymath}
j_1 - j_a + i_1 - i_a + s - (i_1 +1) + t - (j_1 +1) +1 \leq i_1 +1 -j_a + t - i_a  - 1 \leq -j_a +t < n.
\end{displaymath}

This is the desired upper bound.
\end{proof}

\begin{theorem}
Let $\text{Q}$ be a quiver such that the underlying graph is a Dynkin diagram of type $A_n$.  The Orlov spectrum of $\dbmod{kQ}$ is equal to the integer interval $\{0, \ldots , n-1\}$.
\end{theorem}

\begin{proof}
The upper bound is from Corollary~\ref{cor:finitely many objects bound} and Lemma~\ref{lem:nilpotence of A_n quiver}.  The set $\{1, \ldots , n-1\}$ is contained in the Orlov spectrum from Proposition~\ref{prop:interval spectrum}.  Zero is in the Orlov spectrum since the category has finitely many indecomposable objects.
\end{proof}

\section{Isolated singularities: the ungraded case} \label{sec:ungraded}

One can extract a fair bit of information about the structure of the Orlov spectrum for isolated singularities in both the graded and ungraded cases. In this section, we tackle the ungraded case leaving the graded case to the next section. Let us recall the necessary ideas.

Let $S$ be a commutative Noetherian $k$-algebra. 

\begin{definition}\label{def:category of singularities}
 The \textbf{category of singularities}, or \textbf{stable category}, of $S$ is the Verdier quotient of $\dbmod{S}$ by the subcategory consisting of all bounded complexes of finitely-generated projective modules. This is denoted by $\dsing{S}$. 
\end{definition}

Now let us assume that $(S, \mathfrak m_S)$ is a local Noetherian $k$-algebra. We say that $(S, \mathfrak m_S)$ is an \textbf{isolated singularity} if $S_{\mathfrak p}$ is a regular ring for any prime ideal, $\mathfrak p \not = \mathfrak m_S$, of $S$. The following proposition characterizes an isolated singularity purely in terms of its categories of singularities:

\begin{proposition}
 Let $(S, \mathfrak m_S)$ be a local commutative Noetherian $k$-algebra. The following are equivalent:
\renewcommand{\labelenumi}{\emph{\roman{enumi})}}
\begin{enumerate}
 \item $(S, \mathfrak m_S)$ is an isolated singularity
 \item The residue field, $k$, is a generator of $\dsing{S}$.
\end{enumerate}
\end{proposition}

This is the content of Proposition A.2 of \cite{KMV}. The implication $i) \Rightarrow ii)$ also follows immediately from the work in \cite{Sch} or the work in \cite{OrlFC}. A special case of this implication is contained in \cite{Dyc}.

Let us now provide a criterion for when $k$ strongly generates.

\begin{proposition}
 Let $(S, \mathfrak m_S)$ be a local commutative Noetherian $k$-algebra. The following are equivalent:
\renewcommand{\labelenumi}{\emph{\roman{enumi})}}
\begin{enumerate}
 \item $k$ is a strong generator of $\dsing{S}$.
 \item The natural homomorphism $S \to Z(\dsing{S})$ factors through $S/\mathfrak m_S^l$ for some $l$. 
\end{enumerate}
\label{prop:k_gen_m_nilp}
\end{proposition}

\begin{proof}
 Let us assume that $k$ is a strong generator of $\dsing{S}$. From Example~\ref{center ghost}, we see that $s(M)$ is $k$ ghost and $k$ co-ghost for any $M \in \dsing{S}$ and $s \in \mathfrak m_S$. Therefore any element of the form $s_1\cdots s_l \in \mathfrak m_S^l$ gives a ghost sequence for $k$ of length $l$. Since $k$ strongly generates, $\dsing{S} = \langle k \rangle_{l-1}$ for some ${l-1}$, it follows from the Ghost Lemma~\ref{ghost lemma} that $s_1 \cdots s_{l}(M) = s_1(M) \circ \cdots \circ s_l(M) =0$. Therefore, $\mathfrak m_S^l$ lies in the kernel of the map $S \to Z(\dsing{S})$.

 Now, assume that $\mathfrak m_S^l$ lies in the kernel of the map $S \to Z(\dsing{S})$. For an element $s \in S$, let $K(s)$ denote the complex $S \overset{s}{\to} S$. Given a collection of elements $s_1,\ldots,s_m \in S$, consider the Koszul complex associated to this collection,
\begin{displaymath}
K(s_1,\ldots,s_m) = \bigotimes_{i=1}^m K(s_i).
\end{displaymath}
Choose generators, $x_1,\ldots,x_m$, of the maximal ideal $\mathfrak m_S$. For some $l$, t	he cohomology of $K(x^l_1,\ldots,x^l_n)$ is annihilated by $\mathfrak m_S^{nl}$ as every element of $\mathfrak m_S^{nl}$ is divisible by $x_i^l$ for some $i$. Therefore, the cohomology modules of $K(x_1^1,\ldots,x_n^l) \otimes_S M$ are annihilated by $\mathfrak m_S^{nl}$ for any $M$ from $\dbmod{S}$. This implies that $K(x_1^1,\ldots,x_n^l) \otimes_S M$ lies in $\langle k \rangle_{(n+1)(ln+1)-1}$, here taken in $\dbmod{S}$. In $\dsing{S}$, $M$ is a summand of $K(x_1^1,\ldots,x_n^l) \otimes_S M$ and, hence, lies in $\langle k \rangle_{(n+1)(ln+1)-1}$.
\end{proof}

For a general ring (not necessarily of finite-type over $k$), it is unclear whether or not $k$ is always a strong generator of $\dsing{S}$. However, the following proposition covers many examples originating from algebraic geometry.  Recall that $S$ is said to be \textbf{essentially of finite type} if it is the localization of a finitely-generated $k$-algebra.

\begin{proposition}\label{prop:stronggenlocalize}
 Let $S$ be a commutative $k$-algebra that is essentially of finite type. There exists a finitely-generated $S$-module, $E$, and an $l \in \Z_{\geq 0}$ so that
\begin{displaymath}
 \op{D}(\op{Mod }S) = \langle \bar{E} \rangle_l, \hspace{2mm} \op{D}^{\op{b}}(\op{Mod }S) = \langle \tilde{E} \rangle_l, \text{ and} \hspace{2mm} \dbmod{S} = \langle E \rangle_l.
\end{displaymath}
\end{proposition}

%\sidenote{I believe Kuznetsov is suggesting that we can use Rouquier's generator and localize to get this result. I need to think about it for a bit to be sure. After thinking about it, I think you can do this but it would require at least the same amount of work. The advantage would be that you prove if D(R) has finite Rouquier dimension then so does any localization. Although there are completions of isolated singularities which are not isolated so some care should be taken in this approach.}

\begin{proof}
 Recall that Theorem 7.39 of \cite{Ro2} states that such an $E$ exists for the derived categories associated to any finitely-generated $k$-algebra. We will follow and use the proof of Theorem 7.39 in loc.~cit.  The proofs are very similar for $\op{D}(\op{Mod }S)$ and $\op{D}^{\op{b}}(\op{Mod }S)$, so we will only provide the proof of the latter and leave the proof for $\op{D}(\op{Mod }S)$ as an exercise to the reader. The statement for $\dbmod{S}$ is an immediate consequence of Corollary 6.16 and Corollary 3.13 of loc.~cit.
 
 Let $R$ be a finitely-generated $k$-algebra and $I$ a multiplicative subset of $R$ so that $S = R_I$. Let $U$ be a smooth open subset of $\op{Spec} R$ with complement determined by the ideal $J$. Let us proceed by induction on the Krull dimension of $R$. When $R$ has Krull dimension zero, the statement is a consequence of Theorem 7.39 of loc.~cit.~as $R_I$ is finitely-generated over $k$. 

 From the proof of Theorem 7.39 of loc.~cit., one has the following exact triangle in $\dbmod{R^e}$,
\begin{displaymath}
 C \to R \oplus R[1] \to D,
\end{displaymath}
 where $C$ is a perfect $R^e$-module and $D$ is a $R/J^n \otimes_k R$-module. If we localize $C$ and $D$ on the left and right by $I$, we get a triangle,
\begin{equation}
 C_I \to R_I \oplus R_I[1] \to D_I,
\label{eqn:diagdecompsing}
\end{equation}
 where $C_I$ is a perfect $R_I^e$-module and $D_I$ is a $R_I/J^nR_I \otimes_k R_I$-module.

 Let $M$ be any object of $\op{D}^{\op{b}}(\op{Mod-}R_I)$ and apply $-\overset{\mathbf{L}}{\otimes}_{R_I} M$ to Equation \ref{eqn:diagdecompsing}:
\begin{equation*}
 C_I \overset{\mathbf{L}}{\otimes}_{R_I} M \to M \oplus M[1] \to D_I \overset{\mathbf{L}}{\otimes}_{R_I} M.
\end{equation*}
 As $C_I$ is perfect, $C_I \overset{\mathbf{L}}{\otimes}_{R_I} M$ has bounded cohomology. From the long exact sequence of cohomology modules, we see that $D_I \overset{\mathbf{L}}{\otimes}_{R_I} M$ has bounded cohomology.

 From the induction hypothesis, there exists a finitely-generated $R_I/JR_I = (R/J)_I$-module, $E'$, for which $\op{D}^{\op{b}}(\op{Mod }R_I/JR_I) = \langle \tilde{E'} \rangle_l$ for some $l \in \Z_{\geq 0}$. Furthermore, $D_I \overset{\mathbf{L}}{\otimes}_{R_I} M$ lies in,
 \begin{displaymath}
 \op{D}^{\op{b}}(\op{Mod }R_I/J^nR_I) = \langle \tilde{E'} \rangle_{(n+1)(l+1)-1}.
\end{displaymath} 
 As $C_I$ lies in $\langle R_I \otimes_k R_I \rangle_t$ for some $t$, $C_I \overset{\mathbf{L}}{\otimes}_{R_I} M$ lies in $\langle \tilde{R}_I \rangle_t$. This implies that $M \in \langle \widetilde{R_I \oplus E'} \rangle_{(t+1)(n+1)(l+1)-1}$. We can take $E = R_I \oplus E'$. 
\end{proof}

\begin{proposition}
 Let $(S, \mathfrak m_S)$ be a local commutative $k$-algebra that is essentially of finite type over $k$. There exists a finitely-generated $\hat{S}$-module, $E$, and an $l \in \Z_{\geq 0}$ so that
\begin{displaymath}
 \op{D}(\op{Mod }\hat{S}) = \langle \bar{E} \rangle_l, \hspace{2mm} \op{D}^{\op{b}}(\op{Mod }\hat{S}) = \langle \tilde{E} \rangle_l, \text{ and} \hspace{2mm} \dbmod{\hat{S}} = \langle E \rangle_l,
\end{displaymath}
 where $\hat{S}$ is the completion of $S$ at $\mathfrak m_S$.
\end{proposition}

\begin{proof}
 The argument is the same as in the proof of Proposition~\ref{prop:stronggenlocalize} above.
\end{proof}

\begin{corollary}
 If $(S,\mathfrak m_S)$ is a local commutative $k$-algebra essentially of finite type over $k$, then $\dsing{S}$ has finite Rouquier dimension. The same is true for $\dsing{\hat{S}}$.
\end{corollary}

Combining the results above, we get the following characterization of an isolated singularity when the ring is essentially of finite type:

\begin{theorem}
 Let $(S, \mathfrak m_S)$ be a local commutative $k$-algebra essentially of finite type over $k$. The following are equivalent:
\renewcommand{\labelenumi}{\emph{\roman{enumi})}}
\begin{enumerate}
 \item $(S, \mathfrak m_S)$ is an isolated singularity.
 \item $k$ is a strong generator for $\dsing{S}$.
 \item The natural map $S \to Z(\dsing{S})$ factors through $S/\mathfrak m_S^d$ for some $d \in \N$.
\end{enumerate}
\label{thm:charisosing}
\end{theorem}

\begin{proof}
 We know that $ii)$ and $iii)$ are equivalent by Proposition \ref{prop:k_gen_m_nilp}. Since $S$ is essentially of finite type, Proposition \ref{prop:stronggenlocalize} says we have a strong generator. Thus, if $k$ is a generator, $k$ must be a strong generator.
\end{proof}

While Theorem~\ref{thm:charisosing} is an interesting characterization of an isolated singularity, it provides no control over the generation time of $k$ or over the Orlov spectrum of $\dsing{S}$. To get such information, we restrict to the case of an isolated hypersurface singularity.

A local Noetherian $k$-algebra, $(S, \mathfrak m_S)$, is called a \textbf{hypersurface singularity} if $S$ is isomorphic to $R/(w)$ with $(R, \mathfrak m_R)$ a Noetherian, regular local $k$-algebra and $w$ lies in $\mathfrak m_R$. The multiplicity of $w$ will be the minimal $l$ so that $w \in \mathfrak m_R^l$. If $(S, \mathfrak m_S)$ is a hypersurface singularity it is Gorenstein, in particular Cohen-Macaulay. 

There are two additional constructions of $\dsing{S}$ which are useful to consider. Recall that a module, $M$, over $S$ is called a \textbf{maximal Cohen-Macaulay} module, or a MCM module for short, if the depth of $M$ is equal to the Krull dimension of $S$.

For the first construction, let $\op{MCM}(S)$ be the full subcategory of $\op{mod }S$ consisting of MCM modules. $\underline{\op{MCM}}(S)$ is a category with the same objects as $\op{MCM}(S)$ but with
\begin{gather*}
 \op{Hom}_{\underline{\op{MCM}}(S)}(M,N) = \op{Hom}_S(M,N)/\sim
\end{gather*}
where $f \sim g$ if there exists maps $p:M \to P$ and $q:P \to N$ with $f-g = qp$ and $P$ projective.

In the second construction, the objects are sequences of $R$-modules,
\begin{equation*}
 P_0 \overset{A}{\to} P_1 \overset{B}{\to} P_0,
\end{equation*}
with $P_i$ finitely-generated projective $R$-modules, $AB = w \op{id}_{P_1}$, and $BA = w \op{id}_{P_0}$. Such a sequence is a called a \textbf{matrix factorization}. For simplicity, we denote a matrix factorization $(P_0,P_1,A,B)$ by $P$ and let $A_P$ and $B_P$ denote the maps in the matrix factorization. A morphism between two matrix factorizations, $P$ and $Q$, consists of $R$-module maps $f_0:P_0 \to Q_0$ and $f_1: P_1 \to Q_1$ making the following diagram commutative:

\begin{center}
\begin{tikzpicture}[description/.style={fill=white,inner sep=2pt}]
\matrix (m) [matrix of math nodes, row sep=3em, column sep=2.5em, text height=1.5ex, text depth=0.25ex]
{ P_0 & P_1 & P_0  \\
  Q_0 & Q_1 & Q_0 \\ };
\path[->,font=\scriptsize]
(m-1-1) edge node[auto] {$ A_P $} (m-1-2)
        edge node[auto] {$ f_0 $} (m-2-1)
(m-1-2) edge node[auto] {$ B_P $} (m-1-3)
        edge node[auto] {$ f_1 $} (m-2-2)
(m-1-3) edge node[auto] {$ f_0 $} (m-2-3)
(m-2-1) edge node[auto] {$ A_Q $} (m-2-2)
(m-2-2) edge node[auto] {$ B_Q $} (m-2-3);
\end{tikzpicture}
\end{center}

A homotopy between two morphisms $f,g: P \to Q$ is a pair of maps $h_0: P_0 \to Q_1$ and $h_1:P_1 \to Q_0$ so that $f_0 - g_0 = B_Q h_0 + h_1 A_P$ and $f_1 - g_1 = A_Q h_1 + h_0 B_P$. The category of matrix factorization of $w$, $\op{MF}(w)$, has matrix factorizations as objects and has homotopy classes of morphisms between $P$ and $Q$ as morphism sets.

In both of these descriptions, the resulting category is naturally triangulated.  We have the following result, see \cite{Buc86} or \cite{Orl04}:

\begin{theorem}
 For an isolated hypersurface singularity, $S$, the categories $\dsing{S}$, $\underline{\op{MCM}}(S)$, and $\op{MF}(w)$ are all equivalent as triangulated categories.
\label{thm:holytrinity}
\end{theorem}

We draw from this two useful corollaries.

\begin{corollary}
 Every object in $\dsing{S}$ is isomorphic to a MCM module.
\end{corollary}

\begin{proof}
 The equivalence of $\dsing{S}$ and $\underline{\op{MCM}}(S)$ is induced by the inclusion,

\begin{displaymath}
\op{MCM}(S) \hookrightarrow \op{Ch}(\op{mod }S),
\end{displaymath} which sends an MCM module, $M$, to the complex
\begin{displaymath}
 \cdots \to 0 \to M \to 0 \to \cdots 
\end{displaymath}
 with $M$ in degree zero.
\end{proof}

Let $(\partial w) = (\partial_1 w,\ldots, \partial_n w)$.

\begin{corollary}
 The natural map $S \to Z(\dsing{S})$ factors through the projection $S \to S/(\partial w)$. 
\label{cor:kernelJacobian}
\end{corollary}

\begin{proof}
 We consider the category $\op{MF}(w)$. If $P$ is a matrix factorization, then, taking the $i$-th derivatives of $AB = w \id_{P_1}$ and $BA = w \id_{P_0}$, we get $\partial_i A B + A \partial_i B = \partial_i w  \op{id}_{P_1}$ and $\partial_i B A + B \partial_i A = \partial_i w\op{id}_{P_0}$.  This means that $(\partial_iA, \partial_iB)$ is a homotopy between $\partial_1 w \text{ and } 0$.
\end{proof}

Recall the \textbf{Loewy length} of a local Artinian ring, $R$, is the minimal $l$ for which $\mathfrak m_R^l=0$. Denote this as $\op{LL}(R)$. For an isolated hypersurface singularity, the Tjurina algebra, $S/\partial w$, is Artinian.  We can apply the ideas of Proposition~\ref{prop:k_gen_m_nilp} to prove the following: 

\begin{proposition}
 Let $(S, \mathfrak m)$ be an isolated hypersurface singularity. The generation time of $k$ in $\dsing{S}$ is bounded above by $2\op{LL}(S/(\partial w)) - 1$. In particular, $\dsing{S}$ has finite Rouquier dimension.
\label{prop:kstronggen}
\end{proposition}

\begin{proof}
 Let $M$ be any MCM module over $S$ and consider the Koszul complex,
 \begin{displaymath}
K(\partial w) := K(\partial_1 w,\ldots,\partial_n w).
\end{displaymath}
As the Krull dimension of $S/(\partial w)$ is zero and $M$ is MCM-module, there is an $M$-sequence of length $n-1$ in $(\partial w)$. By \cite{Mat89} Theorem 16.8, the cohomology of,
\begin{displaymath}
K(\partial w) \otimes_S M=: K(M,\partial w),
\end{displaymath}
vanishes except for degrees zero and one. Furthermore, $H_i(K(M,\partial w))$ is a module over $S/(\partial w)$.  For any $S/(\partial w)$-module, $L$, we have a filtration:
\begin{displaymath}
0 = \mathfrak m_{S/(\partial w)}^{\op{LL}(S/(\partial w))} L \subseteq \cdots \subseteq \mathfrak m_{S/(\partial w)} L \subseteq L.
\end{displaymath} 
The quotients of this filtration are direct sums of the residue field.
Therefore, we have 
\begin{equation*}
 H_i(K(M,\partial w)) \in \langle k \rangle_{\op{LL}(S/(\partial w))-1} \text{ and } K(M,\partial w) \in \langle k \rangle_{2\op{LL}(S/(\partial w))-1}
\end{equation*}
 in $\dbmod{S}$. In $\dsing{S}$, by Corollary~\ref{cor:kernelJacobian}, the partial derivatives of $w$ vanish.  Hence, $M$ is a summand of $K(M,\partial w)$. Thus, $\dsing{S} = \langle k \rangle_{2\op{LL}(S/(\partial w))-1}$.
\end{proof}

\begin{remark}
 Strong generation of $k$ also follows from work in \cite{Dyc}.
\end{remark}

Our next goal is to study the Orlov spectrum of $\dsing{S}$. Before we wade into the case of a general hypersurface, let us fully analyze the stable category of the ring $A_{n-1} = k[u]/(u^n)$. See also \cite{Orl04}. From the classification of modules over a PID, we know the only indecomposable modules are
\begin{displaymath}
 k[u]/(u^n),k[u]/(u^{n-1}),\ldots,k[u]/(u),0.
\end{displaymath}
Any morphism in $\op{mod }A_{n-1}$ from $k[u]/(u^i)$ to $k[u]/(u^j)$ is a linear combination of the maps
\begin{align*}
 \alpha_{i,j}^l: k[u]/(u^i) & \to k[u]/(u^j) \\
 1 & \mapsto u^l
\end{align*}
for $\max(0,j-i) \leq l < j$. The map, $\alpha_{i,j}^l$, factors through $k[u]/(u^n)$ if and only if $l \geq n-i$. In $\dsing{A_{n-1}}$, we let $V_i$ stand for the image of $k[u]/(u^i)$. The morphism space between $V_i$ and $V_j$ is spanned by the images of $\alpha_{i,j}^l$ with $\max(0,j-i) \leq l < \min(j,n-i)$. Let us compute the cones. We have an exact sequence:
\begin{equation}
 0 \to k[u]/(u^{\max(0,i-j+l)}) \to k[u]/(u^i) \overset{\alpha_{i,j}^l}{\to} k[u]/(u^j) \to k[u]/(u^l) \to 0. 
\label{eqn:A_nses}
\end{equation}
\begin{lemma}
 The extension in Equation~\ref{eqn:A_nses} is trivial.
\label{lem:A_nsplit}
\end{lemma}

\begin{proof}
 We can assume that $i-j+l$ is non-negative. Let us take a free resolution of $k[u]/(u^l)$ and choose a homotopy class of chain maps between the free resolution and the exact sequence~\ref{eqn:A_nses}.
\begin{center}
\begin{tikzpicture}[description/.style={fill=white,inner sep=2pt}]
\matrix (m) [matrix of math nodes, row sep=3em, column sep=2.5em, text height=1.5ex, text depth=0.25ex]
{ k[u]/(u^n) & k[u]/(u^n) & k[u]/(u^n) & k[u]/(u^l) & 0 \\
  k[u]/(u^{i-j+l}) & k[u]/(u^i) & k[u]/(u^j) & k[u]/(u^l) & 0 \\ };
\path[->,font=\scriptsize]
(m-1-1) edge node[auto] {$ \alpha_{n,n}^{n-l} $} (m-1-2)
        edge node[auto] {$ \lambda $} (m-2-1)
(m-1-2) edge node[auto] {$ \alpha_{n,n}^{l} $} (m-1-3)
        edge node[auto] {$ \alpha_{n,i}^0 $} (m-2-2)
(m-1-3) edge node[auto] {$ \alpha_{n,l}^{0} $} (m-1-4)
        edge node[auto] {$ \alpha_{n,j}^0 $} (m-2-3)
(m-1-4) edge (m-1-5)
        edge node[auto] {$ \alpha_{l,l}^0 $} (m-2-4)
(m-2-1) edge node[auto] {$ \alpha_{i-j+l,i}^{j-l} $} (m-2-2)
(m-2-2) edge node[auto] {$ \alpha_{i,j}^{l} $} (m-2-3)
(m-2-3) edge node[auto] {$ \alpha_{j,l}^{0} $} (m-2-4)
(m-2-4) edge (m-2-5);
\end{tikzpicture}
\end{center}
Since $l < n - i$, $\alpha_{n,i}^0 \circ \alpha_{n,n}^{n-l} = \alpha_{n,i}^{n-l}$ is zero. We can take $\lambda$ to be zero which proves the claim.
\end{proof}

In $\dsing{A_{n-1}}$, we get triangles 
\begin{center}
\begin{tikzpicture}[description/.style={fill=white,inner sep=2pt}]
\matrix (m) [matrix of math nodes, row sep=3em, column sep=-0.75em, text height=1.5ex, text depth=0.25ex]
{ V_i & & V_j \\
   & V_{\max(0,i-j+l)}[1] \oplus V_l \\ };
\path[->,font=\scriptsize]
(m-1-1) edge node[auto] {$ \alpha_{i,j}^{l} $} (m-1-3)
(m-1-3) edge (m-2-2)
(m-2-2) edge node[sloped] {$ | $} (m-1-1);
\end{tikzpicture}
\end{center}
We also have isomorphisms, $k[u]/(u^i) \cong k[u]/(u^{n-i})[1]$, coming from the short exact sequences,
\begin{equation*}
 0 \to k[u]/(u^{n-i}) \overset{\alpha^i_{n-i,n}}{\to} k[u]/(u^n) \to k[u]/(u^i) \to 0.
\end{equation*}

\begin{theorem}
 The Orlov spectrum of $\dsing{A_{n-1}}$ is
\begin{displaymath}
 \left\lbrace 0,1, \ldots, \left\lceil \frac{\lfloor n/2 \rfloor}{s}  \right\rceil -1, \dots ,\left\lceil  \frac{\lfloor n/2 \rfloor}{2} \right\rceil -1, \lfloor n/2 \rfloor -1 \right\rbrace
\end{displaymath}
 where $\lfloor \alpha \rfloor$ is the greatest integer less than $\alpha$ and $\lceil \alpha \rceil$ is the least integer greater than $\alpha$.
\label{thm:Anspec}
\end{theorem}

\begin{proof}
 Let $G$ be a generator for $\dsing{A_{n-1}}$. Without loss of generality, we can assume that
\begin{displaymath}
 G = \bigoplus_{i \in I \subset \lbrace 1, \ldots, \lfloor n/2 \rfloor \rbrace} V_i.
\end{displaymath}
 Let
\begin{displaymath}
 \delta(t) = \max\lbrace j | V_j \in \langle G \rangle_t, 0 \leq j \leq \lfloor n/2 \rfloor \rbrace.
\end{displaymath}
 We first show that
\begin{displaymath}
 \tritime(G) \leq \begin{cases} \max\lbrace \lceil  \frac{\lfloor n/2 \rfloor}{\delta(0)}  \rceil - 1,1\rbrace & \langle G \rangle_0 \not = D_{\op{sg}}(A_{n-1}) \\ 0 & \langle G \rangle_0 = D_{\op{sg}}(A_{n-1}). \end{cases}
\end{displaymath}
 Assume that $V_j$, $j \leq \lfloor n/2 \rfloor$, lies in $\langle G \rangle_t$. Without loss of generality we can assume that $j \geq \delta(0)$. To make new indecomposables, the possible cones we could take involve the pairs $(i,j), (i,n-j), (n-i,j), (n-i,n-j)$ with $i \in I$. If we use the pair $(i,j)$, we get indecomposable objects $V_t$ with $\max(0,j-i) \leq t < j$ and $\max(0,i-j) \leq t < i$ in the next step. If we use the pair $(i,n-j)$, we get the indecomposable objects $V_t$ with $n-j-i \leq t < \min(n-j,n-i)$ and $0 \leq t < \min(i,j)$.

 We see that $V_0,\ldots,V_{i+j}$ lies in $\langle G \rangle_{t+1}$. Therefore,
 \begin{displaymath}
 \delta(t+1) \geq \min(\delta(t)+\delta(0),\lfloor n/2 \rfloor),
 \end{displaymath}
and after the zeroth step, if $\langle G \rangle_t$ contains $V_j$ for $j \leq \lfloor n/2 \rfloor$, then it contains $V_s$ for $1 \leq s \leq j$. This gives the claimed upper bound.

 To demonstrate that the lower bound holds, we note that $x^l$ annihilates $G$ when $l \geq \delta(0)$. By Example~\ref{center ghost}, $x^l(V_{\lfloor n/2 \rfloor})$ is $G$ ghost. Furthermore, $(x^l)^{\lceil  \frac{\lfloor n/2 \rfloor}{l}  \rceil -1}(V_{\lfloor n/2 \rfloor})$ is nonzero. Therefore, by the Ghost Lemma~\ref{ghost lemma}, $\lceil  \frac{\lfloor n/2 \rfloor}{l}  \rceil -1$ is a lower bound for the generation time of $G$.
	
 Consequently,
\begin{displaymath}
 \tritime(G) = \begin{cases} \max\lbrace\lceil  \frac{\lfloor n/2 \rfloor}{\delta(0)}  \rceil - 1,1\rbrace & \langle G \rangle_0 \not = D_{\op{sg}}(A_{n-1}) \\ 0 & \langle G \rangle_0 = D_{\op{sg}}(A_{n-1}) \end{cases}
\end{displaymath}
\end{proof}

% As a corollary of the proof, we have the following.
% 
% \begin{corollary}
% The level of the residue field is at most one for any generator of $D_{\op{sg}}(A_m)$.
% \label{cor:lvlofk}
% \end{corollary}

Let us return to the case of a general isolated hypersurface singularity, see also \cite{Tak} Section 5.

\begin{lemma}
 Let $(S, \mathfrak m_S)$ be a hypersurface singularity and let $M$ be a MCM module over $S$. For a generic choice of a regular system of parameters on $R$, $y_1,\ldots,y_n$, the first $n-1$ parameters, $y_1,\ldots,y_{n-1}$, form both an $S$-regular and an $M$-regular sequence and the quotient $S/(y_1,\ldots,y_{n-1})S$ is isomorphic to a zero dimensional hypersurface singularity. Moreover, the multiplicity of $w$ in $R$ is the same as the multiplicity of $\bar{w}$ in $R/(y_1,\ldots,y_{n-1})$.
\label{lem:slicetoAn} 
\end{lemma}

%\sidenote{We can prove that m equals the order of vanishing minus one keeping the coordinate choices generic. Should we include this?}
%\sidenote{It may be interesting, using whatever object you get here, you can control how long it takes $M/y_1/ ldots /y_n$ to generate any other module $N/y_1/ ldots /y_n$.  Then you can build $N$ using nilpotence of $y_i$.  You'll get a different and generally worse bound I imagine.}

\begin{proof}
 Recall that a sequence of elements, $s_1,\ldots,s_i$, is $M$-regular if $s_j$ has zero annihilator in $M/(s_1,\ldots,s_{j-1})M$. Now, $x_1,\ldots,x_n$ is a regular system of parameters for $R$ if $(x_1,\ldots,x_n) = \mathfrak m_R$ with $n$ equal to the Krull dimension of $R$. Recall that $x_1,\ldots,x_n$ is a regular system of parameters for $R$ if and only if the images of $x_1,\ldots,x_n$ form a basis for $\mathfrak m_R/ \mathfrak m_R^2$, see \cite{Mat89} Theorem 14.2.

 We prove the results involving $S$ and then note that the same choices work to establish the result about $M$. We proceed by induction on $n$. The case $n=1$ is clear.

 Assume we know the result below $n-1$ and consider the case of $n$. $w$ has a unique factorization (in $R$) into irreducible elements. Let us denote them by $w_1,\ldots,w_t$. Let $x$ be an element of $R$ that projects to a nonzero vector in $\mathfrak m_R/ \mathfrak m_R^2$. It is clear that $x$ is irreducible and is a zero divisor in $S$ if and only if it equals some $w_i$. The associated graded ring, $\op{gr}_{\mathfrak m_R}(R)$, is isomorphic to a polynomial ring over $k$ in $n$ variables, \cite{Mat89} Theorem 17.10. Let $d$ be the multiplicity of $w$ and denote the image of $w$ in $\op{gr}_{\mathfrak m_R}(R)$ by $w_d$. 

 If $n$ is greater than one, we can choose an element $u$ of $R$ with nonzero image in $\mathfrak m_R/ \mathfrak m_R^2$ so $u$ is a not a zero-divisor in $S$ and the image of $u$ in $\op{gr}_{\mathfrak m_R}(R)$ does not divide $w_d$. Now, complete $u$ to a regular system of parameters for $R$, $u,u_2,\ldots,u_n$. We find that $S/uS$ is another hypersurface singularity to which we can apply the induction hypothesis. 

 Let $\mathfrak{p}$ be an associated prime for $M$. The depth of $M$ is bounded above by the dimension of $A/\mathfrak{p}$, Theorem 17.2 \cite{Mat89}. If $M$ is a MCM module, then the height of $\mathfrak{p}$ cannot be more than zero. By Krull's theorem, $\mathfrak{p}$ cannot contain a non-zerodivisor. Thus, $\mathfrak{p}$ is in the ideal generated by the $w_1,\ldots,w_t$. Since our choices of a regular system of parameters avoids each $w_i$, they also provide an $M$-sequence.
\end{proof}

\begin{lemma}\label{lem:not free}
 Let $(S, \mathfrak m_S)$ be an isolated hypersurface singularity and let $M$ be a module of infinite projective dimension over $S$. If $x \in S$ is a nonunit and $S$ and $M$-regular, then $M/xM$ is a module of infinite projective dimension over $S/(x)$. 
\end{lemma}

\begin{proof}
 Note that $S/(x)$ vanishes in $\dsing{S}$ as it is quasi-isomorphic to the cone of $x(S):S \to S$ and hence perfect. Also note that the morphism, $x(M): M \to M$, in $\dsing{S}$ is a zero divisor by Theorem \ref{thm:charisosing}. 

 Assume that $M/xM$ has finite projective dimension as an $S/(x)$ module. Then, $M/xM$ vanishes in $\dsing{S}$. As $M/xM$ is quasi-isomorphic to the cone of $x(M)$, we see that $x(M)$ is an isomorphism in $\dsing{S}$ and cannot be a zero divisor.
\end{proof}

\begin{lemma}
 Any zero dimensional hypersurface singularity, $S = R/(w)$, is isomorphic to $A_{d-1}$, where $d$ is the multiplicity of $w$.
\label{classification of zero dim}
\end{lemma}

\begin{proof}
 As $S$ is zero dimensional, completion does not change the ring. Thus, $S$ is isomorphic to $\hat{R}/(w)$. Any complete, regular, local, Noetherian ring of dimension one is isomorphic to the formal power series ring in one variable $k[[u]]$ with the uniformizing parameter of $\hat{R}$ getting sent to $u$, \cite{Mat89} Theorem 29.7. A simple change of variables takes $w$ to $u^d$. 
\end{proof}

We now use these lemmas to facilitate a reduction from a general isolated hypersurface singularity to an $A_n$-singularity.

\begin{lemma} 
 Let $(S, \mathfrak m_S)$ be an isolated hypersurface singularity and let $M$ be any non-zero object of $\dsing{S}$. The level of the residue field of $(S, \mathfrak m_S)$ with respect to $M$ is at most $\dim S+1$.  In particular, $M$ is a strong generator of $\dsing{S}$. 
\label{lem: all generate}
\end{lemma}

\begin{proof}
 Let $S$ be isomorphic to $R/(w)$. From Lemmas~\ref{lem:slicetoAn} and \ref{classification of zero dim}, we know we can choose a regular system of parameters, $x_1,\ldots,x_n$, with $x_1,\ldots,x_{n-1}$ a $S$-regular and a $M$-regular sequence and so that $S/(x_1,\ldots,x_{n-1})$ is isomorphic to $A_{d-1}=k[u]/(u^{d})$ where $d$ is multiplicity of $w$. Note that $M/(x_1,\ldots,x_{n-1})M$ cannot be free by Lemma \ref{lem:not free}.

 Let $K(x) = K(x_1,\ldots,x_n)$ and $K(M,x) = K(x) \otimes_S M$.   Notice that $K(M,x)$ is quasi-isomorphic to the complex $M/(x_1,\ldots,x_{n-1})M \overset{x_n}{\to} M/(x_1,\ldots,x_{n-1})M$. Writing $x_n$ as $\alpha_1 u + \cdots + \alpha_m u^{d-1}$, one sees that $M/(x_1,\ldots,x_{n-1})M \overset{x_n}{\to} M/(x_1,\ldots,x_{n-1})M$ is the composition of $M/(x_1,\ldots,x_{n-1})M \overset{u}{\to} M/(x_1,\ldots,x_{n-1})M$ and an automorphism of $M/(x_1,\ldots,x_{n-1})M$. The octahedral axiom tells us that the cone of,
 \begin{displaymath}
 M/(x_1,\ldots,x_{n-1})M \overset{u}{\to} M/(x_1,\ldots,x_{n-1})M,
\end{displaymath}
is isomorphic to the cone of,
\begin{displaymath}
M/(x_1,\ldots,x_{n-1})M \overset{x_n}{\to} M/(x_1,\ldots,x_{n-1})M.
\end{displaymath}
 As $M/(x_1,\ldots,x_{n-1})M$ is nonfree, a direct calculation for the $A_{d-1}$ singularity, see the proof of Theorem~\ref{thm:Anspec}, shows that the cone of,
 \begin{displaymath}
M/(x_1,\ldots,x_{n-1})M \overset{u}{\to} M/(x_1,\ldots,x_{n-1})M
\end{displaymath}
 is quasi-isomorphic to a sum of shifts of $k$. Hence $k$ is a summand of $K(M,x)$ which manifestly lies in $\langle M \rangle_n$.
 
 The above tells us that $M$ generates $k$, and, by Theorem~\ref{thm:charisosing}, $k$ generates $\dsing{S}$.  It follows that $M$ is a strong generator.
\end{proof}

\begin{remark}
 Lemma~\ref{lem: all generate} is not true for complete intersections. For example, consider the ring $S = k[x,y]/(x^2,y^2)$. The module $k[x]/(x^2)$ is nonzero in $\dsing{S}$ but $k[y]/(y^2)$ is orthogonal to it. 
\label{rmk:CI}
\end{remark}

\begin{remark}
 Let $M$ be a MCM module. The arguments in the proof of Lemma~\ref{lem: all generate} give the following statement: $M$ is a generator of $\dsing{S}$ if and only $M \overset{\textbf{L}}{\otimes}_R k \in \langle k \rangle_0$ in $\dbmod{S}$. Does this statement hold for complete intersections? The authors know of no counterexample.
\end{remark}

Combining Proposition~\ref{prop:kstronggen} and Lemma~\ref{lem: all generate} gives us the following theorem:

\begin{theorem}
 Let $(S, \mathfrak m_S)$ be an isolated hypersurface singularity. The ultimate dimension of $\dsing{S}$ is bounded by $2(\dim S + 2)\op{LL}(S/(\partial w)) - 1$.
\label{thm:isohypspecbound}
\end{theorem}

%\sidenote{I'm not sure I like this example because its obviously bounded anyway and I think the number of objects is a better bound for $D_n$}

%For an example, let us consider the $D_n$-singularities $D_n = k[x,y]/(x^2y+y^{n-1})$. While there is classification of the MCM modules over $D_n$ and morphisms between them in the stable category, it is a bit more complicated than in the case of $A_n$. For a quicker and rougher result, we can use the arguments used to prove Theorem~\ref{thm:isohypspecbound} to give an upper on the Orlov spectrum of $\dsing{D_n}$. 

%The function, $x- y$, is a regular element of $D_n$. If we quotient $D_n$ by $x - y$, we get a ring isomorphic to $A_4$. In $\dsing{A_4}$, every object is a shift of the residue field. So the level of $k$ with respect to any generator of $\dsing{D_n}$ is at most one. The Jacobian ideal of $D_n$ is $(xy,x^2+ny^{n-1})$ and we have $x^3 = xy = y^n = 0$ with $y^{n-1} \not = 0$. Thus, the Loewy length of $k[x,y]/(xy,x^2+ny^{n-1},x^2y+y^n)$ is $n$. This gives the following:

%\begin{proposition}
% \begin{displaymath}
% \op{OSpec } \dsing{D_n} \subset \lbrace 0,\ldots,4n-1 \rbrace.
%\end{displaymath}
%\end{proposition}

For an example, let us consider the ring,
\begin{displaymath}
 S_g = k[x,y,z]/(x^{2g+1}+y^{2g+1}+z^{2g+1}-xyz). 
\end{displaymath}
for $g > 1$. Let $w_g = x^{2g+1}+y^{2g+1}+z^{2g+1}-xyz$. We can take $x-z,y-z$ as a regular sequence and $S_g/(x-z,y-z)$ is isomorphic to $A_2$. The level of residue field is at most two for any generator of $\dsing{S_g}$. The Jacobian ideal of $S_g$ is $((2g+1)x^{2g}-yz,(2g+1)y^{2g}-xz,(2g+1)z^{2g}-xy)$. The Loewy length of $S_g/(\partial w_g)$ is $2g+1$. 

There is $\Z/(2g+1)\Z$ action on $S_g$ and it is proven in \cite{Sei08}, for $g=2$, and in \cite{Ef}, for $g \geq 2$, that the idempotent-completion of the $\Z/(2g+1)\Z$-equivariant singularity category, $\dGsing{S_g}{\Z/(2g+1)\Z}$, is equivalent to the idempotent-completion of the derived Fukaya category of a genus $g$ Riemann surface, $\op{D}^{\pi}\op{Fuk}(\Sigma_g)$. In light of Example~\ref{dense equivariant}, we can use our results to control the generation time of certain generators of $\op{D}^{\pi}\op{Fuk}(\Sigma_g)$ (the notation in the following proof can be found in this example).  More precisely, recall that symplectically, the surface, $\Sigma_g$, admits a $\Z/(2g+1)\Z$-branched cover over an orbifold $\mathbb{P}^1$. Let $\psi: \Sigma_g \to \Sigma_g$ be a generator of the covering group. We now have the following result:

\begin{proposition}
 Let $M$ be any nonzero object of $\op{D}^{\pi}\op{Fuk}(\Sigma_g)$. Then, $\bigoplus_{i=0}^{2g} \psi^i(M)$ is a generator of $\op{D}^{\pi}\op{Fuk}(\Sigma_g)$ and its generation time is bounded by $12g+5$. 
\end{proposition}

\begin{proof}
By Lemma~\ref{lem: all generate}, $\text{For}(M)$ generates $\dsing{S_g}$.  By Example~\ref{dense equivariant}, the functor, $\text{Inf}$, is dense and hence $\bigoplus_{i=0}^{2g} \psi^i(M) \cong \text{Inf}(\text{For}(M))$ generates with,
\begin{displaymath}
\tritime(\bigoplus_{i=0}^{2g} \psi^i(M)) = \tritime(\text{For}(\bigoplus_{i=0}^{2g} \psi^i(M))).
\end{displaymath}
The level of $k$ with respect to any object of $D_{\op{sg}}(S_g)$ is at most two and the generation time of $k$ is at most $4g+1$. Thus, $\tritime(\text{For}(\bigoplus_{i=0}^{2g} \psi^i(M))) \leq 12g+5$.
\end{proof}

%If $G$ is a finite group acting on a $k$-algebra, $S$, $D^{b,G}(\op{mod-}S)$ is equivalent to $D^b(\op{mod-}(S \# k[G]))$ where $S \# k[G]$ is %skew-group algebra. We have a functor
%\begin{align*}
% -\otimes_k k[G] : D^b(\op{mod-}S) & \to D^b_G(\op{mod-}S) \\
% M & \mapsto M \otimes_k k[G].
%\end{align*}
%It is left adjoint to the forgetful functor $D^b_G(\op{mod-}S) \to D^b(\op{mod-}S)$.

%\begin{lemma}
% The generation time of $M$ in $D_{\op{sg}}(S)$ equals the generation time of $M \otimes_k k[G]$ in $D_{\op{sg}}^G(S)$. 
%\label{lem:equi_gen}
%\end{lemma}

%\begin{proof}
% As $M \otimes_k k[G]$ is isomorphic to a sum of a copies of $M$ in $D^b(\op{mod-}S)$, we see that the forgetful functor is dense and $\tritime(M %\otimes_k k[G]) \geq \tritime(M)$. If $M$ is a $S \# k[G]$-module, then $M$ is summand of $M \otimes_k k[G]$ in $\op{mod-}(S \# k[G])$ as the %characteristic of $k$ is zero. Thus, the image of $D_{\op{sg}}(S)$ is dense in $D^G_{\op{sg}}(S)$. So $\tritime(M \otimes_k k[G]) \leq \tritime(M)$.
%\end{proof}

\section{Isolated singularities: the graded case} \label{sec:graded}

Most of the results in Section~\ref{sec:ungraded} can be adapted to the graded case in a straightforward manner. When we combine these results with Orlov's results relating derived categories of coherent sheaves to graded categories of singularities, many interesting and nontrivial statements emerge. So, let us begin by recalling Orlov's results from \cite{Orl09}. We let $A = \bigoplus_{n \geq 0} A_n$ be a graded Noetherian $k$-algebra with $A_0 = k$. We write, $\op{gr }A$, for the abelian category of finitely-generated graded $A$-modules.  The morphisms in this category are taken to be degree zero $A$-module homomorphisms. The category has an internal Hom denoted $\underline{\op{Hom}}$. For any graded module, $M$, we can form a new graded module, $M(1)$, with $M(1)_l = M_{l+1}$. Recall that $A$ is AS-Gorenstein if $A$ has finite injective dimension $n$ and $\underline{\op{Ext}}^i_{\op{gr} A}(k,A) = 0$ for $i \not = 0$ and $\underline{\op{Ext}}^n_{\op{gr} A}(k,A) = k(a)$. We call, $a$, the \textbf{Gorenstein parameter} of $A$. We have the maximal ideal $\mathfrak m_A = \bigoplus_{l > 0} A_l$.

Sitting inside of $\op{gr }A$, we have the full subcategory of finite-dimensional modules (over $k$), $\op{tors }A$. Inside of $\dbgr{A}$, we have two thick triangulated subcategories: $\op{perf }A$, the full subcategory consisting of all bounded complexes of finite rank free $A$-modules, and, $\dbtors{A}$, the full subcategory consisting of all complexes quasi-isomorphic to a bounded complex of torsion modules. 

\begin{definition}
 Let $\dbqgr{A}$ denote the Verdier quotient of $\dbgr{A}$ by $\dbtors{A}$. Let $\dgrsing{A}$ denote the Verdier quotient of $\dbgr{A}$ by $\ops{perf}A$. We call, $\dgrsing{A}$, the \textbf{graded category of singularities} of $A$.
\end{definition}

In \cite{Orl09}, Orlov proves the following useful theorem relating $\dbqgr{A}$ and $\dgrsing{A}$:
\begin{theorem}
 For any $i \in \Z$ we have the following statements:
\renewcommand{\labelenumi}{\emph{\roman{enumi})}}
\begin{enumerate}
 \item If $a > 0$, there is a fully-faithful functor, $\Psi_i: \dgrsing{A} \to \dbqgr{A}$, and a semi-orthogonal decomposition,
\begin{displaymath}
 \dbqgr{A} \cong \left\langle A(-i-a+1),\ldots,A(-i),\Psi_i(\dgrsing{A}) \right\rangle.
\end{displaymath}
 \item If $a = 0$, there is an equivalence of triangulated categories,
\begin{displaymath}
 \Phi_i: \dbqgr{A} \to \dgrsing{A}.
\end{displaymath}
 \item If $a < 0$, there is a fully-faithful functor, $\Phi_i: \dbqgr{A} \to \dgrsing{A}$, and a semi-orthogonal decomposition,
\begin{displaymath}
 \dgrsing{A} \cong \left\langle k(-i),\ldots,k(-i+a+1),\Phi_i(\dbqgr{A}) \right\rangle.
\end{displaymath}
\end{enumerate}
\label{thm:orlov}
\end{theorem}

Recall that, in the case $A = k[x_0,\ldots,x_n]/I$, a well-known theorem of Serre states that $\dbqgr{A} \cong \dbcoh{X}$ where $X = \op{Proj}(A)$. If $I$ is generated by an $k[x_0,\ldots,x_n]$-regular sequence, $f_1,\ldots,f_c$, then the algebra, $A$, is AS-Gorenstein with Gorenstein parameter equal to $\sum_{i=1}^c \ops{deg}f_i - (n+1)$. So, if we can control the Orlov spectra of $\dgrsing{A}$, we can also control the Orlov spectra of $\dbcoh{X}$ where $X$ is a complete intersection. To translate statements about the category of singularities into statements about the derived category of coherent sheaves, we first need to understand what the grading shifts corresponds to on either side. We have the following lemma:

\begin{lemma}
 Let $\mathcal T$ be a triangulated category with $\mathcal I$ a thick subcategory. If we have an endofunctor, $F: \mathcal T \to \mathcal T$, so that, for any $I \in \mathcal I$, $F(I)$ is isomorphic to an object in $\mathcal I$, then $F$ descends to an endofunctor, $\bar{F}$, of $\mathcal T/\mathcal I$. Up to natural isomorphism, $\bar{F}$ is the unique functor making the following diagram commute:
\begin{center}
\begin{tikzpicture}[description/.style={fill=white,inner sep=2pt}]
\matrix (m) [matrix of math nodes, row sep=1em, column sep=1em, text height=1.5ex, text depth=0.25ex]
{ \mathcal T & & \mathcal T  \\
  & \circlearrowleft & \\
  \mathcal T / \mathcal I  & & \mathcal T / \mathcal I \\ };
\path[->,font=\scriptsize]
(m-1-1) edge node[auto] {$ F $} (m-1-3)
        edge node[auto] {$ p $} (m-3-1)
(m-3-1) edge node[auto] {$ \bar{F} $} (m-3-3)
(m-1-3) edge node[auto] {$ p $} (m-3-3);
\end{tikzpicture}
\end{center} 
Moreover, if $F$ is an autoequivalence, then $\bar{F}$ is also.
\end{lemma}

\begin{proof}
 This is a direct application of the universal property of the Verdier quotient, \cite{NeeBook} Theorem 2.1.8.
\end{proof}

The autoequivalence, $(1): \dbgr{A} \to \dbgr{A}$, preserves both $\dbtors{A}$ and $\op{perf }A$ and, therefore, descends uniquely to an autoequivalence of both $\dbqgr{A}$ and $\dgrsing{A}$, both of which shall be denoted by $(1)$. However, under the semi-orthogonal decompositions of Theorem~\ref{thm:orlov}, the two distinct versions of $(1)$ do not agree. Our first goal is to identify what operation on $\dbqgr{A}$ corresponds to $(1)$ on $\dgrsing{A}$. To do this, we need to delve a bit deeper into the proof of Theorem~\ref{thm:orlov}.  Let us now recollect the details of Orlov's work.

Let $\pi: \dbgr{A} \to \dbqgr{A}$ and $q: \dbgr{A} \to \dgrsing{A}$ denote the projections coming from Verdier localization. While $\pi$ admits a right adjoint, usually denoted by $\mathbf{R} \omega$, $q$ admits neither a right nor a left adjoint. To fix this, Orlov passes to a subcategory of $\op{gr }A$. Namely he considers, $\op{gr }A_{\geq i}$, the full subcategory of objects, $M$, of $\op{gr-}A$ with $M_j = 0$ for $j < i$. Note that the (stupid) truncation functor, $\sigma_{\geq i}: \op{gr-}A \to \op{gr }A_{\geq i}$, is right adjoint to the natural inclusion $\op{gr }A_{\geq i} \hookrightarrow \op{gr }A$. 

Denote the composition of the natural inclusion, $\dbgr{A_{\geq i}} \hookrightarrow \dbgr{A}$, and the projection, $\pi: \dbgr{A} \to \dbqgr{A}$, by $\pi_i: \dbgr{A_{\geq i}} \to \dbqgr{A}$. The functor, $\pi_i$, has the advantage of admitting a right adjoint, $\sigma_{\geq i} \circ \mathbf{R} \omega =: \omega_i$. For any graded module, $M$, we have an exact sequence,
\begin{displaymath}
 0 \to M_{\geq i} \to M \to M/M_{\geq i} \to 0.
\end{displaymath}
For any object, $X$, of $\dbgr{A}$, this induces an exact triangle,
\begin{displaymath}
 \sigma_{\geq i} X \to X \to C_X,
\end{displaymath}
with $C_X$ lying in $\dbtors{A}$. Thus, the images of $X$ and $\sigma_{\geq i}X$ are isomorphic in $\dbqgr{A}$. 

Denote the composition of the natural inclusion, $\dbgr{A_{\geq i}} \hookrightarrow \dbgr{A}$, and the projection, $q: \dbgr{A} \to \dgrsing{A}$, by $q_i: \dbgr{A_{\geq i}} \to \dgrsing{A}$. For any object, $X$, of $\dbgr{A}$, we can take a minimal graded free resolution $P \to X$. As $P$ is minimal and $A_0 = k$, $P_l$ must be generated by a free basis $e^i_l$ with $\min_i \op{deg}(e_l^i) \geq 1 + \min_i \op{deg}(e^i_{l-1})$. Thus, $P_l$ must be concentrated in degrees above $i$ for large enough $l$. We have an exact sequence of complexes,
\begin{displaymath}
 0 \to P_{< i} \to P \to P_{\geq i} \to 0,
\end{displaymath}
corresponding to splitting the free bases for each $P_l$ into those of degree less than $i$ and those of degree at least $i$. Since $P_{< i} \in \op{perf }A$ it follows that, in $\dgrsing{A}$, we have isomorphisms, $P_{\geq i} \cong P \cong X$.

From this, Orlov deduces that the left orthogonal to $\dbgr{A_{\geq i}}$, in $\dbgr{A}$, is the full subcategory of torsion complexes concentrated in degrees less than $i$, denoted by $\mathcal S_{< i}$, while the right orthogonal to $\dbgr{A_{\geq i}}$ consists of bounded complexes of free modules concentrated in degrees less than $i$, denoted by $\mathcal P_{< i}$. As $\omega_i$ is right adjoint to $\pi_i$, we see that the right orthogonal to the image of $\omega_i$ is the full subcategory of torsion complexes concentrated in degrees at least $i$, $\mathcal S_{\geq i}$. The image of $\omega_i$, denoted by $\mathcal D_i$, is equivalent to $\dbqgr{A}$.  The functor, $\omega_i$, is a quasi-inverse to the functor, $\pi_i|_{\mathcal D_i}$.

Now, the kernel of $q_i$ consists of bounded complexes of graded free modules concentrated in degree at least $i$, denoted by $\mathcal P_{\geq i}$. We now also have a nontrivial right orthogonal to the kernel of $q_i$, denote it by $\mathcal T_i$. The restriction of $q_i$ to $\mathcal T_i$ is an equivalence with $\dgrsing{A}$. The quasi-inverse is the left adjoint to $q_i$. 

From here, Orlov analyzes how the left and right orthogonals to $\mathcal D_i$ and $\mathcal T_i$ compare for different values of $a$ to prove Theorem~\ref{thm:orlov}. He finds that for $a \geq 0$, $T_i \subset D_i$ and, for $a \leq 0$, $D_i \subset T_i$. In the case $a \geq 0$, the left orthogonal to $T_i$ in $D_i$ is generated by objects isomorphic to $A(-i-a+1),\ldots,A(-i)$ in $\dbqgr{A}$. In the case $a \leq 0$, the left orthogonal to $D_i$ in $T_i$ is generated by objects isomorphic $k(-i),\ldots,k(-i+a+1)$ in $\dgrsing{A}$. 

We now begin listing a few observations about the constructions above. The autoequivalence $(1):\dgrsing{A} \to \dgrsing{A}$ admits a nice description as an autoequivalence of $D_i$.

\begin{lemma}
 $L_{A(-i+1)}$ descends to the identity functor on $\dgrsing{A}$.
\label{lem:L_AD_sg}
\end{lemma}

\begin{proof}
 It is clear that $L_{A(-i+1)}$ preserves $\op{perf }A$ and descends to a functor on $\dgrsing{A}$. The cone of the natural transformation, $\eta: \op{Id}_{\dbgr{A}} \to L_{A(-i+1)}$, lies in $\op{perf }A$. Thus, $\bar{\eta}: \op{Id}_{\dgrsing{A}} \to \bar{L}_{A(-i+1)}$ is an isomorphism.
\end{proof}

\begin{lemma}
 $\pi_i \circ L_{A(-i+1)} \circ (1) \circ \omega_i$ is isomorphic to $L_{\pi A(-i+1)} \circ (1)$ on $\dbqgr{A}$.
\label{lem:L_AD-qgr}
\end{lemma}

\begin{proof}
 As $\omega_i$ is right adjoint to $\pi_i$, if we apply $\pi_i$ to the morphism,
\begin{displaymath}
 \bigoplus_{j \in \Z} \op{Hom}_{\dbgr{A}}(A(-i+1),\omega_iF(1)[j]) \otimes_k A(-i+1)[j] \overset{\op{ev}_{\omega_iF(1)}}{\to} \omega_iF(1),
\end{displaymath}
 we get the morphism,
\begin{displaymath}
 \bigoplus_{j \in \Z} \op{Hom}_{\dbqgr{A}}(\pi A(-i+1),F(1)[j]) \otimes_k \pi A(-i+1)[j] \overset{\op{ev}_{F(1)}}{\to} \F(1),
\end{displaymath}
 Thus, $L_{\pi A(-i+1)}(F(1))$ is isomorphic to $\pi_i \circ L_{A(-i+1)} \circ (1) \circ \omega_i(F)$ for each $F$. We can take the dg-enhancements of $\dbgr{A}$ and $\dbqgr{A}$ given by bounded complexes of injectives. On the level of the dg-enhancements, the adjunctions $\pi \vdash \omega$ and $\pi_i \vdash \sigma_{\geq i} \omega$ on the abelian categories give a natural quasi-isomorphism of $L_{\pi A(-i+1)}(F(1))$ and $\pi_i \circ L_{A(-i+1)} \circ (1) \circ \omega_i(F)$. 
\end{proof}

\begin{remark}
 This lemma was first noted independently in \cite{KMV} as Lemma 5.2.1.
\end{remark}

Consider the functors,
\begin{displaymath}
 \tw{1}_i := (-i+1) \circ L_{\pi A} \circ (1) \circ (i-1): \dbqgr{A} \to \dbqgr{A}.
\end{displaymath}
Let $\tw{1} := \tw{1}_1$.

\begin{lemma}
For any $X \in D_i$, one has $L_{A(-i+1)}(X(1)) \in D_i$.
\label{lem:(1)}
\end{lemma}

\begin{proof}
 There is a triangle in $\dbgr{A}$,
\begin{displaymath}
 \op{RHom}_A(A(-i+1),X(1)) \otimes_k A(-i+1) \overset{\op{ev}}{\to} X(1) \to L_{A(-i+1)}(X(1)).
\end{displaymath}
 We know that $X$ lies in the intersection of $\leftexp{\perp}{\mathcal{P}}_{\geq i}$ and $\mathcal{P}_{< i}^{\perp}$ so $X(1)$ lies in the intersection of $\leftexp{\perp}{\mathcal{P}}_{\geq i-1}$ and $\mathcal{P}_{< i-1}^{\perp}$. As $\leftexp{\perp}{\mathcal{P}}_{\geq i} \subset \leftexp{\perp}{\mathcal{P}}_{\geq i-1}$ and $A(-i+1) \in \leftexp{\perp}{\mathcal{P}}_{\geq i}$, we see that $L_{A(-i+1)}(X(1))$ lies in $\leftexp{\perp}{\mathcal{P}}_{\geq i}$. As $A(-i+1)$ is an exceptional object in $\dbgr{A}$, $L_{A(-i+1)}(X(1))$ lies in $\langle A(-i+1)^{\perp},\mathcal{P}_{<i-1}^{\perp} \rangle = \mathcal{P}_{<i}^{\perp}$. 
\end{proof}

We saw above that $\omega_i: \dbqgr{A} \to \dbgr{A_{\geq i}}$ is full and faithful onto $D_i$. Let us denote a quasi-inverse to $q_i: T_i \to \dgrsing{A}$ by $\nu_i: \dgrsing{A} \to T_i$. 

\begin{proposition}\label{prop:whatistwist}
 If $a \geq 0$, then $q_i \circ \omega_i \circ \tw{1}_i \circ \pi_i \circ \nu_i$ is isomorphic to
\begin{displaymath}
 (1):\dgrsing{A} \to \dgrsing{A}.
\end{displaymath}
 If $ a \leq 0$, then $\pi_i \circ \nu_i \circ (1) \circ q_i \circ \omega_i$ is isomorphic to
\begin{displaymath}
 \tw{1}_i: \dbqgr{A} \to \dbqgr{A}.
\end{displaymath}
\end{proposition}

\begin{proof}
 It is easy to see that $\tw{1}_i$ is isomorphic to $L_{\pi A(-i+1)} \circ (1)$. Let us commence the manipulation proper. Assume that $a \geq 0$.
\begin{gather*}
 q_i \circ \omega_i \circ \tw{1}_i \circ \pi_i \circ \nu_i \cong q_i \circ \omega_i \circ \pi_i \circ L_{A(-i+1)} \circ (1) \circ \nu_i
\end{gather*}
 by Lemma~\ref{lem:L_AD-qgr}. By Lemma~\ref{lem:(1)}, the image of $L_{A(-i+1)} \circ (1) \circ \nu_i$ lies in $T_i \subset D_i$. So
\begin{displaymath}
 q_i \circ \omega_i \circ \pi_i \circ L_{A(-i+1)} \circ (1) \circ \nu_i \cong q_i \circ L_{A(-i+1)} \circ (1) \circ \nu_i,
\end{displaymath}
 as $\omega_i \circ \pi_i$ is isomorphic to the identity on $D_i$.
\begin{displaymath}
 q_i \circ L_{A(-i+1)} \circ (1) \circ \nu_i \cong (1)
\end{displaymath}
 by Lemma~\ref{lem:L_AD_sg}.

 Assume that $a \leq 0$.
\begin{displaymath}
 \pi_i \circ \nu_i \circ (1) \circ q_i \circ \omega_i \cong \pi_i \circ \nu_i \circ q_i \circ L_{A(-i+1)} \circ (1) \circ \omega_i
\end{displaymath}
 by Lemma~\ref{lem:L_AD_sg}. As the image of $L_{A(-i+1)} \circ (1) \circ \omega_i$ lies in $T_i$, by Lemma~\ref{lem:(1)} and $\nu_i \circ q_i$ is isomorphic to the identity on $T_i$ , we have
\begin{displaymath}
 \pi_i \circ \nu_i \circ q_i \circ L_{A(-i+1)} \circ (1) \circ \omega_i \cong \pi_i \circ L_{A(-i+1)} \circ (1) \circ \omega_i \cong \tw{1}_i
\end{displaymath}
 where the last isomorphism comes from Lemma~\ref{lem:L_AD-qgr}.
\end{proof}

\begin{remark}
Note that $\pi_i \circ \nu_i$ is $\Phi_i$ and $q_i \circ \omega_i$ is $\Psi_i$ from Theorem~\ref{thm:orlov}. Thus, Proposition~\ref{prop:whatistwist} roughly states that $(1)$ on $\dgrsing{A}$ and $\tw{1}_i$ on $\dbqgr{A}$ correspond under the semi-orthogonal decompositions of Theorem~\ref{thm:orlov}.
\end{remark}

The previous lemma becomes even more useful in the hypersurface case. To see why, we must recall the notion of a graded matrix factorization, see \cite{Orl09}. The definition is a repetition of that of the category of matrix factorization while taking care of the grading. Let $A = k[x_0,\ldots,x_n]/(f)$ with $f$ homogeneous of degree $d$. A graded matrix factorization is pair of graded free $A$-modules is diagram,
\begin{displaymath}
 P_0 \overset{p_0}{\to} P_1 \overset{p_1}{\to} P_0(d),
\end{displaymath}
of morphisms in $\op{gr }A$ so that $p_0 p_1 = f$ and $p_1 p_0 = f$. We such just denote the collection as $P$. Morphisms from $P$ to $Q$ are pairs of maps, $f_0: P_0 \to Q_1$ and $f_1: P_1 \to Q_1$, so that squares in the diagram
\begin{center}
\begin{tikzpicture}[description/.style={fill=white,inner sep=2pt}]
\matrix (m) [matrix of math nodes, row sep=3em, column sep=2.5em, text height=1.5ex, text depth=0.25ex]
{ P_0 & P_1 & P_0(d)  \\
  Q_0 & Q_1 & Q_0(d) \\ };
\path[->,font=\scriptsize]
(m-1-1) edge node[auto] {$ p_0 $} (m-1-2)
        edge node[auto] {$ f_0 $} (m-2-1)
(m-1-2) edge node[auto] {$ p_1 $} (m-1-3)
        edge node[auto] {$ f_1 $} (m-2-2)
(m-1-3) edge node[auto] {$ f_0(d) $} (m-2-3)
(m-2-1) edge node[auto] {$ q_0 $} (m-2-2)
(m-2-2) edge node[auto] {$ q_1 $} (m-2-3);
\end{tikzpicture}
\end{center}
commute. A homotopy between $f:P \to Q$ and $g: P \to Q$ is a pair of maps $h_0: P_0 \to P_1(-d)$ and $h_1:P_1 \to Q_0$ so that $f_0-g_0 = q_1h_0 + h_1p_0$ and $f_1 - g_1 = q_0h_1 + h_0 p_1$. We also have a shift $[1]$ which takes $P$ to matrix factorization
\begin{displaymath}
 P_1 \overset{p_1}{\to} P_0(d) \overset{p_0(d)}{\to} P_1(d).
\end{displaymath}

Let $\op{GrMF}(f)$ denote the homotopy category of the category of graded matrix factorizations. In \cite{Orl09}, Orlov proves the following:

\begin{theorem}
There is an equivalence of triangulated categories between $\op{GrMF}(f)$ and $\dgrsing{A}$.
\end{theorem}

We record the following elementary observations about $\op{GrMF}(f)$:

\begin{lemma}
 Let $A$ be a graded hypersurface. Then, $[2] \cong (d)$.
\end{lemma}

\begin{remark}
Combining the above lemma with Proposition \ref{prop:whatistwist} and Theorem \ref{thm:orlov}, we see that for any smooth hypersurface of degree $n+1$ in $\P^n$, one has:
\begin{displaymath}
\left(L_{\O} \circ (-\otimes_{\mathcal O} \O(1))\right)^{n+1} \cong [2].
\end{displaymath}
This isomorphism was first noticed by M. Kontsevich, \cite{Kon}, based on the relationship with the symplectic monodromy of the mirror Calabi-Yau family, see Example \ref{quintic}. The isomorphism can also be verified without reference to matrix factorizations, see \cite{Asp}.
\end{remark}

\begin{lemma}
 Let $A$ be a graded hypersurface. The natural map from $A$ to the ring of natural transformations, $\bigoplus_{i \in \Z} \op{Nat}(\op{Id}_{\dgrsing{A}}, (i))$, factors through $A/(\partial f)$.
\label{lem:gradedkernelJacobian}
\end{lemma}

\begin{proof}
 This is entirely analogous to the proof of Corollary~\ref{cor:kernelJacobian}.
\end{proof}

We can translate these into more geometric statements. Let $\langle t \rangle := \pi_1 \circ \left(L_A \circ (1)\right)^t \circ \omega_1$.

\begin{proposition}
 Let $X$ be a hypersurface of degree $d$ in $\mathbb{P}^n$ determined by a homogeneous polynomial, $f$, of degree $d$ with $A = k[x_0,\ldots,x_n]/(f)$. If $d \geq n+1$, then the natural map, $A \to \bigoplus_{i \in \Z} \op{Nat}(\op{Id}_{\dbcoh{X}},\tw{i})$, factors through $A/(\partial f)$. If $d < n+1$, then the natural map, $A \to \bigoplus_{i \in \Z} \op{Nat}(\op{Id}_{\dbcoh{X}},\tw{i})$, factors through $A/( \partial f \cdot \mathfrak m_A^a  )$ where $a=n+1-d$.
\label{prop:vanishingJac}
\end{proposition}

\begin{proof}
 Assume that $d \geq n+1$. Choose an $\alpha_i \in A_1$ for $1 \leq i \leq t$. Denote the associated natural transformation in $\dgrsing{A}$ from $\op{Id}_{\dgrsing{A}} \to (1)$ by $\eta_{\alpha_i}$ and the natural transformation from $\op{Id}_{\dbqgr{A}} \to \tw{1}$ in $\dbqgr{A}$ by $\bar{\eta}_{\alpha_i}$. Note that $q_1 \circ \omega_1$ has $\pi_1 \circ \nu_1$ as its left adjoint. To simplify notation, let us set $\Psi = q_1 \circ \omega_1$ and $\Psi^* = \pi_i \circ \nu_i$. Let $Q = \Psi \circ \Psi^*$ and denote the unit of adjunction by $e: \op{Id}_{\dgrsing{A}} \to Q$. The composition
\begin{displaymath}
 \bar{\eta}_{\alpha_t} \circ \cdots \circ \bar{\eta}_{\alpha_1} : \op{Id}_{\dbcoh{X}} \to \Psi^* \circ (1) \circ Q \circ \cdots \circ Q \circ (1) \circ \Psi = \tw{t}
\end{displaymath}
 factors
\begin{center}
\begin{tikzpicture}[description/.style={fill=white,inner sep=2pt}]
\matrix (m) [matrix of math nodes, row sep=0.25em, column sep=10em, text height=1.5ex, text depth=0.25ex]
{  &  \langle t \rangle  \\
  \op{Id}_{\dbcoh{X}} &  \\
   &  \tw{t} \\ };
\path[->,font=\scriptsize]
(m-2-1) edge node[above] {$ \Psi^* \circ \eta_{\alpha_t} \circ \cdots \circ \eta_{\alpha_1} \circ \Psi$} (m-1-2)
        edge node[below] {$ \bar{\eta}_{\alpha_t} \circ \cdots \circ \bar{\eta}_{\alpha_1} $} (m-3-2)
(m-1-2) edge (m-3-2);
\end{tikzpicture}
\end{center}
 where the map $\langle t \rangle \to \tw{t}$ comes from insertions of $e$. If $\eta_{\alpha_t} \circ \cdots \circ \eta_{\alpha_1}$ vanishes, then so does $\bar{\eta}_{\alpha_t} \circ \cdots \circ \bar{\eta}_{\alpha_1}$. The claim now follows from Lemma~\ref{lem:gradedkernelJacobian}.

%\sidenote{It appears that possibly we factor through an ideal larger than the Jacobian but one can show that $\eta_{\alpha}$ vanishes if and only if $\bar{\eta}_{\alpha}$ vanishes.}
 
 When $d < n+1$, we have the semi-orthogonal decompositions,
\begin{displaymath}
 \dbcoh{X} \cong \langle \mathcal O(-a),\ldots,\mathcal O(-1),\dgrsing{A} \rangle.
\end{displaymath}
 Note that, for $1 \leq j \leq a$,
\begin{displaymath}
 \mathcal O(-j)\tw{1} = \begin{cases} 0 & j = 1 \\ \mathcal O(-j+1)  & \text{otherwise}. \end{cases}
\end{displaymath}
 Let $F$ be an object of $\dbcoh{X}$. We can decompose $F$ using the exact triangle,
\begin{displaymath}
 F_s \to F \to F_e
\end{displaymath}
using the semiorthogonal decomposition above. Now let $\alpha$ be a polynomial of degree $i \geq a$ and $\beta$ be a polynomial of degree $j$ in $\partial f$.  We have the following commutative diagram:

%\begin{displaymath}
%\begin{CD}
%F_s @>>> F @>>> F_e\\
%@VV{\alpha}V @VV{\alpha}V @VV{\alpha}V \\
%F_s\{i\} @>>> F\{i\} @>>> F_e\{i\}\\
%@VV{\beta}V @VV{\beta}V @VV{\beta}V \\
%F_s\{i+j\} @>>> F\{i+j\} @>>> F_e\{i+j\}\\
%\end{CD}\end{displaymath}
\begin{displaymath}
\xymatrix{
F_s \ar[r] \ar[d]^{\alpha} &
F \ar[r] \ar[d]^{\alpha} \ar@{.>}[ld]&
F_e \ar[d]^{0}\\
F_s\{i\} \ar[r] \ar[d]^{0} &
F\{i\} \ar[r] \ar[d]^{\beta}&
F_e\{i\} \ar[d]^{\beta}\\
F_s\{i+j\} \ar[r]  &
F\{i+j\} \ar[r] &
F_e\{i+j\} \\
}\end{displaymath}
%\begin{center}
%\begin{tikzpicture}[description/.style={fill=white,inner sep=2pt}]
%\matrix (m) [matrix of math nodes, row sep=3em, column sep=2.5em, text height=1.5ex, text %depth=0.25ex]
%{ F_s & F & F_e  \\
%  F_s\tw{i} & F\tw{i} & F_e\tw{i} \\ };
%\path[->,font=\scriptsize]
%(m-1-1) edge node[auto] {} (m-1-2)
%        edge node[auto] {$ \alpha $} (m-2-1)
%(m-1-2) edge node[auto] {} (m-1-3)
%        edge node[auto] {$ \alpha $} (m-2-2)
%(m-1-3) edge node[auto] {$ \alpha $} (m-2-3)
%(m-2-1) edge node[auto] {} (m-2-2)
%(m-2-2) edge node[auto] {} (m-2-3);
%\end{tikzpicture}
%\end{center}
Let us justify the diagram above.  When $i \geq a$, the functor, $\{i\}$, kills all objects in $\langle \O(-a), \ldots, \O(-1) \rangle$.  Hence, $F_e\tw{i}$ is zero, which gives us the right hand zero.  This tells us that $F \overset{\alpha}{\to} F\tw{i}$ factors through $F_s\tw{i}$, represented by the dotted arrow. Now, $F_s\{i\} \overset{\beta}{\to} F_s\{i+j\}$ vanishes by Lemma~\ref{lem:gradedkernelJacobian}, which gives us the left hand zero.

Now from the diagram, we see that $F \overset{\alpha\beta}{\to} F\{i+j\}$ factors through zero and thus vanishes on the arbitrary object, $F$.  Therefore, $\alpha\beta$ lies in the kernel of the natural map $A \to \bigoplus_{i \in \Z} \op{Nat}(\op{Id}_{\dbcoh{X}},\tw{i})$.  The product ideal, $( \partial f \cdot \mathfrak m_A^a  )$, is generated by elements of this type.
\end{proof}

\begin{theorem}
 Let $f \in k[x_0, \ldots, x_n]$ be a homogeneous polynomial of degree $d$ and $A := k[x_0,\ldots,x_n]/(f)$. Assume that $A$ has an isolated singularity. For any non-zero object, $M$, in $\dgrsing{A}$, the object, $M \oplus M(1) \oplus \cdots \oplus M(d-1)$, is a generator of $\dgrsing{A}$ and
\begin{displaymath}
 \tritime(M \oplus M(1) \oplus \cdots \oplus M(d-1)) \leq 2(n+1)(d(n+1)-2n-1) - 1.
\end{displaymath}

\label{thm:somespecbound}
\end{theorem}

\begin{proof}
 After a change of basis, we can assume that $A/(x_0,\ldots,x_{n-1})$ is isomorphic to $k[u]/(u^d)$ as a graded ring. Replacing $M$ with $\op{grsyz}^n(M)$ if necessary we can assume that $M$ is a MCM module over $A$, in $\op{mod }A$. Here, $\op{grsyz}^n(M)$ is a choice of $n$th graded syzygy. Lemma~\ref{lem:slicetoAn} says that $M/(x_0,\ldots,x_{n-1})$ is nonzero in $\dsing{A_{d-1}}$. Thus, it must be nonzero in $\dgrsing{A_{d-1}}$. The proof of Theorem~\ref{thm:somespecbound} is concluded by Lemma~\ref{lem:A_n twist gen} and Lemma~\ref{lem:Mac thm}.
\end{proof}
%\begin{lemma}
% Let $\mathcal T$ be a triangulated category. If $f: X \to Y$ is a morphism in $\mathcal T$ that is a zero-divisor, then the cone of $f$ is never zero.
%\label{lem:zerodivisor}
%\end{lemma}

%\begin{proof}
% If the cone over $f$ is zero, then $f$ is an isomorphism. If $g: Y \to Z$ is a nonzero map so that $gf = 0$, then $g = 0 \circ f^{-1} = 0$. Thus, the cone over $f$ cannot be an isomorphism. 
%\end{proof}

\begin{lemma}\label{lem:A_n twist gen}
 Let $N$ be any nonzero object of $\dgrsing{A_{d-1}}$. The level of $k(0) \oplus \cdots \oplus k(d-1)$ with respect to $N \oplus N(1) \oplus \cdots \oplus N(d-1)$ is at most one. 
\end{lemma}

\begin{proof}
 Since $(d) \cong [2]$, we can assume that we have $N(i)$ for any $i \in \Z$. For some $l$ and for all $i$, we have $k[u]/(u^l)(i)$ in $\langle \lbrace N(j) \rbrace_{j \in \Z} \rangle_0$. The short exact sequences,
\begin{equation*}
 0 \to k[u]/(u^{l-1})(i) \to k[u]/(u^l)(i) \overset{u}{\to} k[u]/(u^l)(i+1) \to k(i+1) \to 0,
\end{equation*}
 split by Lemma~\ref{lem:A_nsplit}. More precisely, one can choose gradings for all the modules in the proof of Lemma~\ref{lem:A_nsplit} so that the maps are degree zero.  
\end{proof}

\begin{lemma}\label{lem:Mac thm}
 The generation time of $k(0) \oplus \cdots \oplus k(d-1)$ in $\dgrsing{A}$ is bounded above by $2(d(n+1)-2n-1)-1$.
\end{lemma}

\begin{proof}
 The proof is completely analogous to the proof of Proposition~\ref{prop:kstronggen}. Furthermore, by Macaulay's theorem, the nilpotence of $A/(\partial f)$ is $d(n+1)-2n-1$.
\end{proof}

\begin{remark}
 Theorem~\ref{thm:somespecbound} does not hold in the case of a general complete intersection, even if we allow all grading shifts, as we have already seen in Remark~\ref{rmk:CI}.
\end{remark}

We can translate this into a more geometric statement.

\begin{corollary}
 Let $X$ be a smooth hypersurface of degree $d$ in $\mathbb{P}^n$.
\renewcommand{\labelenumi}{\emph{\roman{enumi})}}
\begin{enumerate}
 \item Assume $1 < d < n+1$. Let $F \in \leftexp{\perp}{\left\langle \mathcal O(d-n-1),\ldots, \mathcal O(-1)\right\rangle}$ be nonzero. The object,
\begin{displaymath}
 \mathcal O(d-n-1) \oplus \cdots \oplus \mathcal O(-1) \oplus F \oplus \cdots \oplus F\tw{n+1},
\end{displaymath}
 is a generator of $\dbcoh{X}$ with generation time bounded by $2(n+1)(d(n+1)-2n-1)+n-d$.
 \item Assume $d = n+1$. Let $F$ be a nonzero object of $\dbcoh{X}$. The object,
\begin{displaymath}
 F \oplus \cdots \oplus F\tw{n},
\end{displaymath}
 is a generator of $\dbcoh{X}$ with generation time bounded by $2(n+1)((n+1)^2-2n-1)-1$.
 \item Assume $d > n+1$. Let $F$ be a nonzero object of $\dbcoh{X}$. The object,
\begin{displaymath}
 F \oplus F\langle 1 \rangle \oplus \cdots \oplus F\langle d-1 \rangle,
\end{displaymath}
 is a generator of $\dbcoh{X}$ with generation time bounded by $2(n+1)(d(n+1)-2n-1)-1$.
\end{enumerate}
\label{cor:geospecbound}
\end{corollary}

\begin{proof}
 Both (1) and (2) are straightforward consequences of Theorem~\ref{thm:orlov} and Theorem~\ref{thm:somespecbound} so let us assume that $d \geq n+1$ and take $F \in \dbcoh{X}$ nonzero. From Theorem~\ref{thm:somespecbound}, we know that $\omega_1(F) \oplus \omega_1(F)(1) \oplus \cdots \omega_1(F)(d-1)$ is a generator of $\dgrsing{A}$ of generation time at most $2(n+1)(d(n+1)-2n-1)$. Since $L_A \circ (1)$ is isomorphic to $(1)$ on $\dgrsing{A}$, we have
\begin{displaymath}
 \pi_1\left(\omega_1(F) \oplus \omega_1(F)(1) \oplus \cdots \oplus \omega_1(F)(d-1)\right) \cong F \oplus F\langle 1 \rangle \oplus \cdots \oplus F \langle d-1 \rangle.
\end{displaymath}
\end{proof}

\begin{remark}
Let $i: X \hookrightarrow \P^n$ denote the inclusion of a degree $d \geq n+1$ hypersurface.  One can arrange the equivalence so that, $k(n) \oplus \cdots \oplus k$ corresponds to the restriction of the exceptional collection, $\Omega^n(n),\ldots,\mathcal O$, to $X$ i.e.,
\begin{displaymath}
k(0) \oplus \cdots \oplus k(d-1) =\bigoplus_{j = 0}^n i^*\Omega^j(j).
\end{displaymath}
By getting better bounds on the generation time of this object, one can improve the bound in the corollary above.  For example, as $\bigoplus_{j = 0}^n \Omega^j(j)$ generates $\dbcoh{\P^n}$ in $n$-steps, the restriction, $\bigoplus_{j = 0}^n i^*\Omega^j(j)$, generates all objects in the essential image of $i^*$ in $n$-steps.  In particular, $\O_X(s) \in \langle \bigoplus_{j = 0}^n i^*\Omega^j(j) \rangle_n$ for all $s$.  Although we leave out the details, one can show that all objects are generated in $n-1$ steps using the category consisting of $\O_X(s)$ for all $s$.  In conclusion,
\begin{displaymath}
 \tritime(k(0) \oplus \cdots \oplus k(d-1)) \leq n(n+1) - 1.
\end{displaymath}
 And the bounds in parts $(2)$ and $(3)$ of Corollary~\ref{cor:geospecbound} improve to $n(n-1)(n+1)-1$.
\end{remark}	

\begin{remark}
 We will get a comparison point for the bound in part $(2)$ of Corollary~\ref{cor:geospecbound} in Section~\ref{sec:OSpecE} where we find the the ultimate dimension of a smooth degree three hypersurface in $\P^2$ is $4$. Our bound above is $23$. % To be more careful, any object invariant under $\tw{1}$ has generation time $1$.
\end{remark}

%Note that while there is a natural map $\langle t \rangle \to \tw{t}$ it is not an isomorphism if $t > 1$ or $d > n+1$.
The only obstacle to bounding the Orlov spectrum of $\dbcoh{X}$ is controlling the ultimate dimension under semi-orthogonal decompositions. We state the following hope:

\begin{conjecture}
 Let $\mathcal T$ be a triangulated category with a semi-orthogonal decomposition $\mathcal T = \langle \mathcal A, \mathcal B \rangle$. If the ultimate dimensions of $\mathcal A$ and $\mathcal B$ are finite, then the ultimate dimension of $\mathcal T$ is finite.
\label{conj:SO}
\end{conjecture}

\begin{proposition}
 If Conjecture~\ref{conj:SO} is true, then the ultimate dimension of $\dbcoh{X}$ is bounded for any smooth hypersurface, $X$.
\end{proposition}

\begin{proof}
%\sidenote{change to CI?}
 From Conjecture~\ref{conj:SO} and Theorem~\ref{thm:orlov}, we only have to bound the ultimate dimension of $\dgrsing{A}$ where $A=k[x_0,\ldots,x_n]/(f)$ is the homogeneous coordinate ring of $X$. Let $d$ be the degree of $f$. Since $X$ is smooth, $A$ is an isolated singularity. By the graded version of Proposition~\ref{prop:k_gen_m_nilp}, we know that the natural map, $A \to \bigoplus_{i \in \Z} \op{Nat}(\op{Id}_{\dgrsing{A}},(i))$, factors through $\mathfrak m_A^s$ for some $s$. Let $l$ be divisible by $d$ and larger than $s$. We can change coordinates so that $x_0,\ldots,x_{n-1}$ is an $A$-regular sequence. Then, $x_0^l,\ldots,x_{n-1}^l$ is also an $A$-regular sequence and $A/(x_0^l,\ldots,x_{n-1}^l)$ is graded Artinian complete intersection singularity. Let $M$ be any object of $\dgrsing{A}$. $M/(x_0^l,\ldots,x_{n-1}^l)M$ in $\dgrsing{A}$ lies in $\langle M \rangle_0$. Hence, if $M$ is a generator, then so is $M/(x_0^l,\ldots,x_{n-1}^l)M$. Using Theorem~\ref{thm:orlov} for $A/(x_0^l,\ldots,x_{n-1}^l)$, we see that the category, $\dgrsing{A/(x_0^l,\ldots,x_{n-1}^l)}$, has a full exceptional collection and thus, by Conjecture~\ref{conj:SO}, has bounded Orlov spectrum. So $\op{Lvl}_M(k(i))$ is uniformly bounded for all $i \in \Z$. Since $A$ is isolated, the category consisting of the $k(i)$ generates. This bounds the Orlov spectrum.
\end{proof}

\begin{remark}
 Note that we only need Conjecture~\ref{conj:SO} to hold for case where $\mathcal A$ is equivalent to $\dbmod{k}$. 
\end{remark}

\section{Spherical Collections} \label{sec: Spherical Collections}

In this section we explore the generation time of collections of spherical objects in triangulated categories, specifically the bounded derived category of an elliptic curve and the derived Fukaya category of a genus $g$ surface. By homological mirror symmetry for higher genus curves (see \cite{Sei08} \cite{Ef}), we can compare this to our results from Section~\ref{sec:ungraded}.  However, the method of approach is fairly different from that in Section~\ref{sec:ungraded}.  Here we use the observation of Example~\ref{twist ghost}, that spherical twists induce ghost maps, to produce ghost sequences from certain words in a braid group.  For the reader's convenience we now recall some definitions.

\begin{definition} 
Let $\mathcal T$ be the homotopy category of a triangulated $A_\infty$-category. Assume that $\mathcal T$ possesses a Serre functor, $S$. An object $\mathcal E \in \mathcal T$ is called \textbf{spherical} if,
\begin{itemize}
\item $S(\mathcal E) \cong \mathcal E[n]$
\item $\op{Hom}_{\mathcal T}(\mathcal E,\mathcal E[i]) \cong
\begin{cases} k &  i=0,n \\
0 & \text{otherwise.}
 \end{cases}$
\end{itemize}
\label{def:spherical object}
\end{definition}

\begin{definition}
Let $\mathcal T$ be the homotopy category of a triangulated $A_\infty$-category.  A collection of $m$ spherical objects, $\mathcal E_1, \ldots, \mathcal E_m$, is called an \textbf{$A_m$-configuration} if,
\begin{displaymath}
\op{dim } \left( \bigoplus_{l \in \Z} \op{Hom}_{\mathcal T}(\mathcal E_i, \mathcal E_j[l]) \right) = \begin{cases} 1 & |i-j|=1 \\ 0 & |i-j| \geq 2. \end{cases}
\end{displaymath}

\end{definition}

In definitions~\ref{def:left twist},~\ref{def:right twist}, we already discussed the notions of a left and right twist functors.  When we take the left twist functor with respect to a spherical object, we shall call this a \textbf{spherical twist}. The following result can be found in \cite{ST}:

\begin{theorem}%[Seidel-Thomas]
A spherical twist is an exact autoequivalence.  Moreover, if $\mathcal E_1, \ldots, \mathcal E_m$ is an $A_m$-configuration, then the spherical twists, $L_{\mathcal E_i}$, satisfy the braid relations:
\begin{align*}
 L_{\mathcal E_i} L_{\mathcal E_{i+1}} L_{\mathcal E_i} & \cong L_{\mathcal E_{i+1}} L_{\mathcal E_i} L_{\mathcal E_{i+1}} &i = 1, \ldots,m - 1,\\
 L_{\mathcal E_i} L_{\mathcal E_{j}} & \cong L_{\mathcal E_{j}} L_{\mathcal E_i} &|i-j| \geq 2.
\end{align*}
\end{theorem}

The following proposition will allow us to control the generation times of spherical collections:
\begin{proposition} \label{spherical upper bound}
Let $S_1,\ldots, S_n$ be spherical objects in the homotopy category, $\mathcal T$, of a triangulated cohomologically-finite $A_\infty$-category and assume we have $\emph{HH}^0(\mathcal T) = k$.  Suppose there exists a relation:
\begin{displaymath}
L_{S_{a_1}} \cdots L_{S_{a_r}} \cong \op{Id}_{\mathcal T}
\end{displaymath}
with $1 \leq a_i \leq n$.  Then $S_1 \oplus \cdots \oplus S_n$ strongly generates $\mathcal T$ with generation time at most $r-1$.  Furthermore, if we partition the relation into intervals containing mutually orthogonal spherical objects, then the generation time is at most the number of intervals minus one.
\end{proposition}

\begin{proof}
For any object, $X$, in $\mathcal T$, the left twist by $X$ comes equipped with a natural transformation, $\op{Id}_{\mathcal T} \to L_X$, which descends from a morphism of $A_\infty$-bimodules.  Composing these natural transformations yields a natural transformation,
$\zeta: \op{Id}_{\mathcal T} \to L_{S_{a_1}} \cdots L_{S_{a_r}} = \op{Id}_{\mathcal T}$.  As this descends from a morphism of $A_\infty$-bimodules, we have $\zeta \in \text{HH}^0(\mathcal T)$.  By assumption, $\zeta$ must be a scalar multiple of the identity natural transformation.  Since $\zeta$ vanishes on $S_{a_1}$ it must be zero. Hence, for any object, $X \in \mathcal T$ we get a sequence of $r$ morphisms,
\begin{displaymath}
X \to L_{S_{a_1}}(X) \to \cdots \to L_{S_{a_2}} \cdots L_{S_{a_r}}(X) \to X.
\end{displaymath}
The total map must be zero and the cones of each map lie in $\langle S_1 \oplus \cdots \oplus S_n \rangle_0$.  Repeated application of the octahedral axiom reveals that $X$ is constructed in at most $r-1$ steps.

Now if $S_1, \ldots, S_l$ are mutually orthogonal, then $L_{S_1}  \cdots  L_{S_l} = L_{S_1 \oplus \cdots \oplus S_l}$.  Thus the sequence of $l$ cones can be replaced by a single cone.  The result follows.
\end{proof}

\begin{remark}
More generally, one can assume that the relation is of the form,
\begin{displaymath}
L_{S_{a_1}}  \cdots  L_{S_{a_r}} \cong [s],
\end{displaymath}
and that any nonzero element of $\text{HH}^s(\mathcal T)$ does not vanish on $S_i$ for some $i$ (in particular one can take $\text{HH}^s(\mathcal T) =0$).
\end{remark}

\subsection{The Orlov spectrum of an elliptic curve}\label{sec:OSpecE}
In this section we study the Orlov spectrum of a smooth proper curve of genus one over an algebraically closed field.  Although this is a slight abuse of terminology, we refer to such a curve simply as an elliptic curve.  Our goal will be to prove the following theorem in a series of lemmas:
\begin{theorem}
The Orlov spectrum of the bounded derived category of coherent sheaves on an elliptic curve is $\{ 1,2,3,4 \}$.
\end{theorem}
\begin{proof}
This follows from Lemma~\ref{at most 4} and Lemma~\ref{1234} proven below.
\end{proof}

\begin{lemma}
Let $E$ be an elliptic curve and $G$ be a generator of $\dbcoh{E}$.  Then up to shifting summands, $G$ is either a vector bundle which is not semi-stable or a vector bundle plus a torsion sheaf.
\label{lem:generate curve}
\end{lemma}
\begin{proof}
Since $\text{Coh }E$ is hereditary, all complexes are isomorphic to their cohomology, (see for example \cite{H}). Hence, after shifting the summands, any generator, $G$, is a sheaf.  From Atiyah's classification of vector bundles on an elliptic curve, we may assume $G = V_1 \oplus \cdots \oplus V_n$, where the $V_i$ are indecomposable sheaves of slope $\mu_i$. If $V_i$ is torsion, then we say it has infinite slope.  If $\mu_1 = \cdots = \mu_n \neq \infty$, then by a well-known result of Faltings \cite{Fa}, there exists a vector bundle which is orthogonal to $G$. If $\mu_1 = \cdots = \mu_n = \infty$, then clearly they cannot generate $\dbcoh{E}$ as all the objects generated by this object must be torsion sheaves.  Therefore, we may assume $\mu_1 \neq \mu_2$.  As there exists an autoequivalence, $F$, of $\dbcoh{E}$ such that the slope of $F(V_2)$ has infinite slope, we may assume that $V_2$ is a torsion sheaf. Let $D$ be the support of $V_2$.  From $V_2$ and a the vector bundle $V_1$ we can get $V_1(nD)$ for all $n$.  Since the full subcategory consisting of the objects $\lbrace \O(nD) \rbrace|_{n \in \Z}$ generates $\dbcoh{E}$ and $- \otimes_{\mathcal O} V_1$ is dense (see Example~\ref{dense bundle}), it follows that $V_1 \oplus V_2$ generates.
\end{proof}

\begin{lemma} 
Let $E$ be a smooth curve of genus one. Let $V$ be a vector bundle on $E$ and $T$ be a torsion sheaf.  Then the generation time of $V \oplus T$ is bounded above by the generation time of $\O \oplus T$.
\label{lem: change out vector bundle}
\end{lemma}

\begin{proof}
The functor, $- \otimes_{\mathcal O} V$, is dense (see Example~\ref{dense bundle}). By Lemma~\ref{density lemma}, for any generator, $G$, one has:
\begin{displaymath}
\tritime(G \otimes_{\mathcal O} V) \leq \tritime(G).
\end{displaymath}
Letting $G = \O \oplus T$ we obtain:
\begin{displaymath}
\tritime((\O \oplus T) \otimes_{\mathcal O} V) = \tritime(V \oplus T^{\oplus \text{rk }(V)})= \tritime(V \oplus T) \leq \tritime(\O \oplus T).
\end{displaymath}
as desired.
\end{proof}

\begin{lemma} 
Let $E$ be an elliptic curve with identity element $e$. Let $G$ be a generator of $\dbcoh{E}$.  Then the generation time of $G$ is bounded above by the generation time of $\O \oplus \O_e$.
\label{at most 4}
\end{lemma}

\begin{proof}
Write $G = V_1 \oplus \cdots \oplus V_n$ where the $V_i$ are indecomposable sheaves of slope $\mu_i$, if $V_i$ is torsion then we say it has infinite slope.  By Lemma~\ref{lem:generate curve} at least two of these objects have different slope, by reordering we may assume $\mu_1 \neq \mu_2$.  Notice then that $V_1 \oplus V_2$ also generates and we have:
\begin{displaymath}
\tritime(G) \leq \tritime(V_1 \oplus V_2).
\end{displaymath}
As there exists an autoequivalence, $F$, of $\dbcoh{E}$ such that the slope of $F(V_2)$ is infinite, we may also assume that $V_2$ is a torsion sheaf. 

Let $\mathcal P$ be the Poincar\'e line bundle on $E \times E$.  Now we have the following inequalities:
\begin{displaymath}
\tritime(G) \leq \tritime(V_1 \oplus V_2) \leq \tritime(\O_E \oplus V_2) = \tritime(\Phi_{\mathcal P}(V_2) \oplus \O_e) \leq \tritime(\O_E \oplus \O_e).
\end{displaymath}
The first inequality is above.  The second is from Lemma~\ref{lem: change out vector bundle}.  The equality in the middle comes from applying the autoequivalence $\Phi_{\mathcal P}$ given by the Fourier-Mukai transform through the Poincar\'e line bundle.  The last inequality is achieved by applying Lemma~\ref{lem: change out vector bundle} once again.
\end{proof}

In order to calculate the generation time of objects on an elliptic curve, we appeal to Proposition 4.3 of \cite{Opp}. Let $\op{coh}_I E$ denote the subcategory of $\op{coh }X$ consisting of sheaves of slope, $\mu \in I \subset \R$. Following Oppermann, for an indecomposable vector bundle, $V$, on an elliptic curve, we define,
\begin{displaymath}
\delta(V) = \frac{q(V)}{(\text{rk }(V))^2}
\end{displaymath}
where $q(V)$ is the number of terms in the Jordan-Holder filtration of $V$.

\begin{proposition} 
 Let $V_1$ and $V_2$ be semi-stable vector bundles of slope $\mu_1$ and $\mu_2$ respectively.  Suppose $\mu_1 < \mu_2$ and $\Delta = \mu_2 - \delta(V_2) - (\mu_1 + \delta(V_1)) > 0$. Then, any coherent sheaf in
\renewcommand{\labelenumi}{\emph{\roman{enumi})}}
\begin{enumerate}
\item $\emph{coh}_{\leq \mu_1 - \delta(V_1) -\frac{\delta{(V_1)}}{\Delta}} E$,
\item $\emph{coh}_{>\mu_2 + \delta(V_2) +\frac{\delta{(V_2)}}{\Delta}} E$, or
\item $\emph{coh}_{>\mu_1 + \delta(V_1)} E \cap \emph{coh}_{\leq\mu_2 - \delta(V_2)} E$,
\end{enumerate}
 is a summand of the cone over a map from an object of $\langle V_1 \rangle_0$ to an object of $\langle V_2 \rangle_0$.
\label{prop: oppermann}
\end{proposition}

It is proven in \cite{O4}, that $\{1,2\} \subsetneq \dbcoh{C}$ for a smooth proper curve of genus at least one.  However, for completeness, let us give explicit generators of $\dbcoh{E}$ achieving the set, $\{1,2,3,4 \}$.

\begin{lemma} 
We have the following:
\begin{enumerate}
\renewcommand{\labelenumi}{\emph{\roman{enumi})}}
\item $\tritime(\O(-3e)\oplus \O \oplus \O(3e))=1$,
\item $\tritime(\O \oplus \O(3e)) = 2$,
\item $\tritime(\O \oplus \O_{2e})= 3$, and
\item $\tritime(\O \oplus \O_e)=4$.
\end{enumerate}
\label{1234}
\end{lemma}

\begin{proof}
The fact that $\tritime(\O(-3e)\oplus \O \oplus \O(3e))=1$ follows directly from Proposition~\ref{prop: oppermann} and the fact that all torsion sheaves are obtained from $\O$ and $\O(3e)$ (see also \cite{Opp} Example 4.6 and \cite{O4} Lemma 7). 

To show $\tritime(\O \oplus \O(3e)) = 2$ first note that $\O(-3e) \in \langle \O \oplus \O(3e) \rangle_1$.  As we have already shown that $\tritime(\O(-3e)\oplus \O \oplus \O(3e))=1$, we obtain, $\tritime(\O \oplus \O(3e)) \leq 2$.  For the lower bound, note that, if $p \not = q$, $\O(p-q)$ is both left and right orthogonal to $\O$. Hence, $\O(p-q)$ can not be obtained in one step as it can not be obtained from $\O(3e)$ alone (this is the argument from \cite{O4}).

To prove $\tritime(\O \oplus \O_{2e})= 3$, begin by noting that $\O(-2e), \O(2e) \in \langle \O \oplus \O_{2e} \rangle_1$ and $\O(-4e), \O(4e) \in \langle \O \oplus \O_{2e} \rangle_2$.  Applying Proposition~\ref{prop: oppermann} part $(iii)$ with $V_1 = \O(-2e)$ and $V_2 = \O(2e)$, we obtain all semi-stable bundles of slope $-1 < \mu \leq 1$  in three steps. Using $V_1 = \O$ and $V_2 = \O(4e)$, from part $(iii)$, we get all semi-stable bundles of slope $ 1 < \mu \leq 3$ and from part $(i)$ we get all semi-stable bundles with slope $\mu \leq 2$.  Now as the generator is self dual, we see that we get all possible slopes are achieved in three steps.  The torsion sheaves are obtained in one step using $\O$ and $\O(4e)$.  Hence, $\tritime(\O \oplus \O_{2e}) \leq 3$.  For the lower bound, let $q$ be a point of order two and consider the following sequence:
\begin{displaymath}
\O_q \to \O(-q)[1] \to \O(2e-q)[1] \to \O_q[1]. 
\end{displaymath}
One easily verifies that all these maps are ghost for $\O \oplus \O_{2e}$, hence by Lemma~\ref{ghost lemma} we obtain the lower bound.

Finally, to show $\tritime(\O \oplus \O_e) =4$, we use the same methods.  Note that  $\O(-e) \in \langle \O \oplus \O_e \rangle_1$ and $\O(2e) \in \langle \O \oplus \O_e \rangle_2$.  Therefore, by Proposition ~\ref{prop: oppermann} part $(iii)$ with $V_1 = \O(-e)$ and $V_2 = \O(2e)$, all semi-stable bundles of slope $\mu$ with $0 < \mu \leq 1$ are obtained in four steps.  As above, all torsion sheaves are also achieved by these two objects.  Since the generator is self dual we see that all objects of slope $-1 \leq \mu < 0$ are achieved in four steps as well.  Furthermore, as the generator is fixed under the autoequivalence given by the Poincar\'e line bundle which inverts the slopes, we see that all objects of slope $\mu =0$ or $|\mu| \geq 1$ are obtained in four steps as well.  This covers all possible slopes.  For the lower bound, let $q$ be a point of order two and consider the following sequence:
\begin{displaymath}
\O_q \to \O(-q)[1] \to \O(e-q)[1] \to \O(2e-q)[1] \to \O_q[1]. 
\end{displaymath}
One easily verifies that all these maps are ghost for $\O \oplus \O_e$ (see also Proposition~\ref{A_n gen time} below), hence by Lemma~\ref{ghost lemma} we obtain the lower bound.

\end{proof}

\subsection{The Orlov spectrum of the Fukaya category of a Riemann surface of higher genus}

In the previous section, we showed that the Orlov spectrum of the bounded derived category of coherent sheaves on an elliptic curve is $\{1,2,3,4\}$.  Via homological mirror symmetry, we could equally well view this category as the derived Fukaya category of an elliptic curve, see \cite{PZ}.  In this case, the generator with maximal generation time can be described by two loops on a torus which generate the fundamental group.

Let $\Sigma_g$ be a symplectic surface of genus g. $\Sigma_g$ admits a double branched over $\P^1$. Let $\tau: \Sigma_g \to \Sigma_g$ denote the corresponding hyperelliptic involution. Let $S_1,\ldots,S_{2g}$ be a choice of an $A_{2g}$-configuration of Lagrangian spheres, generating $H_1(\Sigma_g,\Z)$ and anti-invariant under $\tau$, up to a Hamiltonian isotopy, i.e. $\tau(S_i)$ is $S_i$ with orientation reversed up to Hamiltonian isotopy. Figure $1$ describes the situation.

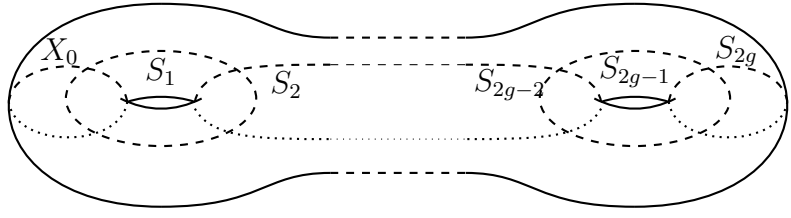
\begin{figure}[ht]
\begin{tikzpicture}[scale=0.9]
	\draw[thick] (-3.5,-1.5) .. controls (-6.5,-1.5) and (-6.5,1.5)  ..  (-3.5,1.5);
	\draw[thick] (3.5,-1.5) .. controls (6.5,-1.5) and (6.5,1.5)  ..  (3.5,1.5);
	\draw[thick] (-3.5,1.5) .. controls (-2,1.5) and (-2,1) .. (-1,1);
	\draw[thick] (-3.5,-1.5) .. controls (-2,-1.5) and (-2,-1) .. (-1,-1);
	\draw[thick] (3.5,1.5) .. controls (2,1.5) and (2,1) .. (1,1);
	\draw[thick] (3.5,-1.5) .. controls (2,-1.5) and (2,-1) .. (1,-1);
	\draw[dashed,thick] (-1,-1) to (1,-1);
	\draw[dashed,thick] (-1,1) to (1,1);
	\draw[thick] (-4.1,.1) .. controls (-3.7,-.1) and (-3.3,-.1) .. (-2.9,.1);
	\draw[thick] (-4,.05) .. controls (-3.7,.15) and (-3.3,.15) .. (-3,.05);
	\draw[thick] (4.1,.1) .. controls (3.7,-.1) and (3.3,-.1) .. (2.9,.1);
	\draw[thick] (4,.05) .. controls (3.7,.15) and (3.3,.15) .. (3,.05);
	\draw[dashed,thick] (-3,.05) .. controls (-2.9,.6) and (-2,.6) .. (-1,.6);
	\draw[dotted,thick] (-3,-.05) .. controls (-2.75,-.5) and (-2,-.5) .. (-1,-.5);
	\draw[dashed,thick] (3,.05) .. controls (2.9,.6) and (2,.6) .. (1,.6);
	\draw[dotted,thick] (3,-.05) .. controls (2.75,-.5) and (2,-.5) .. (1,-.5);
	\draw[dotted] (-1,-.5) to (1,-.5);
	\draw[dashed] (-1,.6) to (1,.6);
	\draw[dashed,thick] (-3.5,.1) ellipse (40 pt and 20 pt);
	\draw[dashed,thick] (3.5,.1) ellipse (40 pt and 20 pt);
	\draw[dashed,thick] (4,.05) .. controls (4.1,.75) and (5.65,.75) .. (5.75,.05);
	\draw[dotted,thick] (4,.05) .. controls (4.1,-.65) and (5.65,-.65) .. (5.75,.05) ;
	\draw[dashed,thick] (-4,.05) .. controls (-4.1,.75) and (-5.65,.75) .. (-5.75,.05);
	\draw[dotted,thick] (-4,.05) .. controls (-4.1,-.65) and (-5.65,-.65) .. (-5.75,.05) ;
    \node at (-3.5,.5) {$S_1$};
    \node at (-1.65,.3) {$S_2$};
    \node at (1.65,.3) {$S_{2g-2}$};
    \node at (3.5,.5) {$S_{2g-1}$};
    \node at (5,.8) {$S_{2g}$};
    \node at (-5,.8) {$X_0$};
\end{tikzpicture}
\caption{A choice of the $S_i$ and $X_0$ from the proof of Proposition~\ref{prop:dfuk}.}
\label{fig}
\end{figure}

To fix notation, we denote the morphism space from $X$ to $Y$ in $\op{DFuk}(\Sigma_g)$ by $\op{Hom}_{\Sigma_g}(X,Y)$. We can shift the gradings of the $S_i$ so that $\op{Hom}_{\Sigma_g}(S_i,S_j[1]) = 0$ for $i < j$. Let $L_i = L_{S_i}$ and let $s_i$ denote the symplectic Dehn twist about $S_i$. As mentioned in Example~\ref{twist ghost}, the endofunctor on $\op{DFuk}(\Sigma_g)$ induced by $s_i$ is $L_i$. 

In order to proceed we will need to use the following relation in the mapping class group due to Matsumoto (see \cite{M} Theorem 1.5):
\begin{theorem}
 In the mapping class group of $\Sigma_g$, we have the equality,
 \begin{displaymath}
 \left(s_1 \cdots s_{2g}\right)^{2g+1} = \tau.
 \end{displaymath}
\end{theorem}

This gives the following corollary:

\begin{corollary} \label{relation}
 We have an isomorphism of endofunctors, $(L_1 \cdots L_{2g+1})^{2g+1} \cong \tau$, of $\op{DFuk}(\Sigma_g)$.
\end{corollary}

We shall also need to know $\op{HH}^0(\op{DFuk}(\Sigma_g))$.

\begin{lemma}
 $\op{HH}^i(\op{DFuk}(\Sigma_g)) = \begin{cases}
 k   &  i=0\\
 k^{\oplus 2g} & i=1. \\
 \end{cases}$
\label{HH Curve}
\end{lemma}

\begin{proof}
 It is sufficient to compute the Hochschild homology of $\dGsing{S_g}{\Z/(2g+1)\Z}$ which can be done using the formula in Theorem 2.5.4 of \cite{PV}. The computation is straightforward. We leave the details to the reader.
\end{proof}

\begin{proposition}\label{prop:dfuk}
Let $G = S_1 \oplus \cdots \oplus S_{2g}$. Then, $4g \leq \tritime(G) \leq 8g+3$.  \label{A_n gen time}
\end{proposition}

\begin{proof}
To prove the lower bound, we construct a ghost sequence for $G$ of length $4g$.  Namely, consider a simple loop, $X_0 \in \text{DFuk}(\Sigma_g)$, which is orthogonal to $S_2, \ldots , S_{2g}$ and is anti-invariant under the hyperelliptic involution, see Figure 1.  For $0 < i  \leq 2g$ define $X_i$ inductively by $X_{i} = L_{i}(X_{i-1})$ and for $2g < i \leq 4g$ by $X_i = L_{4g+1-i}(X_{i-1})$. We also have a map, $f_i: X_i \to X_{i+1}$, given by the exact triangle:

\begin{displaymath}
 \op{Hom}_{\Sigma_g}(S_j[1],X_i) \otimes_k S_j[1] \oplus \op{Hom}_{\Sigma_g}(S_j, X_i) \otimes_k S_j \to X_i \to X_{i+1}.
\end{displaymath}
with $j= i+1$ for $0 < i  \leq 2g$ and $j=4g+1-i$ for $2g < i  \leq 4g$.

Our ghost sequence for $G$ will be the following:
\begin{displaymath}
X_0 \overset{f_0}{\to} X_1 \overset{f_1}{\to} \cdots \overset{f_{4g-1}}{\to} X_{4g}.
\end{displaymath}
In order to apply Lemma~\ref{ghost lemma} we will need to show that the total map is non-zero and $f_i$ is ghost for $G$ for all $i$.

We begin our proof by showing that for all $i$, the map $f_i$ is $G$ ghost. Equivalently, we must show that $f_i$ is ghost for $S_j$ for all $j$ and all $i$.  For notational simplicity, we will consider the case where $0 < i  \leq 2g$, though the proof is the same for $2g < i \leq 4g$.

The first step is to consider the triangle:
\begin{displaymath}
 \op{Hom}_{\Sigma_g}(S_{i+1}[1],X_i)\otimes_k S_{i+1}[1] \oplus \op{Hom}_{\Sigma_g}(S_{i+1}, X_i) \otimes_k S_{i+1} \to X_i \to X_{i+1}.
\end{displaymath}
Since any map from $S_{i+1}$ to $X_i$ factors through,
\begin{displaymath}
\op{Hom}_{\Sigma_g}(S_{i+1}[1],X_i)\otimes_k S_{i+1}[1] \oplus \op{Hom}_{\Sigma_g}(S_{i+1}, X_i) \otimes_k S_{i+1},
\end{displaymath}
it follows that $f_i$ is ghost for $S_{i+1}$.

To show that $f_i$ is ghost for $S_j$ with $j \neq {i+1}$ we will show that
\begin{equation} \label{orthogonality}
S_j \in \leftexp{\perp}{\langle} X_i \rangle \text{ unless } j=i \text{ or } j=i+1,
\end{equation}
From this equation it follows that $f_i$ is ghost for $S_j$ for $j \neq i+1$ because in this case either $\op{Hom}_{\Sigma_g}(S_j, X_i)=0$ or $\op{Hom}_{\Sigma_g}(S_j, X_{i+1})=0$.

Hence, in order to finish showing that all the $f_i$ are ghost for $G$, we must prove the orthogonality conditions of Equation~\eqref{orthogonality}.  To achieve this, we proceed by induction on $i$.  Assume Equation~\eqref{orthogonality} holds for $i-1$.  Now consider the triangle:
\begin{displaymath}
\op{Hom}_{\Sigma_g}(S_i, X_{i-1}) \otimes_k S_i \oplus \op{Hom}_{\Sigma_g}(S_i, X_{i-1}[1]) \otimes_k S_i[1]  \to X_{i-1} \to X_i.
\end{displaymath}
Let \[H_1:= \op{Hom}_{\Sigma_g}(S_i, X_{i-1}) \otimes_k \op{Hom}_{\Sigma_g}(S_j, S_i) \oplus  \op{Hom}_{\Sigma_g}(S_i, X_{i-1}[1]) \otimes_k \op{Hom}_{\Sigma_g}(S_j, S_i[1]),\]
\[H_2:=
\op{Hom}_{\Sigma_g}(S_i, X_{i-1}) \otimes_k \op{Hom}_{\Sigma_g}(S_j, S_i[1]) \oplus  \op{Hom}_{\Sigma_g}(S_i, X_{i-1}[1]) \otimes_k \op{Hom}_{\Sigma_g}(S_j, S_i).
\]
Applying the functor $\op{Hom}_{\Sigma_g}(S_j, -)$, one obtains a long exact sequence:
\begin{displaymath}
\xymatrix{
            \cdots \ar@{->}[r] & H_1\ar@{->}[r] & \op{Hom}_{\Sigma_g}(S_j,X_{i-1})\ar@{->}[r]
                   &  \op{Hom}_{\Sigma_g}(S_j,X_i) \ar `r[d] `_l[lll] `^d[dlll] `^r[dll] [dll]& \\
            & H_2 \ar@{->}[r] & \op{Hom}_{\Sigma_g}(S_j,X_{i-1}[1])\ar@{->}[r]
                   & \op{Hom}_{\Sigma_g}(S_j,X_i[1]) \ar `r[d] `_l[lll] `^d[dlll] `^r[dll] [dll]& \\
           & H_1\ar@{->}[r] & \op{Hom}_{\Sigma_g}(S_j,X_{i-1})\ar@{->}[r]
                   &  \op{Hom}_{\Sigma_g}(S_j,X_i)\ar@{->}[r]& \cdots \\}
\end{displaymath}
For $j \neq i-1, i, \text{ or } i+1$, $S_j$ is orthogonal to $S_i$ hence $H_1=H_2=0$. By the induction hypothesis, $S_j$ is orthogonal to $X_{i-1}$, hence the terms in the middle in the long exact sequence above vanish as well.  Therefore, $S_j$ is orthogonal to $X_{i}$.

The only case which remains to show is when $j=i-1$.  In this case,
\begin{displaymath}
\op{Hom}_{\Sigma_g}(S_{i-1}, S_i[1]) = 0.
\end{displaymath}
Hence, \[H_1 = \op{Hom}_{\Sigma_g}(S_i, X_{i-1}) \otimes \op{Hom}_{\Sigma_g}(S_{i-1}, S_i)\] and \[H_2 =\op{Hom}_{\Sigma_g}(S_i, X_{i-1}[1]) \otimes_k \op{Hom}_{\Sigma_g}(S_{i-1}, S_i).\]
Therefore, applying the long exact sequence, one must show that the maps
\begin{displaymath}
\alpha: \op{Hom}_{\Sigma_g}(S_i, X_{i-1}) \otimes_k \op{Hom}_{\Sigma_g}(S_{i-1}, S_i) \to \op{Hom}_{\Sigma_g}(S_{i-1},X_{i-1})
\end{displaymath}
 and
\begin{displaymath}
\beta: \op{Hom}_{\Sigma_g}(S_i, X_{i-1}[1]) \otimes_k \op{Hom}_{\Sigma_g}(S_{i-1}, S_i) \to \op{Hom}_{\Sigma_g}(S_{i-1},X_{i-1}[1])
\end{displaymath}
 are isomorphisms.  To this end, consider the following exact triangle:
\begin{displaymath}
S_{i-1} \to S_i \to L_{i-1}(S_i).
\end{displaymath}
Notice that 
\begin{displaymath}
\op{Hom}_{\Sigma_g}(L_{i-1}(S_i), X_{i-1}[k]) = \op{Hom}_{\Sigma_g}(L_{i-1}(S_i), L_{i-1}(X_{i-2})[k]) = \op{Hom}_{\Sigma_g}(S_i, X_{i-2}[k]).
\end{displaymath}
Hence, this morphism space vanishes by the induction hypothesis.  Therefore, when we apply the functor $\op{Hom}_{\Sigma_g}( - , X_{i-1})$ we get two isomorphisms:
\begin{displaymath}\op{Hom}_{\Sigma_g}( S_i , X_{i-1}) \to \op{Hom}_{\Sigma_g}( S_{i-1} , X_{i-1})
\end{displaymath} and 
\begin{displaymath}
\op{Hom}_{\Sigma_g}( S_i , X_{i-1}[1]) \to \op{Hom}_{\Sigma_g}( S_{i-1} , X_{i-1}[1]).
\end{displaymath}
 Since $\op{Hom}_{\Sigma_g}(S_{i-1},S_i)$ is one dimensional, these two isomorphisms can be identified with $\alpha$ and $\beta$.  

In summary, we have proven the validity of equation~\eqref{orthogonality} and from this we were able to deduce that all maps in this sequence are ghost for $G$.

Next, we would like to show that the total map, $X_0 \to X_{4g}$ is non-zero.  To get this result for the map from $X_0$ to $X_{4g-1}$, we proceed once again by induction.  To establish the base case, notice that as $X_0$ is not a summand of $S_1$, the map, $X_0 \to X_1$, is non-zero.    Now, consider the triangle:
\begin{displaymath}
\op{Hom}_{\Sigma_g}(S_j, X_i) \otimes_k S_j \oplus \op{Hom}_{\Sigma_g}(S_j, X_i[1]) \otimes_k S_j[1] \to X_i \to X_{i+1}.
\end{displaymath}
Applying the functor $\op{Hom}_{\Sigma_g}(X_0, - )$ to the above triangle and using the fact that $S_j$ is orthogonal to $X_0$ for $j \ge 2$, we obtain that this map in non-zero until $i=4g-1$ i.e. $X_0 \to  X_{4g-1}$ is non-zero.

Now we have:
\begin{displaymath}
X_{4g} = L_1\cdots L_{2g} L_{2g} \cdots L_1(X_0) \cong (L_1 \cdots L_{2g})^{2g+1}(X_0) \cong \tau(X_0) \cong X_0[1].
\end{displaymath}
The second equality follows from equation~\eqref{orthogonality}, the third equality comes from the relation in Corollary~\ref{relation} and the last equality comes from the fact that $X_0$ was chosen to be $\tau$-anti-invariant.  From the above equation, it follows that $X_{4g-1} = L_1^{-1}(X_0)[1]$.  This allows us to easily calculate morphisms from $S_1$ to $X_{4g-1}$.  Namely,
\begin{displaymath}
\op{Hom}_{\Sigma_g}(S_1, L_1^{-1}(X_0)) = \op{Hom}_{\Sigma_g}(S_1, X_0)
\end{displaymath}
 is one dimensional and $\op{Hom}_{\Sigma_g}(S_1, L_1^{-1}(X_0)[1]) = \op{Hom}_{\Sigma_g}(S_1, X_0[1])=0$.  Hence the map from $X_{4g-1}$ to $X_{4g}$ fits into the following triangle:
\begin{displaymath}
S_1 \to X_{4g-1}[1] \to X_{4g}[1].
\end{displaymath}
Applying the functor $\op{Hom}_{\Sigma_g}(X_0, -)$ one obtains a long exact sequence,
%\begin{gather*} 
%\op{Hom}_{\Sigma_g}(X_0, X_{4g}[1]) \to \op{Hom}_{\Sigma_g}(X_0, S_1[1]) \to %\op{Hom}_{\Sigma_g}(X_0, X_{4g-1}) \to \\ \op{Hom}_{\Sigma_g}(X_0, X_{4g})  \to %\op{Hom}_{\Sigma_g}(X_0,S_1)
%\end{gather*}
\begin{displaymath}
\xymatrix@C=2em{ 
           &  \cdots \ar[r] & \op{Hom}_{\Sigma_g}(X_0, X_{4g}[1]) \ar[r] & \op{Hom}_{\Sigma_g}(X_0, S_1[1])\ar `r[d] `_l[lll] `^d[dlll] `^r[dll] [dll] &\\
           & \op{Hom}_{\Sigma_g}(X_0, X_{4g-1}) \ar[r] & \op{Hom}_{\Sigma_g}(X_0, X_{4g}) \ar[r] & \op{Hom}_{\Sigma_g}(X_0,S_1) &}
\end{displaymath}
Since $X_{4g} = X_0[1]$, we deduce that $\op{Hom}_{\Sigma_g}(X_0, X_{4g}[1]) = \op{Hom}_{\Sigma_g}(X_0,X_0)$. Thus, the first map in the above sequence is nonzero because the identity cannot lie in the kernel. Furthermore, as $\op{Hom}_{\Sigma_g}(X_0,S_1[1])$ is one dimensional, the first map must be a surjection.  We conclude that the map, $ \op{Hom}_{\Sigma_g}(X_0, X_{4g-1}) \to  \op{Hom}_{\Sigma_g}(X_0, X_{4g}) $, is an inclusion.  As we have already deduced that our map $X_0 \to X_{4g-1}$ is nonzero, it follows that the total map, $X_0 \to X_{4g}$ is nonzero.

Ultimately, we have produced a nonzero map which factors as $4g$ ghost maps for $G$.  By Lemma~\ref{ghost lemma} we get $4g \leq \tritime(G)$.

For the upper bound one notes that,
\begin{displaymath}
(L_1 L_3 \cdots L_{2g-1} L_2 L_4 \cdots L_{2g})^{4g+2} \cong (L_1 \cdots L_{2g})^{4g+2} \cong \op{Id}_{\op{DFuk}(\Sigma_g)}.
\end{displaymath}
The first equality is just a formal relation in the braid group, the second comes from squaring the relation in Corollary~\ref{relation}.  By Lemma~\ref{HH Curve}, we can apply Proposition~\ref{spherical upper bound} which yields the upper bound, $\tritime(G) \leq 8g+3$.
\end{proof}

\begin{remark}
The beginning of the proof above works in the abstract setting.  That is, if $S_0, \ldots, S_n$ is an $A_{n+1}$-configuration of spherical objects in $\mathcal T$ such that the $A_n$-configuration, $S_1 \oplus  \cdots \oplus S_n$, generates, then $2n-1 \leq \tritime(S_1 \oplus \cdots \oplus S_n)$. 
\end{remark}

\begin{remark}
The lower bound of $4$ for the generator $\O \oplus \O_e$ on an elliptic curve is a special case of the proposition above when $\dbcoh{E}$ is viewed as a derived Fukaya category via mirror symmetry.  We notice that in this case, using various algebraic techniques, we were able to achieve an upper bound of $4$ as well.  It is believed by the authors, that on a general curve, the lower bound we prove above is in fact an equality, i.e. this generator has generation time $4g$.
\end{remark}

\vspace{2.5mm}
\noindent \textbf{Acknowledgments:}
The authors would like to warmly and heartily thank Denis Auroux, Umut Isik,  Alexander Kuznetsov, Dmitri Orlov, Tony Pantev, Paul Seidel, and Junwu Tu for valuable insight and guidance. The first author thanks the Erwin Schroedinger Institute for its hospitality. A significant portion of this paper came into focus while he was visiting the second and third authors at the ESI. The first author was funded by NSF Research Training Group Grant DMS 0636606.  The second and third authors were funded by NSF Grant DMS0600800, NSF FRG Grant DMS-0652633, by FWF Grant P20778 and by an ERC Grant.

\end{document}